\documentclass[11pt]{amsart}
\usepackage{amssymb,amsmath,amsthm}
\usepackage[all]{xy}
\xyoption{arc}
\usepackage[dvips]{graphicx}

\usepackage{a4wide}

\newcommand{\alg}{\mathrm{alg}}

\newcommand{\CH}{\operatorname{CH}}
\newcommand{\Chow}{\operatorname{CH}}

\newcommand{\Corr}{\operatorname{Corr}}
\newcommand{\cow}{\mathrm{cw}}

\newcommand{\ev}{\operatorname{ev}}

\newcommand{\Fil}{\operatorname{Fil}}

\newcommand{\Fourier}{\mathcal{F}}

\newcommand{\gen}{\mathrm{gen}}

\newcommand{\gr}{\operatorname{gr}}
\newcommand{\GrAb}{\operatorname{GrAb}}

\newcommand{\Hom}{\operatorname{Hom}}
\newcommand{\id}{\mathrm{id}}
\newcommand{\im}{\operatorname{Im}}
\renewcommand{\Im}{\operatorname{Im}}
\let\Image=\im
\newcommand{\Ind}{\operatorname{Ind}}
\newcommand{\indcatlim}{\text{``}\lim\text{''}}
\newcommand{\IndMot}{\Ind\text{-}\Mot}

\renewcommand{\ker}{\operatorname{Ker}}
\newcommand{\Ker}{\operatorname{Ker}}
\newcommand{\Kp}{\mathbb{K}}

\newcommand{\limproj}{\operatorname{proj.lim}}
\newcommand{\Lp}{\mathbb{L}}

\newcommand{\M}{\Kp\langle u\rangle}
\newcommand{\Mot}{\mathcal{M}}
\newcommand{\mot}{{\operatorname{mot}}}

\newcommand{\opp}{{\operatorname{opp}}}

\newcommand{\Pic}{\operatorname{Pic}}
\newcommand{\pr}{\mathrm{pr}}

\newcommand{\Spec}{\operatorname{Spec}}

\newcommand{\Sym}{\operatorname{Sym}}
\newcommand{\Tate}{\unitmot(1)}

\newcommand{\unitmot}{\mathbf{1}}
\newcommand{\Var}{\mathcal{V}}

\newcommand{\CC}{C}
\newcommand{\DD}{\mathcal{D}}

\newcommand{\EE}{E}

\newcommand{\JJ}{J}

\newcommand{\LL}{\mathcal{L}}

\newcommand{\OO}{\mathcal{O}}
\newcommand{\PP}{\mathcal{P}}
\newcommand{\RR}{\mathcal{R}}

\newcommand{\Proj}{\mathbb{P}}

\renewcommand{\P}{\mathbb{P}}
\newcommand{\Q}{\mathbb{Q}}

\newcommand{\Z}{\mathbb{Z}}

\newcommand{\Cbul}{\CC^{[\bullet]}}
\newcommand{\Cinf}{\CC^{[\infty]}}

\newcommand{\ov}{\overline}
\newcommand{\ol}{\overline}

\newcommand{\pa}{\partial}

\newcommand{\sub}{\subset}

\newcommand{\wt}{\widetilde}

\newcommand{\dirsum}{\mathop{\oplus}}
\newcommand{\kbar}{\bar{k}}

\newcommand{\longhookleftarrow}{\longleftarrow\joinrel\rhook}

\renewcommand{\a}{\alpha}

\newcommand{\Ga}{\Gamma}

\newcommand{\De}{\Delta}
\newcommand{\eps}{\epsilon}

\newcommand{\si}{\sigma}

\newcommand{\TCH}{\mathcal{T}\!\Chow}

\newcommand{\TKp}{\mathcal{T}\! \Kp_\Q}
\newcommand{\TL}{\mathcal{T}\! \Lp_\Q}
\newcommand{\Tbeta}{\mathcal{T}\beta}
\newcommand{\Tgamma}{\mathcal{T}\gamma}
\newcommand{\Tbetatilde}{\mathcal{T}\tilde\beta}
\newcommand{\Tgammatilde}{\mathcal{T}\tilde\gamma}

\newtheorem{introthm}{Theorem}
\numberwithin{equation}{subsection}
\newtheorem{thm}[subsection]{Theorem}
\newtheorem{prop}[subsection]{Proposition}
\newtheorem{lem}[subsection]{Lemma}
\newtheorem{cor}[subsection]{Corollary}

\theoremstyle{definition}
\newtheorem{ex}[subsection]{Example}
\newtheorem{exas}[subsection]{Examples}

\newtheorem{defi}[subsection]{Definition}

\newtheorem{rem}[subsection]{Remark}
\newtheorem{rems}[subsection]{Remarks}
\newtheorem{caution}[subsection]{Caution}

\newdimen\normalparindent       
\parindent=\normalparindent

\let\geq=\geqslant
\let\ge=\geqslant
\let\leq=\leqslant
\let\le=\leqslant
\begin{document}

\title{Algebraic cycles on the relative symmetric powers and on the
relative Jacobian of a family of curves. II}

\author{Ben Moonen}
\address{Department of Mathematics, University of Amsterdam, P.O.~Box 94248, 1090~GE Amsterdam, The Netherlands}
\email{b.j.j.moonen@uva.nl}

\author{Alexander Polishchuk}
\address{Department of Mathematics, University of Oregon, Eugene, OR 97403, USA}
\email{apolish@uoregon.edu}
\thanks{The research of Alexander Polishchuk was partially supported by NSF grant DMS-0601034}

\begin{abstract} 
Let $\CC$ be a family of curves over a non-singular variety~$S$. We study algebraic cycles on the relative symmetric powers $\CC^{[n]}$ and on the relative Jacobian~$\JJ$. We consider the Chow homology $\CH_*(\Cbul/S) := \oplus_n\, \CH_*(\CC^{[n]}/S)$ as a ring using the Pontryagin product. We prove that $\CH_*(\Cbul/S)$ is isomorphic to $\CH_*(\JJ/S)[t]\langle u\rangle$, the PD-polynomial algebra (variable:~$u$) over the usual polynomial ring (variable:~$t$) over $\CH_*(\JJ/S)$. We give two such isomorphisms that over a general base are different. Further we give precise results on how $\CH_*(\JJ/S)$ sits embedded in $\CH_*(\Cbul/S)$ and we give an explicit geometric description of how the operators $\partial_t^{[m]}$ and $\partial_u$ act. This builds upon the study of certain geometrically defined operators $P_{i,j}(a)$ that was undertaken by one of us in Part~1 of this work,~\cite{P-sym}.

Our results give rise to a new grading on $\CH_*(\JJ/S)$. The associated descending filtration is stable under all operators $[N]_*$ and $[N]_*$ acts on $\gr^m_{\Fil}$ as multiplication by~$N^m$. Hence, after $-\otimes \Q$ this filtration coincides with the one coming from Beauville's decomposition. The grading we obtain is in general different from Beauville's.

Finally we give a version of our main result for tautological classes, and we show how our methods give a simple geometric proof of some relations obtained by Herbaut and van der Geer-Kouvidakis, as later refined by one of us in~\cite{Mo}.
\end{abstract}
\maketitle

\setcounter{tocdepth}{1}
\tableofcontents


\section*{Introduction}

Let $S$ be a quasi-projective variety that is smooth of dimension~$d$ over a field. Let $\CC/S$ be a smooth curve over $S$ that has a section $p_0 \colon S \to C$. Let $\CC^{[n]}$ denote the $n$th relative symmetric power of $\CC$ over~$S$, and let $\JJ$ denote the relative Jacobian of $\CC/S$. Further let $\psi := p_0^*K \in \Pic(S)$, where $K \in \Pic(\CC)$ is the relative canonical class.
\medskip

\noindent
{\bf The central result.} In this paper we study algebraic cycles on the symmetric powers $\CC^{[n]}$ and on the Jacobian~$\JJ$, and in particular we study the relations between the corresponding Chow groups. It is of course classical that $\CC^{[n]} \cong \P(E_n)$ over~$J$, where $E_n$ is the Fourier transform of $\OO_C(n\cdot p_0)$ which for $n \geq 2g-1$ is a vector bundle on~$J$. This gives one type of relation between $\CH(\JJ)$ and the groups $\CH(\CC^{[n]})$.

What we do in this paper is something entirely different. The starting point is the remark that the disjoint union of all $\CC^{[n]}$ is a monoid in the category of (graded) schemes such that the product maps are proper. Correspondingly we have a Pontryagin product on the Chow homology $\CH_*(\Cbul/S) := \oplus_{n \geq 0}\, \CH_*(\CC^{[n]}/S)$ making it into a commutative bigraded ring. The first grading is the one by (relative) dimension of cycles, putting $\CH_i(\CC^{[n]}/S) := \CH_{d+i}(\CC^{[n]})$. The second grading is the one obtained by placing $\CH_*(\CC^{[n]}/S)$ in degree~$n$.

One of the main results of the paper is the following. (Thms.~\ref{CHJtuThm} and \ref{KtuThm} in the text.)

\begin{introthm}\label{intro-Thm1}
The ring $\CH_*(\Cbul/S)$ is isomorphic to $\CH_*(\JJ/S)[t]\langle u\rangle$.
\end{introthm}

Here we consider $\CH_*(\Cbul/S)$ and $\CH_*(\JJ/S)$ as rings with the Pontryagin product, and $\CH_*(\JJ/S)[t]\langle u\rangle$ is the PD-polynomial algebra in the variable~$u$ over the polynomial algebra in the variable~$t$ over $\CH_*(\JJ/S)$. 

Note that this result holds with integral coefficients (not just modulo torsion or after $- \otimes \Q$, in which case it would be pointless to consider PD-rings anyway) and works over any smooth and quasi-projective base. If we work over a field then by the results in~\cite{DPCRIFT}, Section~1, we have natural PD-structures on the ideal $\CH_{>0}(\Cbul/S) \subset \CH_*(\Cbul/S)$ and on the ideal in $\CH_*(\JJ/S)[t]\langle u\rangle$ generated by $\CH_{>0}(\JJ/S)$ together with all classes $u^{[m]}$ for $m \geq 1$. In Thm.~\ref{PD-thm} we prove that the isomorphism in Thm.~\ref{intro-Thm1} is compatible with these PD-structures.

Working over a general base we in fact find \emph{two} natural isomorphisms 
\begin{equation}\label{intro-tildebetagamma}
\tilde{\beta}, \tilde{\gamma} \colon \CH_*(\JJ/S)[t]\langle u\rangle \xrightarrow{\sim} \CH_*(\Cbul/S)\, .
\end{equation} 
These isomorphisms are equal modulo~$\psi$ but in general they are different. Under both isomorphisms $t$ maps to the class $[p_0] \in \CH_0(\CC/S)$ and $u^{[m]}$ maps to $\bigl[\CC^{[m]}\bigr]$; the difference---and the most nontrivial point---lies in the way that $\CH_*(\JJ/S)$ is embedded into $\CH_*(\Cbul/S)$ as a subring.

Using a slight generalization of the ``Manin principle'', the main result can also be interpreted motivically. Here we work in the ind-category $\IndMot(S)$ of the category $\Mot(S)$ of Chow motives over~$S$ (with integral coefficients) with respect to graded correspondences. We write $\unitmot$ for the identity motive and $\unitmot(1)$ for the Tate motive. The ring structures on Chow groups are encoded in the motives as multiplicative structures.

\begin{introthm}
We have isomorphisms
$$
R_*(J/S)\bigl[\unitmot\bigr]\bigl\langle \unitmot(1)\bigr\rangle \xrightarrow{\ \sim\ } R_*(\Cbul/S)
$$
in $\IndMot(S)$, compatible with the Pontryagin multiplicative structures.
\end{introthm}

See Thm.~\ref{tilde-beta-gamma-motiv} in the text.
\smallskip

Let us briefly explain how Thm.~\ref{intro-Thm1} is proven. We consider an ind-scheme~$\Cinf$, called the ``infinite symmetric power'' of~$\CC$. It is defined as the inductive limit of the system $S = \CC^{[0]} \to \CC \to \CC^{[2]} \to \cdots$, where the maps $i_n\colon \CC^{[n-1]} \to \CC^{[n]}$ are given by $D \mapsto p_0 + D$. The natural map $\tilde{\sigma} \colon \Cbul \to \JJ$ factors as $\Cbul \xrightarrow{q} \Cinf \xrightarrow{\sigma} \JJ$. Over a field, and working with $\Q$-coefficients, $\Cinf$ and its Chow homology
$$
\CH_*(\Cinf/S) := \oplus_i\, \CH_i(\Cinf/S)
\qquad\text{with}\quad
\CH_i(\Cinf/S) := \varinjlim_n\; \CH_{i+d}(\CC^{[n]})
$$
($d=\dim(S)$) have been studied by Kimura and Vistoli in \cite{KV}. We generalise and refine their results, working with integral coefficients and over more general base schemes. For a certain natural class $\Gamma \in \CH_g(\Cinf/S)$ the map $r \colon \CH_*(\JJ/S) \to \CH_*(\Cinf/S)$ given by $x \mapsto \sigma^*(x) \cap \Gamma$ is a (homomorphic) section of $\sigma_*$. Then we have the following result; see Thm.~\ref{KV-thm} and Cor.~\ref{[u]-cor} in the text.

\begin{introthm}\label{intro-beta}
We have an isomorphism $\beta \colon \CH_*(\JJ/S)\bigl\langle u\bigr\rangle \xrightarrow{\sim} \CH_*(\Cinf/S)$, given by $r$ on $\CH_*(\JJ/S)$ and with $u^{[m]} \mapsto \bigl[\CC^{[m]}\bigr]$.
\end{introthm}

As a second main ingredient for Thm.~\ref{intro-Thm1} we use that on $\Cbul/S$ we have a large collection of geometrically defined operators. Given a class $\alpha \in \CH(\CC)$ and integers $i$,~$j$, we have an operator $P_{i,j}(\alpha) \colon \CH_*(\CC^{[N]}/S) \to \CH(\CC^{[N+i-j]}/S)$ given by
$$
P_{i,j}(\alpha)\bigl(x\bigr) = (s_{i,N+i-j})_* \bigl(\pr_1^*(\alpha) \cdot s_{j,N}^*(x)\bigr)\, ,
$$
where $s_{a,b} \colon \CC \times_S \CC^{[b-a]} \to \CC^{[b]}$ is the map given by $(p,D) \mapsto a\cdot p_0+D$. Also for some of these operators there are naturally defined divided powers $P_{i,j}(\alpha)^{[m]}$. In Part~1 of this work, \cite{P-sym}, one of us has undertaken a systematic study of these operators. In particular, it was proved there that $\CH_*(\Cbul/S)$ has the structure of a module under the ring $\Z\bigl[t,u^{[\bullet]},\partial_t^{[\bullet]},\partial_u \bigr] = \Z[t,\partial_u]\langle \partial_t,u\rangle$, via 
\begin{align*}
t \mapsto P_{1,0}\bigl([p_0]\bigr)\, ,\quad &u^{[m]} \mapsto P_{1,0}(\CC)^{[m]}\, ,\\
\partial_t^{[m]} \mapsto P_{0,1}(\CC)^{[m]}\, ,\quad  &\partial_u \mapsto P_{0,1}\bigl([p_0]+\psi\bigr)\, .
\end{align*}
{}From this it readily follows (ibid., Prop.~3.9) that if we define $\Kp \subset \CH_*(\Cbul/S)$ by
$$
\Kp := \ker\bigl(P_{0,1}([p_0])\bigr) \cap\; \bigcap_{n\geq 1} \ker\bigl(P_{0,1}(\CC)^{[n]}\bigr)
$$
then we have an isomorphism 
\begin{equation}\label{Ktu=CHCinf}
\Kp\bigl[t\bigr]\bigl\langle u\bigr\rangle \xrightarrow{\sim} \CH_*(\Cbul/S)\, .
\end{equation}

The first isomorphism $\tilde{\beta}$ of \eqref{intro-tildebetagamma} is constructed as a lift of the isomorphism~$\beta$ of Thm.~\ref{intro-beta}. Namely, using \eqref{Ktu=CHCinf} we define a section $r \colon \CH_*(\Cinf/S) \to \CH_*(\Cbul/S)$ of~$q_*$; then $\tilde{\beta}$ is given by $\tilde{s} := r \circ s$ on $\CH_*(\JJ/S)$, with $t \mapsto [p_0]$ and $u^{[m]} \mapsto \bigl[\CC^{[m]}\bigr]$. Next we prove that the restriction of $\tilde{\sigma}_*$ to~$\Kp$ gives an isomorphism $\Kp \xrightarrow{\sim} \CH_*(\JJ/S)$. The second isomorphism, $\tilde{\gamma}$, is then obtained using the inverse map $\tilde{s}^\prime \colon \CH_*(\JJ/S) \to \Kp \subset \CH_*(\Cbul/S)$ on the coefficients.

Once we have the isomorphisms $\tilde{\beta}$ and $\tilde{\gamma}$, we prove (Thm.~\ref{beta-gamma-compare}) that they are equal modulo~$\psi$. (In particular, they are equal if $S = \Spec(k)$ with $k$ a field.) For both isomorphisms we have an explicit description of $\partial_t^{[m]}$ and~$\partial_u$ in terms of operators $P_{i,j}(\alpha)$. Also, the isomorphism~$\tilde{\gamma}$ descends to an isomorphism $\gamma \colon \CH_*(\JJ/S)\bigl\langle u\bigr\rangle \xrightarrow{\sim} \CH_*(\Cinf/S)$ that is different, in general, from the isomorphism~$\beta$ of Thm.~\ref{intro-beta}.
\medskip

\noindent
{\bf A new grading on $\CH_*(\JJ/S)$, and its relation with Beauville's decomposition.} As already mentioned, one of the most intriguing aspects of Thm.~\ref{intro-Thm1} is to understand how $\CH_*(\JJ/S)$ sits embedded into $\CH_*(\Cbul/S)$. In particular, we can bring into play the natural grading on $\CH_*(\Cbul/S)$ obtained by placing $\CH_*(\CC^{[n]}/S)$ in degree~$n$. We prove that~$\Kp$, which is the image of $\CH_*(\JJ/S)$ under~$\tilde{\gamma}$, is homogeneous for this grading. The result we obtain is as follows; see Thm.~\ref{filtr-Z-thm} and Cor.~\ref{filtr-Z-cor}.

\begin{introthm}
We have a decomposition
\begin{equation}\label{CHJ/S-grading}
\CH_*(\JJ/S) = \bigoplus_{m=0}^{2g+d}\, \CH_*^{[m]}(\JJ/S)\, ,
\end{equation}
where $x \in \CH_*^{[m]}(\JJ/S)$ if and only if $\tilde{\gamma}(x) \in \CH_*(\CC^{[m]}/S) \subset \CH_*(\Cbul/S)$. The associated descending filtration $\Fil^\bullet$ on $\CH_*(\JJ/S)$ is stable under the operators $[N]_*$, and $[N]_*$ acts on $\gr^m_{\Fil}$ as multiplication by~$N^m$. The subspace $\oplus_{m>2g}\, \CH_*^{[m]}(\JJ/S)$ is torsion.
\end{introthm}

The grading in \eqref{CHJ/S-grading} is compatible with the grading by (relative) dimension, so we effectively obtain a bigrading $\CH_*(\JJ/S) = \oplus\, \CH_i^{[m]}(\JJ/S)$. We show that $\CH_i^{[m]}(\JJ/S)$ can be nonzero only if $m \leq \min\{g+2d+i,2g+d\}$, and $\CH_i^{[m]}(\JJ/S)$ is torsion if $m > \min\{g+d+i,2g\}$. 

Working with $\Q$-coefficients we have Beauville's decomposition 
$$
\CH_*(\JJ/S)_\Q = \bigoplus_{i,j}\, \CH_{i,(j)}(\JJ/S)_\Q\, ,
$$ 
such that $[N]_*$ acts on $\CH_{i,(j)} \subset \CH_i(\JJ/S)_\Q$ as multiplication by $N^{2i+j}$. We say that $\CH_{i,(j)}$ has coweight $2i+j$. Our result says that $\Fil^\bullet \otimes \Q$ coincides with the descending filtration on $\CH_*(\JJ/S)_\Q$ by coweight. (This is in agreement with the general conjectures on filtrations on Chow groups; see \cite{Murre}, \cite{Jannsen}.) However, the grading given by~\eqref{CHJ/S-grading} is \emph{not}, in general, the one given by Beauville's decomposition, not even over a field. A further difference is that our grading \eqref{CHJ/S-grading} is defined integrally.
\medskip

\noindent
{\bf Tautological classes.} In the last two sections of the paper, we prove some results about the most manageable classes, the so-called tautological classes. On all three levels, $\Cbul/S$, $\Cinf/S$ and $\JJ/S$, we define a subalgebra $\TCH_*$ of tautological classes in the Chow homology. (Here we work with $\Q$-coefficients.) We prove (Thm.~\ref{Gamma-thm} and Cor.~\ref{KV-taut-thm}) that $\TCH_*(\Cbul/S)_\Q \cong \TCH_*(\JJ/S)_\Q\bigl[t,u\bigr]$ and $\TCH_*(\Cinf/S)_\Q \cong \TCH_*(\JJ/S)_\Q\bigl[u\bigr]$. As an easy application of our theory we obtain, and lift, some relations of Herbaut~\cite{Her} and van der Geer-Kouvidakis \cite{GK}, as later refined by one of us in~\cite{Mo}. Here we work over a field and we assume that the curve~$\CC$ admits a $g^r_d$. This assumption means that a certain class $[D] * L^{[r]}$ in $\CH_*(\Cinf)$ can be realised in $\CC^{[d]}$. This implies that the image of this class under the section~$r$ has no components in $\CH_*(\CC^{[m]})$ for $m>d$. Writing out what this means gives us concrete relations between classes on $\Cbul$, and pushing down to the Jacobian we recover the relations of \cite{Her}, \cite{GK}, in their refined form, working modulo rational equivalence.  

\subsection{Notation}
\label{NotConv}
Let $S$ be a non-singular quasi-projective variety of dimension~$d$ 
over a field. Let $\pi \colon \CC \to S$ be a smooth curve of 
genus~$g$ over~$S$, and let $p_0 \colon S \to \CC$ be a section 
of~$\pi$. Let
$\CC^{[n]}$ be the $n$th symmetric power of~$\CC$ relative to~$S$.
We write $\DD_n \subset \CC \times_S \CC^{[n]}$ for the universal
divisor. Further we define $\RR_n \subset \CC^{[n]}$ as the image 
of the closed embedding $i_n\colon \CC^{[n-1]} \to \CC^{[n]}$ given by
addition of the point~$p_0$, with the convention that $\RR_0$ is 
the zero divisor. Note that $\RR_n$ can also be described as the
pull-back of~$\DD_n$ under $p_0 \times \id$.

Let $\rho\colon \JJ \to S$ be the Jacobian of $\CC$ over~$S$, and
let $0\colon S \to \JJ$ be the zero section. We write $\sigma_n 
\colon \CC^{[n]} \to \JJ$ for the morphism that sends the class of 
a relative effective Cartier divisor~$D$ of degree~$n$ to the class 
of $\OO_\CC(D-n\cdot p_0)$. We usually write $\iota\colon \CC \to 
\JJ$ for the map~$\sigma_1$.

Let $\PP_{\JJ}$ be the Poincar\'e bundle on $\JJ \times_S \JJ$. Let 
$\LL_\CC$ be the Poincar\'e bundle on $\CC \times_S \JJ$, normalised 
such that its pullbacks via $p_0 \times \id$ and $\id \times e$ are 
trivial. Note that $\LL_\CC$ is the inverse of $(\iota \times \id)^* 
\PP_\JJ$.

We denote $\psi := p_0^*K\in\Pic(S)$, where $K \in 
\Pic(\CC)$ is the relative canonical class. If $h \colon X \to S$ 
is a scheme over~$S$ (such as $X = \CC^{[n]}$ or $X = \JJ$) then 
we also write $\psi$ for the class $h^*(\psi)$ in $\Pic(X)$ or 
$\CH^1(X)$. 

If $R$ is a commutative ring, and $R[t]$ is the polynomial ring in one variable over~$R$ then we write $\partial_t^{[m]}$ for the differential operator with $t^r \mapsto \binom{r}{m} t^{r-m}$. These operators are divided powers of $\partial_t = \partial_t^{[1]}$, in the sense that $\partial_t^{[m]} \partial_t^{[n]}  = \binom{m+n}{m} \partial_t^{[m+n]}$; in particular, $m!\cdot \partial_t^{[m]} = \partial_t^m$.


\section{Chow homology and Chow cohomology of infinite 
symmetric powers}

In \cite{KV}, Kimura and Vistoli consider Chow 
homology and Chow cohomology of the infinite symmetric 
power of a curve over a field. The main purpose of this 
section is to generalise and to refine their results. We
consider the relative situation of a smooth curve~$C$ over a
quasi-projective base variety~$S$, as in \ref{NotConv}, and 
we obtain results with integral coefficients.

The infinite symmetric power $\Cinf$ of~$\CC$ is defined as the direct limit of the symmetric powers $\CC^{[n]} := \Sym^n(\CC/S)$ via the inclusion maps $i_n \colon \CC^{[n-1]} \to \CC^{[n]}$ associated with the point~$p_0$. We define the Chow homology $\CH_*(\Cinf/S)$ and Chow cohomology $\CH^*(\Cinf/S)$ of~$\Cinf$. The main result of this section, see Thm.~\ref{KV-thm} and Cor.~\ref{[u]-cor}, is that $\CH_*(\Cinf/S) \cong \CH_*(\JJ/S)\bigl\langle u\bigr\rangle$, the PD-polynomial algebra in one variable over the Chow ring (with Pontryagin product) of the Jacobian. This isomorphism is made very explicit. 
\bigskip

We first define the main objects that we want to study. We consider the situation as in~\ref{NotConv}. We define
$$
\Cbul := \amalg_{n \geq 0}\; \CC^{[n]}\, ,
$$
the disjoint union of all symmetric powers of~$\CC$. Next we define
$\Cinf$, which is an ind-scheme, to be the inductive limit
$$
\lim\left(S= \CC^{[0]} \to \CC \to \CC^{[2]} \to \CC^{[3]} \to \cdots \right)
$$
where the transition maps are the morphisms $i_n\colon \CC^{[n-1]}
\to\CC^{[n]}$ associated with the point $p_0 \in\CC(S)$. Write
$\pi^{[\bullet]} \colon \Cbul \to S$ and $\pi^{[\infty]} \colon 
\Cinf \to S$ for the structural morphisms, and let $q \colon \Cbul 
\to \Cinf$ be the natural map. The collection of morphisms $\sigma_n
\colon \CC^{[n]}\to\JJ$ gives rise to morphisms $\sigma \colon \Cinf 
\to \JJ$ and $\tilde\sigma \colon \Cbul \to \JJ$ over~$S$, with $\tilde\sigma = \sigma \circ q$.

The addition maps $\alpha_{m,n} \colon \CC^{[m]} \times_S \CC^{[n]} 
\to \CC^{[m+n]}$ give rise to morphisms $\alpha \colon \Cbul \times_S 
\Cbul \to \Cbul$ and $\alpha \colon \Cinf \times_S \Cinf \to \Cinf$, 
making $\Cbul$ and~$\Cinf$ into monoid (ind-)schemes over~$S$. Further, 
the diagonal maps $\Delta_n \colon \CC^{[m]} \to \CC^{[mn]}$ given by 
$D \mapsto n\cdot D$ give rise to morphisms $[n] \colon \Cbul \to 
\Cbul$ and $[n]\colon \Cinf \to \Cinf$ over~$S$ that lift the 
``multiplication by~$n$'' maps on the Jacobian.

We define the Chow homology of $\Cbul$ over~$S$ as $\CH_*(\Cbul/S) 
:= \oplus_{i \geq 0}\, \CH_i(\Cbul/S)$, with
$$
\CH_i(\Cbul/S) := \oplus_{n \geq 0}\, \CH_{i+d}(\CC^{[n]})\, .
$$ 
On it we have a convolution product (or Pontryagin product), given by
$x * y := \alpha_*\bigl(\pr_1^*(x) \cdot \pr_2^*(y)\bigr)$, and this
gives $\CH_*(\Cbul/S)$ the structure of a graded $\CH(S)$-algebra.
Here we take $\CH(S)$ to be the usual intersection ring of~$S$ but
with grading given by putting $\CH^i(S)$ in degree~$-i$. (Formally we
should use some notation like $\CH_*(S/S)$, but as this will play no 
important role in what follows, we simply write $\CH(S)$.)  
The structural map $\CH(S) \to \CH_*(\Cbul/S)$ is given by push-forward via
the inclusion $S \to \Cbul$.

The Chow homology of $\Cinf$ over~$S$ is defined by $\CH_*(\Cinf/S) 
:= \oplus_i \CH_i(\Cinf/S)$, with 
$$
\CH_i(\Cinf/S) := \lim\left(\CH_{i+d}(\CC) \to\CH_{i+d}(\CC^{[2]}) 
\to \CH_{i+d}(\CC^{[3]}) \to\cdots \right)
$$
where the transition maps are the maps $i_{n,*}$. Again the 
convolution product defines the structure of a graded $\CH(S)$-algebra 
on $\CH_*(\Cinf/S) := \oplus_i \CH_i(\Cinf/S)$.

\begin{lem}\label{in*inj-rem}
\textup{(i)} The maps $i_{n,*} \colon \CH_*(\CC^{[n-1]}/S) \to \CH_*(\CC^{[n]}/S)$ are injective. 

\textup{(ii)} If $y$ is an element of $\CH_*(\CC^{[n]}/S) \subset \CH_*(\Cbul/S)$ then $[p_0] * y = i_{n+1,*}(y)$. 
The quotient map $q_* \colon \CH_*(\Cbul/S) \to \CH_*(\Cinf/S)$ induces an isomorphism
\begin{equation}\label{CH-modp-1}
\CH_*(\Cbul/S)/([p_0] - 1) \xrightarrow{\ \sim\ } \CH_*(\Cinf/S)\, ,
\end{equation}
where $[p_0]-1 \in \CH_*(\Cbul/S)$ is not a zero divisor.
\end{lem}

\begin{proof}
Part~(i) is proved in \cite{P-sym}, Prop.~3.10. Alternatively, this follows from the existence of correspondences $\Psi_n \colon \CC^{[n+1]} \vdash \CC^{[n]}$ with $i_n \circ \Psi_{n-1} = \id$, as in \cite{KV}, section~1. In~(ii) the fact that $[p_0]-1$ is not a zero divisor in $\CH_*(\Cbul/S)$ follows from~(i). The rest is straightforward.
\end{proof}
\medskip

The Chow cohomology of $\Cinf$ over~$S$ is defined 
by
$$
\CH^i(\Cinf/S) := \limproj \left( \CH^i(\CC) \leftarrow 
\CH^i(\CC^{[2]}) \leftarrow \CH^i(\CC^{[3]}) \leftarrow \cdots 
\right)
$$
where the transition maps are the maps~$i_n^*$. The intersection 
product induces a natural ring structure on $\CH^*(\Cinf/S) := 
\oplus_i\; \CH^i(\Cinf/S)$, and the structure of a 
$\CH^*(\Cinf/S)$-module on $\CH_*(\Cinf/S)$. We view
$\CH^*(\Cinf/S)$ as a graded algebra over $\CH^*(S)$ 
via~$\pi^{[\infty],*}$.

Note that we shall not consider Chow cohomology of~$\Cbul$.

The multiplication by~$n$ morphisms on $\Cbul$ and~$\Cinf$ (for $n 
\geq 0$) give rise to ring
endomorphisms $[n]^* \colon \CH^*(\Cinf/S) \to \CH^*(\Cinf/S)$ in 
cohomology and $[n]_* \colon \CH_*(\Cbul/S) \to \CH_*(\Cbul/S)$ 
and $[n]_* \colon \CH_*(\Cinf/S) \to \CH_*(\Cinf/S)$ in homology.

We shall also consider the Chow homology and cohomology of the
Jacobian~$\JJ$. Though they coincide as ungraded objects, even here 
it is useful to keep the distinction between homology and cohomology. 
We write $\CH^*(\JJ/S)$ for the usual Chow ring of~$\JJ$, which is a graded
$\CH^*(S)$-algebra by intersection product. The structural homomorphism
$\CH_*(S) \to \CH^*(\JJ/S)$ is the map~$\rho^*$. The Chow homology
$\CH_*(\JJ/S) :=  \CH_{*+d}(\JJ)$ is a graded $\CH(S)$-algebra by 
Pontryagin product, where the structural homomorphism $\CH(S) \to 
\CH_*(\JJ/S)$ is the map~$0_*$. With our notation we have $\CH^i(\JJ/S) 
= \CH_{g-i}(\JJ/S)$. Further, $\CH_*(\JJ/S)$ is a module over $\CH^*(\JJ/S)$
by cap-product.

For $n \in \Z$ we again have endomorphisms $[n]^*$ of $\CH^*(\JJ/S)$ and 
$[n]_*$ of $\CH_*(\JJ/S)$. The map $\sigma \colon \Cinf \to J$ gives 
rise to homomorphisms of $\CH(S)$-algebras $\sigma^*\colon 
\CH^*(\JJ/S) \to \CH^*(\Cinf/S)$ and $\sigma_* \colon \CH_*(\Cinf/S) 
\to \CH_*(\JJ/S)$. Similarly we have $\tilde\sigma_* \colon \CH_*(\Cbul/S) 
\to \CH_*(\JJ/S)$.

\begin{rem}\label{psi-times-rem}
It easily follows from the definitions that for classes $\alpha \in \CH(S)$ and $x \in \CH_*(\JJ/S)$ we have $0_*(\alpha) * x = \rho^*(\alpha) \cap x$. Similarly, for $y \in \CH_*(\Cinf/S)$ and $z \in \CH_*(\Cbul/S)$ we have $p_{0,*}(\alpha) * y = \tilde{\sigma}^*(\alpha) \cap y$ and $p_{0,*}(\alpha) * z = \sigma^*(\alpha) \cap z$.

In particular, with our convention (see \ref{NotConv}) to simply write $\psi$ for the pull-back of~$\psi$ to the scheme on which we work, we have that $\psi * x$ (meaning: $0_*(\psi) * x$, or $p_{0,*}(\psi) * x$) equals $\psi \cdot x$ (meaning: $\rho^*(\psi) \cap x$, or $\sigma^*(\psi) \cap x$, or $\tilde{\sigma}^*(\psi) \cap x$). In practise we simply write $\psi x$ or $\psi \cdot x$. 
\end{rem}

The first main goal of this section is to give a description of the 
Chow cohomology of~$\Cinf$. We begin by introducing a class~$\xi \in \CH^1(\Cinf)$. Define 
$\OO_{\CC^{[n]}}(1) := \OO_{\CC^{[n]}}(\RR_n+n\cdot \psi)$, and let 
$\xi_n := c_1\bigl(\OO_{\CC^{[n]}}(1)\bigr)$. Then 
$i_n^*\OO_{\CC^{[n]}}(1) \cong \OO_{\CC^{[n-1]}}(1)$. Hence $\xi := 
(\xi_1,\xi_2,\xi_3,\ldots)$ is a well-defined element of $\CH^1(\Cinf)$.

For a different description of the classes~$\xi_n$, consider the 
Poincar\'e bundle~$\LL_\CC$ on $\CC \times_S \JJ$. 
For $n \geq 0$ define $\EE_n := \pr_{2,*}
\bigl(\pr_1^* \OO_\CC(n\cdot p_0) \otimes \LL_\CC\bigr)$, the 
Fourier transform of~$\OO_\CC(n\cdot p_0)$. Then $\EE_n$ is a 
vector bundle on~$\JJ$ if $n > 2g-2$, and we have an isomorphism 
$\CC^{[n]} \cong \Proj(\EE_n)$ over~$\JJ$. Under this isomorphism 
$\OO_{\CC^{[n]}}(1)$ corresponds to the standard line bundle~$\OO(1)$
on~$\Proj(\EE_n)$. Note that we have the formula
$$
(\id\times\si_n)^*c_1(\LL_{\CC})= [\DD_n]- n\cdot \pr_1^*[p_0]- 
\pr_2^*[\RR_n]- n\cdot \psi
$$
for the pull-back of the Poincar\'e line bundle under the morphism
$\id\times\si_n \colon \CC\times\CC^{[n]} \to \CC\times\JJ$. This 
formula shows that 
$$
\si_n^*\EE_n \cong \pr_{2,*}\bigl(\OO_{\CC\times\CC^{[n]}}(\DD_n)
\bigr) \otimes \OO_{\CC^{[n]}}(-\RR_n - n\cdot \psi)\, ;
$$
hence, we get an embedding as a subbundle $\OO_{\CC^{[n]}}(-\RR_n - n\cdot \psi)\to\si_n^*\EE_n$, which is the standard embedding $\OO_{\CC^{[n]}}(-1) \to \si_n^*\EE_n$.

\begin{lem}\label{add-effect-lem}

\textup{(i)} We have $\alpha^* \xi = \pr_1^*\xi + \pr_2^*\xi$.

\textup{(ii)} For all $N \geq 0$ we have $[N]^* \xi = N \cdot \xi$.
\end{lem}

\begin{proof}
For (i) we use that under the addition maps $\a_{m,n}$ we have 
$\a_{m,n}^*\RR_{m+n} = \pr_1^*\RR_m + \pr_2^*\RR_n$; this follows 
for example from formula (1.1) in the proof of Lemma 1.1 
of~\cite{P-sym}. Now pass to the limit.

For (ii) we argue by induction. The case $N=1$ is trivial. For the 
induction step, note that $[N+1] \colon \Cinf \to 
\Cinf$ equals the composition
$$
\Cinf \xrightarrow{\ \Delta\ } \Cinf \times_S 
\Cinf \xrightarrow{\id \times [N]} \Cinf 
\times_S \Cinf \xrightarrow{\ \alpha\ } \Cinf\, ;
$$
now use~(i).
\end{proof}

With the lemma at our disposal the theory of \cite{KV}, \S~2 goes 
through with the same arguments. The result we obtain is the 
following. Here we view $\CH^*(\Cinf/S)$ as an algebra over 
$\CH^*(\JJ/S)$ via~$\sigma^*$.

\begin{thm}\label{Chow-cohom-thm}
The Chow cohomology $\CH^*(\Cinf/S)$ is isomorphic to the 
polynomial algebra $\CH^*(\JJ/S)\bigl[\xi\bigr]$.
\end{thm}
\medskip

Next we want to study the Chow homology of~$\Cinf$. Recall 
that $\CH_*(\Cinf/S)$ has the structure of a module over 
$\CH^*(\Cinf/S)$ by cap-product. Part~(i) of 
Lemma~\ref{add-effect-lem} implies (cf.\ \cite{KV}, Lemma~3.4) 
that the map $y\mapsto \xi \cap y$ is a derivation of 
$\CH_*(\Cinf/S)$, of degree~$-1$.

The following fact follows from the standard results on Chow groups of projective bundles.

\begin{lem}\label{proj-bun-lem}
A class $y\in\CH_*(\Cinf/S)$ is uniquely determined by 
the collection of classes $\si_*(\xi^i\cdot y)$ for $i \geq 0$. 
\end{lem}

The next proposition generalizes \cite{KV}, 
Lemma~3.6 and Proposition~3.8.

\begin{prop}\label{Gamma-L-Prop}
\textup{(i)} There is a unique element $\Gamma \in 
\CH_g(\Cinf/S)$ such that $\si_*(\Ga)=[\JJ]$ and 
$\xi\cap \Ga = 0$. The map
\begin{equation}\label{alb-section-eq}
s\colon \CH_*(\JJ/S) \to \CH_*(\Cinf/S)
\quad\text{given by}\quad 
s(y) = \si^*(y)\cap \Ga
\end{equation}
is a section of $\si_*$ and is a ring homomorphism.

\textup{(ii)} For any $m \geq 0$ there is a unique element $L^{[m]} \in 
\CH_m(\Cinf/S)$ such that
$$
\si_*(\xi^j \cap L^{[m]})=
\begin{cases}
0 & \text{if $j \neq m$}\\
[0] & \text{if $j=m$,}
\end{cases}
$$
where $[0]\in\CH_0(\JJ/S)$ is the class of the zero section $0\colon 
S\to \JJ$. Furthermore, one has $L^{[0]} = [S]$ and if $m > 0$ then $\xi\cap L^{[m]} = L^{[m-1]}$.
\end{prop}

\begin{rems}\label{Gamma-L-rems}
(i) In the last assertion, $[S] \in \CH_0(\Cinf/S)$ denotes 
the image of $[S] \in \CH(S)$ under the structural homomorphism 
$\CH(S) \to \CH_*(\Cinf/S)$. It is the unit element for the
Pontryagin product.

(ii) We shall write $L := L^{[1]} \in \Chow_1(\Cinf/S)$. The classes $L^{[m]}$ are divided powers of~$L$, i.e., $m! \cdot L^{[m]} = L^{*m}$ and $L^{[m]} * L^{[n]} = \binom{m+n}{m} L^{[m+n]}$. This readily follows from Lemma~\ref{proj-bun-lem}, the fact that $\sigma_*$ is compatible with $*$-products, plus the fact that $\xi \cap -$ is a derivation.  
\end{rems}

\begin{proof}
(i) Choose an integer $n$ with $n > 2g+d-1$. (Recall that $d = \dim(S)$.) 
Let $Q$ be the universal quotient bundle of $\sigma_n^* \EE_n$ on 
$\CC^{[n]} = \P(\EE_n)$; so $Q$ is a vector bundle of rank $n-g$ and 
we have an exact sequence
\begin{equation}\label{exactseq}
0 \longrightarrow \OO_{\CC^{[n]}}(-1) \longrightarrow \sigma_n^* \EE_n 
\longrightarrow Q \longrightarrow 0\, .
\end{equation}
Define $\Gamma \in \CH_g(\Cinf/S)$ as the class represented 
by $c_{n-g}(Q) \cap \bigl[\CC^{[n]}\bigr]$. By \cite{Fulton}, 
Example~3.3.3 we have $\sigma_*(\Gamma) = [\JJ]$. The class $\xi \cap 
\Gamma$ is represented by the element $\xi_n \cdot c_{n-g}(Q)$, which 
is the top Chern class of $\sigma_n^* \EE_n$. But $\EE_n$ is a bundle 
on~$\JJ$ of rank $n+1-g$, so given our choice of~$n$, the top Chern 
class of~$\EE_n$ vanishes for dimension reasons. In particular we find 
that for all $y \in \CH_*(\JJ/S) = \CH^*(\JJ/S)$ we have
$$
\sigma_*\bigl(\xi^k \cap (\sigma^*(y) \cap \Gamma) \bigr) = 
\begin{cases}
y & \text{if $k=0$,}\\
0 & \text{else.}
\end{cases}
$$
This shows that $s$ is a section of~$\sigma_*$.
Because $\xi \cap \Gamma =0$, the image of $s$ is contained in $\ker(\xi \cap -) \subset \CH_*(\Cinf/S)$, which is a subalgebra because $\xi \cap -$ is a derivation.
Furthermore, by Lemma~\ref{proj-bun-lem}, the restriction of $\sigma_*$ to this
subalgebra is injective. This immediately implies that $s$ is a ring homomorphism. 
Finally, the uniqueness of~$\Gamma$ follows also from Lemma~\ref{proj-bun-lem}; 
in particular, the element~$\Gamma$ that we have defined is independent 
of the choice of~$n$.

(ii) Fix $m \geq 0$ and take $n > 2g+d+m-1$. Let $Q$ again be the universal quotient bundle of $\sigma_n^* \EE_n$. Define $L^{[m]} \in \CH_1(\Cinf/S)$ 
as the class represented by $\bigl( \sigma_n^*[0] \cdot c_{n-g-m}(Q)\bigr) 
\cap \bigl[\CC^{[n]}\bigr]$. Using the projection formula we find that 
$\sigma_*(\xi^j \cap L^{[m]}) = 0$ if $j < m$, for dimension reasons. Next, $\xi^m \cap L^{[m]}$ is the class represented by $\bigl(\sigma_n^*[0] \cdot \xi^m 
\cdot c_{n-g-m}(Q)\bigr) \cap \bigl[\CC^{[n]}\bigr]$. But by our choice 
of~$n$ we have $c_{n-g-i}(\EE_n) = 0$ for all $i \in \{0,1,\ldots,m-1\}$, for dimension reasons, which gives the relation $\xi \cdot c_{n-g-i-1}(Q) = c_{n-g-i}(Q)$. Hence $\xi^m \cdot c_{n-g-m}(Q) = c_{n-g}(Q)$. So we find that $\xi^m \cap L^{[m]} = 
\sigma^*[0] \cap \Gamma$, and it follows that
$$
\sigma_*(\xi^k \cap L^{[m]}) = 
\begin{cases}
[0] & \text{if $k=m$,}\\
0 & \text{if $k > m$.}
\end{cases}
$$
The uniqueness again follows from Lemma~\ref{proj-bun-lem}.

To see that $L^{[0]} = [S]$ it suffices, again by 
Lemma~\ref{proj-bun-lem}, to show that $\sigma_*[S] = [0]$ and 
$\sigma_*\bigl(\xi^k \cap [S]\bigr) =0$ if $k > 0$. The first is clear. For the second identity, note that $\sigma_*\bigl(\xi^k \cap [S]\bigr) = 0_*\bigl(p_0^*(\xi)^k\bigr)$ and remark that $p_0^*(\xi) = 0$. Finally, the relation $\xi \cap L^{[m]} = L^{[m-1]}$ for $m>1$ readily follows from the defining property of $L^{[m]}$, once again using the lemma.
\end{proof}

\begin{rem}\label{L-over-field}
Suppose the base variety $S$ is a point. Taking $n$ sufficiently big, as in the proof of~(ii), write $F$ for the fibre of $\sigma_n \colon \CC^{[n]} \to \JJ$ over the origin $0 \in J$, and let $j \colon F \hookrightarrow \CC^{[n]}$ be the inclusion. Then $F$ is a projective space of dimension $n-g$ and we have $L^{[m]} = j_*[V_m]$, where $[V_m] = (j^*\xi)^{n-g-m} \in \CH(F)$ is the class of an $m$-plane.
\end{rem}

\begin{rems}\label{Gamma-mod-psi-rem}
(i) Modulo $\psi$ we can realise the class $\Gamma$ in $\CH_g(\CC^{[2g]}/S)$. Indeed, if we write $Q_n$ for the universal quotient bundle of $\sigma_n^*(E_n)$ then we claim that $\Gamma \equiv c_{\mathrm{top}}(Q_{2g}) = c_g(Q_{2g})$ modulo~$\psi$. The argument is the same as in (i) of the Proposition; all we need is that the top Chern class of~$E_{2g}$ vanishes modulo~$\psi$. This follows from the exact sequence
$$
0 \longrightarrow E_{2g-1} \longrightarrow E_{2g} \longrightarrow \OO_J(-2g \cdot \psi) \longrightarrow 0\, ,
$$
which we get from the sequence $0 \to \OO_C\bigl((2g-1)\cdot p_0\bigr) \to \OO_C(2g\cdot p_0) \to p_{0,*} p_0^* \OO_C(2g\cdot p_0) \to 0$, taking into account that $p_0^* \OO_C(-p_0) \cong p_0^* \Omega^1_{C/S} = \OO_S(\psi)$ by adjunction. 

(ii) With notation as in the proof, writing $L := L^{[1]}$, we have
\begin{alignat}{2}
&\Gamma = \sum_{k\geq 0}\; \xi^k \cap \biggl(\sigma_n^*\bigl(c_{n-g-k}(\EE_n)\bigr) \cap \bigl[\CC^{[n]}\bigr]\biggr)&\qquad &\text{for $n \geq 2g+d$}\, , \label{formula-for-Gamma}\\
&L = \sum_{k\geq 0}\; \xi^k \cap \biggl(\sigma_n^*\bigl([0] \cdot 
c_{n-g-1-k}(\EE_n)\bigr) \cap \bigl[\CC^{[n]}\bigr]\biggr)&\qquad &\text{for $n > 2g+d$}\, .\notag
\end{alignat}
Indeed, the exact sequence \eqref{exactseq} gives the relation $c_j(Q) = \sum_{k \geq 0}\; \xi^k \cdot c_{j-k}\bigl(\sigma_n^*(\EE_n)\bigr)$.
\end{rems}

\begin{lem}\label{sN*=N*s}
For all $N \geq 0$ we have $s \circ [N]_* = [N]_* \circ s$.
\end{lem}

\begin{proof}
For $n > 2g+d-1$, let $\Gamma(n) \in \CH_*(\CC^{[n]}/S)$ be the class representing~$\Gamma$, as constructed in the proof of Prop.~\ref{Gamma-L-Prop}. Consider the diagram
$$
\begin{matrix}
\CC^{[n]} & \xrightarrow{\ f\ } & F & \xrightarrow{\ h\ } & \CC^{[Nn]} \\
&& \llap{$\scriptstyle g$}\Big\downarrow && \Big\downarrow\rlap{$\scriptstyle \sigma_{Nn}$} \\
&& J & \xrightarrow{\ [N]\ } & J\\
\end{matrix}
$$
in which the square is cartesian, and where $f$ is the morphism such that $h \circ f = \Delta_N$ and $g \circ f = \sigma_n$. It suffices to show that $f_* \Gamma(n) = h^* \Gamma(Nn)$, as this gives
\begin{align*}
[N]_* s(y) &= h_* f_*\bigl(f^* g^*(y) \cap \Gamma(n) \bigr)\\
&= h_*\bigl(g^*(y) \cap f_*\Gamma(n)\bigr)\\
&= h_*\bigl(g^*(y) \cap h^*\Gamma(Nn)\bigr) = \bigl(h_*g^*(y)\bigr) \cap \Gamma(Nn) = s\bigl([N]_* y\bigr)\, .
\end{align*}
Now $F$ is a projective bundle over~$J$ with relatively ample class $h^*\xi$; hence to prove that $f_* \Gamma(n) = h^* \Gamma(Nn)$ it suffices to show that $g_*\bigl((h^*\xi)^i \cdot h^*\Gamma(Nn)\bigr) = g_*\bigl((h^*\xi)^i \cdot f_*\Gamma(n)\bigr)$ for all~$i$. We have
$$
g_*\bigl((h^*\xi)^i \cdot h^*\Gamma(Nn)\bigr) = g_*h^*\bigl(\xi^i \cdot \Gamma(Nn)\bigr) = [N]^* \sigma_{Nn,*}\bigl(\xi^i \cdot \Gamma(Nn)\bigr) = \begin{cases} [J] & \text{if $i=0$,}\\ 0 & \text{if $i >0$,}\end{cases} 
$$
and
\begin{align*}
g_*\bigl((h^*\xi)^i \cdot f_*\Gamma(n)\bigr) = g_*f_*\bigl(f^*h^*\xi^i \cdot \Gamma(n) \bigr) &= \sigma_{n,*}\bigl(\Delta_N^* \xi^i \cdot \Gamma(n) \bigr)\\
& = N^i \sigma_{n,*}\bigl(\xi^i \cdot \Gamma(n)\bigr) = \begin{cases} [J] & \text{if $i=0$,}\\ 0 & \text{if $i >0$,}\end{cases} 
\end{align*}
as desired.
\end{proof}

We now come to the main result of this section, which is a generalisation and refinement of \cite{KV}, Theorem~3.9. Before we can state the theorem we need to introduce some notation.

If $A$ is a ring then by $A\langle x\rangle$ we denote the PD-polynomial algebra over~$A$ in the variable~$x$; so $A\langle x\rangle = \oplus_{m \geq 0}\, A \cdot x^{[m]}$ with $x^{[0]} = 1$ and $x^{[m]} \cdot x^{[n]} = \binom{m+n}{n} x^{[m+n]}$. We have a unique PD-structure on the ideal $A\langle x\rangle_+ := \oplus_{m > 0}\, A \cdot x^{[m]}$ such that $\gamma_m(x) = x^{[m]}$ for all~$m$.

We view $\CH_*(\Cinf/S)$ as an algebra over $\CH_*(\JJ/S)$ via the 
homomorphism~$s$ in~\eqref{alb-section-eq}. Further let us introduce, for $0 \leq l \leq n$, the notation $i_{l,n} \colon \CC^{[l]} \to \CC^{[n]}$ for the composition
$$
\CC^{[l]} \xrightarrow{i_{l+1}} \CC^{[l+1]} \xrightarrow{i_{l+2}} \cdots \xrightarrow{i_{n}} \CC^{[n]}\, .
$$
In particular, $i_{n,n}$ is the identity on~$\CC^{[n]}$ and $i_{n-1,n} = i_n$. Also, $i_{0,n} \colon S = \CC^{[0]} \to \CC^{[n]}$ is the composition
$$
S \xrightarrow{p_0} \CC \xrightarrow{\Delta_n} \CC^{[n]}\, .
$$ 

\begin{thm}\label{KV-thm} 
The algebra homomorphism $h \colon\CH_*(\JJ/S)\bigl\langle x\bigr\rangle \to 
\CH_*(\Cinf/S)$ that extends the section~$s$ and sends $x^{[m]}$ to~$L^{[m]}$ is an isomorphism. 
Under this isomorphism the derivation $d/dx$ corresponds to the 
operator $\xi \cap -$, and $\sigma_*$ corresponds to the evaluation map $\CH_*(\JJ/S)\bigl\langle x\bigr\rangle \to \CH_*(\JJ/S)$ at $x = 0$. The inverse isomorphism sends a class $a \in \CH_*(\CC^{[n]}/S)$ to the polynomial $f_a \in \CH_*(\JJ/S)\bigl\langle x\bigr\rangle$ given by
\begin{equation}\label{fa-formula}
f_a = \sum_{m=0}^n\; \sigma_*\bigl(i_{n-m,n}^*(a)\bigr) \cdot \exp(\psi x)^{n-m} \cdot \left(\frac{\exp(\psi x) -1}{\psi}\right)^{[m]}\, .
\end{equation}
\end{thm}

Note that the class $\psi \in \CH^1(S)$ is nilpotent; hence $\exp(\psi x)$ is a PD-polynomial in~$x$. Further note that $\bigl(\exp(\psi x) -1\bigr)/\psi$ lies in the ideal $\CH_*(\JJ/S)\langle x\rangle_+$, so its divided powers are well-defined.

\begin{proof}
By Remark \ref{Gamma-L-rems}(ii) the map~$h$ is a homomorphism. It follows from Proposition~\ref{Gamma-L-Prop} that for any $f \in \CH_*(\JJ/S)\bigl\langle x\bigr\rangle$ we have $f(0) = \sigma_*\bigl(h(f)\bigr)$ and $h(df/dx) = \xi \cap h(f)$. This immediately implies that $h$ is injective. It remains to be shown that for $a \in \CH_*(\CC^{[n]}/S)$ the polynomial $f_a$ given in \eqref{fa-formula} indeed maps to~$a$ under~$h$. For $n=0$ this is clear. So by induction on~$n$ we may assume that $h(f_b) = b$ for all $b \in \CH_*(\CC^{[n-1]}/S)$.

Write $\alpha := h(f_a)$. Direct calculation shows that $df_a/dx = f_{i_n^*(a)} + n\cdot \psi \cdot f_a$. Applying~$h$ this gives the relation $\xi \cap \alpha = i_n^*(a) + n\cdot \psi \cdot \alpha$, where we use the induction assumption with $b = i_n^*(a)$. On the other hand, $\xi \cap a = i_n^*(a) + n \cdot \psi \cdot a$, as $\xi_n = \RR_n + n \cdot \psi$. So we find that $\xi^j \cap (\alpha - a) = (n\psi)^j \cdot (\alpha -a)$ for all $j \geq 0$. Because $\sigma_*(\alpha) = f_a(0) = \sigma_*(a)$ this implies that $\sigma_*\bigl(\xi^j \cap (\alpha -a)\bigr) = 0$ for all $j \geq 0$; hence $\alpha = a$ by Lemma~\ref{proj-bun-lem}. 
\end{proof}

\begin{ex}\label{KV-thm-exa}
For $a \in \CH_*(\CC/S)$ we have $f_a = \iota_*(a) \cdot \exp(\psi x) + p_0^*(a) \cdot \frac{\exp(\psi x) - 1}{\psi}$. (Recall that $\iota = \sigma_1 \colon \CC \to \JJ$.) In particular, the class $[\CC]$ corresponds to the polynomial 
\begin{align}\label{[C]-formula}
f_{[C]} &= \bigl[\iota(C)\bigr]\cdot \exp(\psi x) + \frac{\exp(\psi x)-1}{\psi} \\
&= \bigl[\iota(C)\bigr] + \bigl(1+\psi\cdot [\iota(C)]\bigr) \cdot \bigl(x+\psi\, x^{[2]} + \psi^2\, x^{[3]} + \cdots\bigr)\, .\nonumber
\end{align}
To avoid any misunderstanding let us note again that the ring multiplication in $\CH_*(\Cinf/S)$ is the Pontryagin product and that we view $\CH_*(\Cinf/S)$ as an algebra over $\CH_*(\JJ/S)$ via the homomorphism~$s$. Thus, for instance, the formula for $f_{[C]}$ just given means that
$$
[C] = s\bigl[\iota(C)\bigr] + \bigl(1+\psi\cdot s[\iota(C)]\bigr) * \bigl(L + \psi\, L^{[2]} + \psi^2\, L^{[3]} + \cdots\bigr)\, .
$$

More generally, for $n \geq 1$ the class of $\CC^{[n]}$ corresponds to the polynomial
\begin{equation}\label{C[n]Polynomial}
f_{[\CC^{[n]}]} = \sum_{m=0}^n\; \sigma_*\bigl[\CC^{[n-m]}\bigr]  \cdot \exp(\psi x)^{n-m} \cdot \left(\frac{\exp(\psi x) - 1}{\psi}\right)^{[m]}\, .
\end{equation}
\end{ex}
\medskip

As a corollary of the preceding results we obtain expression for the classes $L$ and~$\Gamma$ as linear combinations of explicit geometric classes with coefficients in the ring $\Z[\psi]$.

\begin{cor}\label{L-Gamma-form}
In $\CH_*(C^{[\infty]}/S) \otimes \Q$ we have the identities
\begin{equation}\label{L-formula}
\begin{split}
L & = \frac{\log\bigl(1+\psi\cdot [C]\bigr) - \log\bigl(1+\psi\cdot s[\iota(C)]\bigr)}{\psi} \\
& = \sum_{n\ge 1} (-\psi)^{n-1} (n-1)!\, \bigl([C^{[n]}] - s[\sigma_{n,*}C^{[n]}]\bigr)^{[n]}\, ;\\
\end{split}
\end{equation}
\begin{equation}\label{Gamma-L-eq}
\Ga = \left(\frac{\exp(-\psi\cdot L)-1}{\psi} + [C] * \exp(-\psi\cdot  L)\right)^{[g]}\, ;
\end{equation}
and for any $n \geq 2$:
\begin{equation}\label{Gamma-taut-eq}
\Ga = \frac{1}{(n-1)^g} \left(\frac{\log\bigl(1+\psi\cdot [n]_*[C]\bigr) - n\log\bigl(1+\psi\cdot [C]\bigr)}{n\psi}\right)^{[g]}\, .
\end{equation}
\end{cor}

Over a field the first two identities actually hold in $\CH_*(C^{[\infty]}/S)$; see Remark~\ref{Gamma-form-int}.

\begin{proof}
Identity \eqref{L-formula} is just the inverse of formula \eqref{[C]-formula}. For the second identity we first rewrite \eqref{[C]-formula} as
$$
s\bigl([\iota(C)]\bigr) = \frac{\exp(-\psi L)-1}{\psi} + [C] * \exp(-\psi L)\, .
$$
But also we have
\begin{equation}\label{siotaC*g}
s([\iota(C)])^{* g} = s\bigl([\iota(C)]^{*g}\bigr) = g!\cdot s\bigl([J]\bigr) = g!\cdot \Ga\, .
\end{equation}
This gives~\eqref{Gamma-L-eq}.

To deduce \eqref{Gamma-taut-eq} we start with the relation $[n]_* L = n \cdot L$. Using \eqref{L-formula} this gives
$$
\frac{\log\bigl(1+\psi\cdot [n]_*[C]\bigr)}{n\cdot \psi} - \frac{\log\bigl(1+\psi\cdot [C]\bigr)}{\psi}
=
s\left(
\frac{\log\bigl(1+\psi \cdot [n]_*[\iota(C)]\bigr)}{n\cdot \psi} - \frac{\log\bigl(1+\psi \cdot [\iota(C)]\bigr)}{\psi}
\right)\, .
$$
It remains to be shown that the $g$th power of the right-hand side
is equal to $g!(n-1)^g\cdot \Ga$. But when calculating the $g$th Pontryagin power
of a class in $\CH_1(J/S)$ we can work modulo~$\psi$ and modulo homological equivalence. Now use that
$$
\frac{\log\bigl(1+\psi \cdot [n]_*[\iota(C)]\bigr)}{n\cdot \psi} - \frac{\log\bigl(1+\psi \cdot [\iota(C)]\bigr)}{\psi} \equiv \frac{[n]_* \bigl[\iota(\CC)\bigr]}{n} - \bigl[\iota(\CC)\bigr] \qquad \bmod (\psi)
$$
and that $[n]_* \bigl[\iota(\CC)\bigr]$ is homologically equivalent to $n^2 \cdot \bigl[\iota(\CC)\bigr]$. 
\end{proof}

\begin{cor}\label{[u]-cor}
The homomorphism of $\CH(S)$-algebras
$$
\beta \colon \CH_*(\JJ/S)\bigl\langle u\bigr\rangle \to\CH_*(\Cinf/S)$$
restricting to $s$ on $\CH_*(\JJ/S)$ and sending $u^{[m]}$ to~$\bigl[\CC^{[m]}\bigr]$, is an isomorphism.
\end{cor}

\begin{proof}
Write $A := \CH_*(\JJ/S)$. We have a nilpotent element $\psi \in A$ and, via the isomorphism of Theorem~\ref{KV-thm}, a homomorphism of $A$-algebras $\beta \colon A\bigl\langle u\bigr\rangle \to A\bigl\langle x\bigr\rangle$. Write $\ol{A} := A/(\psi)$ and consider the induced map $\ol{\beta} \colon \ol{A}\bigl\langle u\bigr\rangle \to \ol{A}\bigl\langle x\bigr\rangle$. By \eqref{C[n]Polynomial} we have
$$
\ol{\beta}\bigl(u^{[m]}\bigr) \in x^{[m]} + \sum_{i=0}^{m-1} \ol{A} \cdot x^{[i]}
$$
for all $m \geq 0$. {}From this it follows by an elementary argument that $\beta$ is an isomorphism.
\end{proof}

\begin{rem}
The relations we have obtained greatly simplify if we calculate modulo~$\psi$. (This applies for instance if we work over a field.) Theorem~\ref{KV-thm} gives
$$
a \equiv \sum_{m=0}^n s\bigl(\sigma_*(i_{n-m,n}^*(a))\bigr) * L^{[m]} \quad\bmod (\psi)
$$
for $a \in \CH_*(\CC^{[n]}/S)$. Modulo $\psi$ the identities \eqref{Gamma-L-eq} and~\eqref{Gamma-taut-eq}
take the form
$$
\Ga \equiv \frac{\bigl([C]-L\bigr)^{* g}}{g!} 
\equiv \frac{1}{g!} \cdot \left(\frac{[n]_*[C] - n\cdot [C]}{n(n-1)}\right)^{* g}
\mod(\psi) \, .
$$
\end{rem}


\section{Operators on the Chow (co)homology}\label{Operators}

In this section we study several geometrically defined operators on $\CH_*(\Cbul/S)$. These operators, which have been studied in detail in \cite{P-sym}, play a key role in the proofs of the main results of the paper. After recalling the definitions, we give some examples of the operators that are most important for us, and we prove some identities that are used later. The main result is Cor.~\ref{deriv-cor}.
\bigskip 

Given integers $0 \leq m \leq n$, let $s_{m,n} \colon \CC \times_S \CC^{[n-m]} \to \CC^{[n]}$ be the morphism given by $(p,D) \mapsto D + m\cdot p$. In \cite{P-sym}, one of us has defined and studied a family of operators~$P_{i,j}(a) \colon \CH_*(\CC^{[n]}/S) \to \CH_*(\CC^{[n+i-j]}/S)$, for $i$, $j \geq 0$ and $a \in \CH(\CC)$, given by
$$
P_{i,j}(a)\bigl(x\bigr) := (s_{i,n+i-j})_*\bigl(\pr_1^*(a) \cdot s_{j,n}^*(x)\bigr)\, .
$$
These give rise to operators $P_{i,j}(a)$ on $\CH_*(\Cbul/S)$. Also in \cite{P-sym} the commutation relations between the various operators $P_{i,j}(a)$ were calculated.

\begin{exas}\label{Pij-exa} 
We write $[p_0]$ for the class of $p_0(S)$ in $\CH_0(\CC/S)$.

(a) The operator $P_{0,0}(a) \colon \CH_*(\CC^{[n]}/S) \to \CH_*(\CC^{[n]}/S)$ is given by $x \mapsto \pi_*(a) \cdot x$, where we recall that $\pi\colon \CC \to S$ is the structural morphism. (To avoid confusion, note that in the expression $\pi_*(a) \cdot x$ we view $\pi_*(a)$ as an element of $\CH(\CC^{[n]})$ via pullback, and the dot denotes intersection product.)
 
(b) The operator $P_{1,0}(a) \colon \CH_*(\CC^{[n]}/S) \to \CH_*(\CC^{[n+1]}/S)$ is given by $x \mapsto a * x$. In particular, $P_{1,0}\bigl([p_0]\bigr) = i_{n+1,*}$.
 
(c) Let $\pr_1 \colon \CC \times_S \CC^{[n-1]} \to \CC$ and $\pr_2 \colon \CC \times_S \CC^{[n-1]} \to \CC^{[n-1]}$ be the projections and let $\alpha = \alpha_{1,n-1} \colon \CC \times_S \CC^{[n-1]} \to \CC^{[n]}$. Then the operator $P_{0,1}(a) \colon \CH_*(\CC^{[n]}/S) \to \CH_*(\CC^{[n-1]}/S)$ is given by $x \mapsto \pr_{2,*}\bigl(\pr_1^*(a) \cdot \alpha^*(x)\bigr)$. In particular, $P_{0,1}\bigl([p_0]\bigr) = i_n^*$. Using this last identity it is not difficult to show that for $a \in \CH(\CC)$ and $n \geq 0$, we have the relation $P_{0,1}([p_0])\bigl(\Delta_{n,*}(a)\bigr) = n\cdot p_0^*(a) * [p_0]^{*(n-1)}$. We shall use this later.

(d) The operator $P_{1,1}(a) \colon \CH_*(\CC^{[n]}/S) \to \CH_*(\CC^{[n]}/S)$ is given by $x \mapsto \alpha_*\bigl(\pr_1^*(a)\bigr) \cdot x = \bigl(a * [\CC^{[n-1]}] \bigr) \cdot x$. For instance, for $x \in \CH_*(\CC^{[n]}/S)$ we have $P_{1,1}(C)\bigl(x\bigr) = n\cdot x$, and $P_{1,1}\bigl([p_0]+\psi\bigr)(x) = \xi_n \cap x$. Also, $P_{1,1}\bigl([p_0]\bigr) = P_{1,0}\bigl([p_0]\bigr) \circ P_{0,1}\bigl([p_0]\bigr)$, which is the operator $i_{n,*} i_n^*$ sending $x$ to $R_n \cdot x$.
\end{exas} 
\medskip

As we have seen in the above examples, $P_{1,0}\bigl([p_0]\bigr)$ is the operator $x \mapsto [p_0] * x$. Recall from Lemma~\ref{in*inj-rem}(ii) that $\CH_*(\Cinf/S)$ is the quotient of $\CH_*(\Cbul/S)$ modulo the ideal $\bigl([p_0]-1\bigr)$. Hence any operator on $\CH_*(\Cbul/S)$ that commutes with $P_{1,0}\bigl([p_0]\bigr)$ induces an operator on $\CH_*(\Cinf/S)$. In particular, it follows from the commutation relations in Thm.~0.1 of~\cite{P-sym} that this applies to all $P_{i,j}(a)$ for $a \in \CH(\CC)$ with $p_0^*(a) = 0$, so for all such classes~$a$ we get induced operators $\ol{P}_{i,j}(a)$ on $\CH_*(\Cinf/S)$. 

As an example, the derivation $\xi \cap -$ considered before is the 
operator $\ol{P}_{1,1}\bigl([p_0] + \psi\bigr)$.
\medskip

Next we recall that in \cite{P-sym}, section~3, also divided powers of the operators $P_{0,n}(\CC)$ and $P_{n,0}(\CC)$ were introduced. Concretely, we define
$$
P_{n,0}(\CC)^{[m]} \bigl( x\bigr) := \delta_n^{[m]} * x
\qquad\text{and}\qquad
P_{0,n}(\CC)^{[m]} \bigl( x\bigr) := \nu_{\delta_n^{[m]}}(x)\, ,
$$
where $\delta_n^{[m]} := [n]_* \bigl([\CC^{[m]}]\bigr) \in \CH_m\bigl(\CC^{[nm]}/S\bigr)$, and where for a class $a \in \CH_*\bigl(\CC^{[k]}/S\bigr)$ we define maps $\nu_a \colon \CH_*(\CC^{[j]}/S) \to \CH_*(\CC^{[j-k]}/S)$ by $\nu_a(x) := \pr_{2,*}\bigl(\pr_1^*(a) \cdot \alpha_{k,j-k}^*(x) \bigr)$. These operators are indeed divided powers, in the sense that $P_{0,n}(\CC)^{[1]} = P_{0,n}(\CC)$ and $P_{0,n}(\CC)^{[l]} \circ P_{0,n}(\CC)^{[m]} = \binom{l+m}{l} \cdot P_{0,n}(\CC)^{[l+m]}$, and likewise for the $P_{n,0}$. In particular, $m! \cdot P_{n,0}(\CC)^{[m]} = P_{n,0}(\CC)^{m}$ and $m! \cdot P_{0,n}(\CC)^{[m]} = P_{0,n}(\CC)^{m}$.

The most relevant for this paper is $P_{0,1}(\CC)^{[m]} \colon \CH_*(\CC^{[j]}/S) \to \CH_*(\CC^{[j-m]}/S)$, which is given by
$$
P_{0,1}(\CC)^{[m]}\bigl(x\bigr) = \pr_{2,*} \alpha_{m,j-m}^*(x)\, .
$$


To prove the properties that we shall need, it is useful to work with some of these operators directly on the level of cycles. Given a scheme $X$ over~$S$, write $\mathcal{Z}_i(X/S)$ for the group of cycles on~$X$ of dimension $d+i$. (So $i$ is the relative dimension over~$S$ and may be negative.) Let $\mathcal{Z}_*(X/S) := \oplus_i\, \mathcal{Z}_i(X/S)$.

Let $t_{m,N}\colon \CC\times \CC^{[N-m]} \to \CC\times \CC^{[N]}$ be the morphism sending $(x,D)$ to $(x,m\cdot x+D)$. Write $\mathcal{Z}_*(\Cbul/S)=\oplus_N\; \mathcal{Z}_*(\CC^{[N]}/S)$ and $\mathcal{Z}_*(\CC\times\Cbul/S) := \oplus_N\; \mathcal{Z}_*(\CC\times\CC^{[N]}/S)$. Then we can define operators $\mathcal{P}_{m,1}\colon \mathcal{Z}_*(\Cbul/S) \to \mathcal{Z}_*(\CC\times\Cbul/S)$ by
\begin{equation}\label{calPm1-def}
\mathcal{P}_{m,1}(\zeta) := (t_{m,N+m-1})_* s_{1,N}^*(\zeta)
\quad \text{for $\zeta\in \mathcal{Z}_*(\CC^{[N]}/S)$.}
\end{equation}
Note that the map $s_{1,N} = \alpha_{1,N-1}\colon \CC\times_S \CC^{[N-1]} \to \CC^{[N]}$ is flat (see Remark~1.2 in~\cite{DPCRIFT}), so $\mathcal{P}_{m,1}$ is well-defined on the level of cycles. It respects rational equivalence and therefore induces an operator $\mathcal{P}_{m,1}\colon \CH_*(\Cbul/S) \to \CH_*(\CC\times\Cbul/S)$.

The Pontryagin product makes $\mathcal{Z}_*(\Cbul/S)$ a commutative ring, and makes $\mathcal{Z}_*(\CC\times\Cbul/S)$ into a $\mathcal{Z}_*(\Cbul/S)$-module. Finally, let us observe that $P_{0,1}(C)^{[m]}$ is also well-defined on the level of cycles since the maps $\a_{m,N-m}\colon C^{[m]}\times C^{[N-m]}\to C^{[N]}$ are flat.

\begin{lem}\label{deriv-lem-NEW} 
\textup{(i)} For every $m\ge 0$ the map $\mathcal{P}_{m,1} \colon \mathcal{Z}_*(\Cbul/S) \to \mathcal{Z}_*(C\times_S \Cbul/S)$ is a derivation.

\textup{(ii)} For all $n \geq 0$ and all $\zeta_1$, $\zeta_2 \in \mathcal{Z}_*(\Cbul/S)$ we have 
\begin{equation}\label{P01[n]-zeta1zeta2}
P_{0,1}(C)^{[n]}\bigl(\zeta_1 * \zeta_2) = \sum_{\nu=0}^n\, P_{0,1}(C)^{[\nu]}\bigl(\zeta_1) * P_{0,1}(C)^{[n-\nu]}\bigl(\zeta_2)\, .
\end{equation}
\end{lem}

\begin{proof}
(i) Let $Z_1\sub \CC^{[M]}$ and $Z_2\sub \CC^{[N]}$ be closed subvarieties.
We need to check the equality of cycles 
\begin{equation}\label{Pm1-cycle-deriv-eq}
\mathcal{P}_{m,1}\bigl([Z_1]*[Z_2]\bigr) = \mathcal{P}_{m,1}(Z_1) * Z_2 + \mathcal{P}_{m,1}(Z_2) * Z_1
\end{equation}
on $C\times C^{[M+N-1+m]}$. Consider the diagram with cartesian squares
$$
\begin{matrix}
\Pi & \xrightarrow{\quad\tau\quad} & \CC \times_S \CC^{[M]} \times_S \CC^{[N]} & \xrightarrow{\ \pr_{23}\ } & C^{[M]}\times_S C^{[N]} \\
\llap{$\scriptstyle h$}\Big\downarrow & & \Big\downarrow\rlap{$\scriptstyle \id \times \alpha_{M,N}$} & & \Big\downarrow\rlap{$\scriptstyle \alpha_{M,N}$} \\[6pt]
C\times_S C^{[M+N-1]} & \xrightarrow{t_{1,M+N}} & \CC \times_S \CC^{[M+N]} & \xrightarrow{\ \pr_2\ } & C^{[M+N]}
\end{matrix}
$$
Let $\pi := \pr_{23} \circ \tau \colon \Pi \to \CC^{[M]} \times_S \CC^{[N]}$, and note that $\pr_2 \circ t_{1,M+N} = s_{1,M+N}$. The left-hand side of \eqref{Pm1-cycle-deriv-eq} is equal to the push-forward of the 
cycle of the subscheme $\pi^{-1}(Z_1\times_S Z_2)$ under the map $t_{m,M+N-1+m} \circ h \colon \Pi \to C\times C^{[M+N-1]}\to C\times C^{[M+N-1+m]}$. 

Recall that $\mathcal{D}_m \sub C\times_S C^{[m]}$ denotes the universal divisor over $C^{[m]}$. Further note that $t_{1,M+N}$ is a closed immersion whose image is precisely $\mathcal{D}_{M+N}$. Hence, $\Pi$, viewed as a closed subscheme of $C\times_S C^{[M]}\times_S C^{[N]}$, is the effective Cartier divisor $(\id\times\a_{M,N})^{-1} \mathcal{D}_{M+N}$. It follows that $\pi^{-1}(Z_1\times_S Z_2)$ is equal to the pull-back of this divisor to $C\times_S Z_1\times_S Z_2$. Now the required formula follows from the equality of Cartier divisors
$$
(\id\times\a_{M,N})^{-1} \mathcal{D}_{M+N}= \pr_{12}^{-1} \mathcal{D}_M+ \pr_{13}^{-1} \mathcal{D}_N
$$
that holds by definition of the map $\a_{M,N}$.

(ii) We have $P_{0,1}(C) = \pr_{2,*} \circ \mathcal{P}_{0,1}$, so it follows from (i) that $P_{0,1}(C)$ is a derivation, too. This implies that \eqref{P01[n]-zeta1zeta2} holds after multiplication on both sides by~$n!$. But the group of cycles has no torsion; hence \eqref{P01[n]-zeta1zeta2} holds.
\end{proof}

\begin{cor}\label{deriv-cor} 
\textup{(i)} For any $a\in\CH(\CC)$ and $n\ge 0$ the operator $P_{n,1}(a)$ on $\CH_*(\Cbul/S)$ is a derivation. If $p_0^*(a) = 0$ then $\ol{P}_{n,1}(a)$ is a derivation on $\CH_*(\Cinf/S)$.

\textup{(ii)} For all $n \geq 0$ and all $x$, $y \in \CH_*(\Cbul/S)$ we have 
$$
P_{0,1}(C)^{[n]}\bigl(x * y) = \sum_{\nu=0}^n\, P_{0,1}(C)^{[\nu]}\bigl(x) * P_{0,1}(C)^{[n-\nu]}\bigl(y)\, .
$$
\end{cor}

\begin{proof} The operator $P_{m,1}(a)$ is related to the operator $\mathcal{P}_{m,1}$ by the identity
\begin{equation}\label{Pm1-spec-eq}
P_{m,1}(a)\bigl(x\bigr) = \pr_{2*}\bigl(\pr_1^*a\cdot \mathcal{P}_{m,1}(x)\bigr)\, ,
\end{equation}
where $\pr_1 \colon \CC \times_S \Cbul \to \CC$ and $\pr_2 \colon \CC \times_S \Cbul \to \Cbul$ are the projections. It is easy to see that the map $\CH_*(\CC\times_S \Cbul/S) \to \CH_*(\Cbul/S)$ given by $y\mapsto \pr_{2*}(\pr_1^*a\cdot y)$ is a homomorphism of $\CH_*(\Cbul)$-modules. Hence (i) follows from (i) of the lemma, and (ii) is immediate from (ii) of the lemma.
\end{proof}


\section{The Chow homology of $\Cbul$}\label{CHHomCbul}

In this section we prove the core result of the paper, namely that the Chow homology of $\Cbul$ over~$S$ is isomorphic to $\CH_*(\JJ/S)\bigl[t\bigr]\bigl\langle u\bigr\rangle$. Somewhat surprisingly, we find two natural ways to define such an isomorphism. The isomorphisms $\CH_*(\JJ/S)\bigl[t\bigr]\bigl\langle u\bigr\rangle \xrightarrow{\sim} \CH_*(\Cbul/S)$ that we obtain are equal modulo~$\psi$ but in general they are different. The images of the variables $t$ and all~$u^{[m]}$ in $\CH_*(\Cbul/S)$ are the same in both cases; the difference lies in the way that $\CH_*(\JJ/S)$ is embedded into $\CH_*(\Cbul/S)$ as a subring. Taking quotients, we also obtain a second isomorphism $\CH_*(\JJ/S)\bigl\langle u\bigr\rangle \xrightarrow{\sim} \CH_*(\Cinf/S)$ that in general is different from the isomorphism~$\beta$ of Corollary~\ref{[u]-cor}.

We also give a description of the various differential operators $\partial_t^{[m]}$ and~$\partial_u$ in terms of the geometrically defined operators that we studied in Section~\ref{Operators}. Further we prove some results about how $\CH_*(\JJ/S)$ sits embedded into $\CH_*(\Cbul/S)$, which is one of the most intriguing aspects of our result.

Let us give an overview of the most important notation that we use. (The precise details are given later, in a different order than we introduce notation here.) In addition to the section~$s$ defined in~\eqref{alb-section-eq} we shall introduce another section~$s^\prime$ of the map~$\sigma_*$. Also we shall define a section~$r$ of the map~$q_*$.
$$
\xymatrix{
\CH_*(\Cbul/S) \ar[r]_-{q_*} & \CH_*(\Cinf/S) \ar@/_1pc/[l]_{r} \ar[r]_-{\sigma_*} & \CH_*(\JJ/S) \ar@/_1pc/[l] \ar@/_1pc/@<-1ex>[l]_(.4){s, s^\prime}
}
$$
Then we shall have isomorphisms $\beta$, $\gamma \colon \CH_*(\JJ/S)\langle u\rangle \xrightarrow{\sim} \CH_*(\Cinf/S)$ given by $s$ and~$s^\prime$, respectively, on $\CH_*(\JJ/S)$, and by $u^{[m]} \mapsto \bigl[\CC^{[m]}\bigr]$. Next we have $\tilde{s} := r \circ s$ and $\tilde{s}^\prime := r \circ s^\prime$, which are sections of $\tilde{\sigma}_* = \sigma_* \circ q_*$.
$$
\xymatrix{
\CH_*(\Cbul/S) \ar[r]_-{\tilde{\sigma}_*} & \CH_*(\JJ/S) \ar@/_1pc/[l] \ar@/_1pc/@<-1ex>[l]_(.4){\tilde{s}, \tilde{s}^\prime}
}
$$
We shall consider the subrings $\Kp$, $\Lp \subset \CH_*(\Cbul/S)$ with $\Kp = \Im(\tilde{s}^\prime)$ and $\Lp = \Im(\tilde{s})$. Finally we have isomorphisms $\tilde{\beta}$, $\tilde{\gamma} \colon \CH_*(\JJ/S)[t]\langle u\rangle \xrightarrow{\sim} \CH_*(\Cbul/S)$ given by $\tilde{s}$ and~$\tilde{s}^\prime$, respectively, on the coefficients, and with $t \mapsto [p_0]$ and $u^{[m]} \mapsto \bigl[\CC^{[m]}\bigr]$.
\bigskip 

Now we turn to the actual work. We start with two easy lemmas.

\begin{lem}\label{hr-isom}
Write $[p_0] \in \CH_0(\CC/S)$ for the class of $p_0(S) \subset \CC$. If $r \colon \CH_*(\Cinf/S) \to \CH_*(\Cbul/S)$ is a (homomorphic) section of~$q_*$ then the homomorphism 
$$
h_r \colon \CH_*(\Cinf/S)\bigl[t\bigr] \to \CH_*(\Cbul/S)
$$ 
given by~$r$ on the coefficients and sending $t$ to~$[p_0]$, is an isomorphism of $\CH(S)$-algebras. 
\end{lem}

\begin{proof}
It is a priori clear that the kernel of~$h_r$ is contained in the ideal $(t-1)$. Now use that $[p_0]-1$ is not a zero divisor to conclude that $h_r$ is injective. To see that $h_r$ is also surjective, we note that for any $y \in \CH_*(\Cbul/S)$ we can write
$$
y = h_r\bigl(q_*(y)\bigr) + \bigl([p_0]-1\bigr) * z
$$
for a unique $z \in \CH_*(\Cbul/S)$, as $y - h_r\bigl(q_*(y)\bigr)$ is in the kernel of~$q_*$. Further, if $y \in \oplus_{n \leq N}\, \CH_*(\CC^{[n]}/S)$ then $z \in \oplus_{n \leq N-1}\, \CH_*(\CC^{[n]}/S)$, so by induction we find that $y \in \Im(h_r)$.
\end{proof}

\begin{lem}\label{eval-lem-new}
Let $R$ be a commutative ring. Let $I \subset R$ be a nilpotent ideal; so $I^n = (0)$ for some $n > 0$. Let $M$ be an $R$-module, and let $N$ and~$N^\prime$ be direct summands of~$M$ that have the same image in $M/IM$. Then any projection $p \colon M \to N$ (with $p\big|_{N} = \id_N$) restricts to an isomorphism $N^\prime \xrightarrow{\sim} N$.
\end{lem}

\begin{proof}
Let $\alpha \colon N^\prime \to N$ be the restriction of~$p$. The assumption that $N$ is a direct summand implies that $N/IN$ maps isomorphically to the image of $N$ in $M/IM$; likewise for~$N^\prime$. Hence $\alpha$ is the identity modulo~$I$, which implies that it is surjective. Similarly, if we choose a projection $p^\prime \colon M \to N^\prime$ then $\beta := p^\prime\big|_{N} \colon N \to N^\prime$ is the identity modulo~$I$; hence $\beta \circ \alpha \colon N^\prime \to N^\prime$ differs from the identity on~$N^\prime$ by a nilpotent map. This implies that $\beta \circ \alpha$ is invertible; so $\alpha$ is injective.  
\end{proof}
\medskip

Following \cite{P-sym} we consider
$$
\Kp := \ker\bigl(P_{0,1}([p_0])\bigr) \cap\; \bigcap_{n\geq 1} \ker\bigl(P_{0,1}(\CC)^{[n]}\bigr)\sub \CH_*(\Cbul/S)\, .
$$
Note that by Cor.~\ref{deriv-cor}, $\Kp$ is a subring of $\CH_*(\Cbul/S)$.
As proven in \cite{P-sym}, Cor.~3.11, we have an isomorphism of $\CH^*(S)$-algebras
\begin{equation}\label{tu-isom-0}
\Kp\bigl[t\bigr]\bigl\langle u\bigr\rangle \xrightarrow{\ \sim\ } \CH_*(\Cbul/S)
\end{equation}
sending $t$ to~$[p_0]$ and $u^{[m]}$ to~$\bigl[\CC^{[m]}\bigr]$. Moreover, under this isomorphism the operators $\partial_t^{[n]}$ (see \ref{NotConv}) and~$\partial_u$ correspond to $P_{0,1}(\CC)^{[n]}$ and $P_{0,1}\bigl([p_0]+\psi\bigr)$, respectively. (See the proof of Proposition 3.10 in~\cite{P-sym}.) In what follows we shall use \eqref{tu-isom-0} as an identification. {}From the given description of the operators~$\partial_t^{[n]}$ we get
$$
\Kp\bigl\langle u\bigr\rangle = \cap_{n\geq 1} \ker\bigl(\partial_t^{[n]}\bigr) = \cap_{n\geq 1} \ker \bigl(P_{0,1}(\CC)^{[n]}\bigr)\, .
$$
Further, as $\CH_*(\Cinf/S)$ is the quotient of $\CH_*(\Cbul/S)$ modulo $[p_0]-1$, the composition
\begin{equation}\label{K<u>-CHCinf}
\Kp\bigl\langle u\bigr\rangle \hookrightarrow \CH_*(\Cbul/S) \xrightarrow{\ q_*\ } \CH_*(\Cinf/S)
\end{equation}
is an isomorphism. We define the section
$$
r\colon \CH_*(\Cinf/S) \to \CH_*(\Cbul/S)
$$
of the map $q_*$ as the inverse of this isomorphism. (So $\im(r) = \Kp\langle u\rangle \subset \CH_*(\Cbul/S)$.) We shall give an alternative description of~$r$ in Remark~\ref{r-alt-descr} below.

\begin{thm}\label{CHJtuThm}
Let $\tilde{s} := r \circ s\colon \CH_*(\JJ/S) \to \CH_*(\Cbul/S)$, where $s$ is the homomorphism given in~\eqref{alb-section-eq}. Then $\tilde{s}$ is the unique lifting of~$s$ with the property that $\im(\tilde{s})$ is contained in $\cap _{n\geq 1}\, \ker\bigl(P_{0,1}(\CC)^{[n]}\bigr)$. Furthermore, the map 
$$
\tilde\beta\colon \CH_*(\JJ/S)\bigl[t\bigr]\bigl\langle u\bigr\rangle \xrightarrow{\ \sim\ } \CH_*(\Cbul/S)
$$ 
restricting to~$\tilde{s}$ on $\CH(\JJ/S)$ and with $t \mapsto [p_0]$ and $u^{[m]} \mapsto \bigl[\CC^{[m]}\bigr]$, is an isomorphism of $\CH(S)$-algebras.
\end{thm}

\begin{proof}
The first assertion follows directly from the definitions. The assertion that $\tilde\beta$ is an isomorphism follows from Lemma~\ref{hr-isom} together with Cor.~\ref{[u]-cor}.
\end{proof}

Now we give the second isomorphism between $\CH_*(\JJ/S)[t]\langle u\rangle$ and $\CH_*(\Cbul/S)$.

\begin{thm}\label{KtuThm}
The homomorphism $\tilde\sigma_* \colon \CH_*(\Cbul/S) \to \CH_*(\JJ/S)$ restricts to an isomorphism $\Kp \xrightarrow{\sim} \CH_*(\JJ/S)$. Denoting by 
$\tilde{s}^\prime:\CH_*(\JJ/S)\to \Kp$ its inverse, we obtain an isomorphism 
$$
\tilde\gamma\colon \CH_*(\JJ/S)\bigl[t\bigr]\bigl\langle u\bigr\rangle \xrightarrow{\sim} \CH_*(\Cbul/S)
$$ 
restricting to $\tilde{s}^\prime$ on $\CH_*(\JJ/S)$ and with $t \mapsto [p_0]$ and $u^{[m]} \mapsto \bigl[\CC^{[m]}\bigr]$. Under this isomorphism the operators $\partial_t^{[n]}$ and~$\partial_u$ correspond to $P_{0,1}(\CC)^{[n]}$ and $P_{0,1}\bigl([p_0]+\psi\bigr)$, respectively.
\end{thm}

\begin{proof}
All we need to prove is that $\tilde\sigma_* \colon \CH_*(\Cbul/S) \to \CH_*(\JJ/S)$ restricts to an isomorphism $\Kp \xrightarrow{\sim} \CH_*(\JJ/S)$; the remaining assertions then follow from \eqref{tu-isom-0}.

Write $A \subset \CH_*(\Cinf/S)$ for the image of the section~$s$. Write $B \subset \CH_*(\Cinf/S)$ for the image of $\Kp$ under the isomorphism \eqref{K<u>-CHCinf}.  The map $s \circ \sigma_* \colon \CH_*(\Cinf/S) \to A$ is a projection, and we are done if we can show that it restricts to an isomorphism $B \xrightarrow{\sim} A$. For this we use Lemma~\ref{eval-lem-new}. Note that $A$ and~$B$ are direct summands of $\CH_*(\Cinf/S)$ as a $\CH(S)$-module. Hence we are done if we can show that $A$ and~$B$ are equal modulo~$\psi$. But we know that $B$ is the kernel of $\ol{P}_{0,1}\bigl([p_0]+\psi\bigr)$, whereas by Thm.~\ref{KV-thm} $A$ is the kernel of the operator $(\xi \cap -) = \ol{P}_{1,1}\bigl([p_0]+\psi\bigr)$. Now use that $P_{1,1}\bigl([p_0]\bigr) = P_{1,0}\bigl([p_0]\bigr) \circ P_{0,1}\bigl([p_0]\bigr)$ and note that the operator $P_{1,0}\bigl([p_0]\bigr)$ (given by the maps~$i_{n,*}$, see~\ref{Pij-exa}) induces the identity on $\CH_*(\Cinf/S)$.
\end{proof}

\begin{rem}
{}From $\tilde\gamma$ we obtain, passing to quotients modulo the ideals generated by $t-1$, respectively $[p_0] - 1$, an isomorphism 
\begin{equation}\label{gamma-isom}
\gamma\colon \CH_*(\JJ/S)\bigl\langle u\bigr\rangle \xrightarrow{\ \sim\ } \CH_*(\Cinf/S)\, .
\end{equation} 
It turns out that $\gamma$ is \emph{not}, in general, the same isomorphism as the isomorphism~$\beta$ that we obtained in Cor.~\ref{[u]-cor}. (See the next theorem.) 

If $j \colon \CH_*(\JJ/S)\langle u\rangle \hookrightarrow \CH_*(\JJ/S)[t]\langle u\rangle$ is the inclusion map and $\ev_1 \colon \CH_*(\JJ/S)[t]\langle u\rangle \to \CH_*(\JJ/S)\langle u\rangle$ is the map given by $t \mapsto 1$ then in the diagrams
$$
\begin{matrix}
\CH_*(\JJ/S)\bigl[t\bigr]\bigl\langle u\bigr\rangle & \xrightarrow[\tilde{\gamma}]{\ \sim\ } & \CH_*(\Cbul/S) \\
\llap{$\scriptstyle j$}\Big\uparrow\; \Big\downarrow\rlap{$\scriptstyle \ev_1$} && \llap{$\scriptstyle r$}\Big\uparrow\; \Big\downarrow\rlap{$\scriptstyle q_*$} \\[6pt]
\CH_*(\JJ/S)\bigl\langle u\bigr\rangle & \xrightarrow[\gamma]{\ \sim\ } & \CH_*(\Cinf/S)
\end{matrix}
\qquad\qquad
\begin{matrix}
\CH_*(\JJ/S)\bigl[t\bigr]\bigl\langle u\bigr\rangle & \xrightarrow[\tilde{\beta}]{\ \sim\ } & \CH_*(\Cbul/S) \\
\llap{$\scriptstyle j$}\Big\uparrow\; \Big\downarrow\rlap{$\scriptstyle \ev_1$} && \llap{$\scriptstyle r$}\Big\uparrow\; \Big\downarrow\rlap{$\scriptstyle q_*$} \\[6pt]
\CH_*(\JJ/S)\bigl\langle u\bigr\rangle & \xrightarrow[\beta]{\ \sim\ } & \CH_*(\Cinf/S)
\end{matrix}
$$
both the squares with upward vertical arrows and those with downward vertical arrows are commutative.

The restriction of~$\gamma$ to $\CH_*(\JJ/S) \subset \CH_*(\JJ/S)\bigl\langle u\bigr\rangle$ defines a homomorphism 
$$
s^\prime \colon \CH_*(\JJ/S) \to \CH_*(\Cinf/S)
$$ 
that is a section of~$\sigma_*$. We have $\tilde{s}^\prime = r \circ s^\prime$. If we apply Lemma~\ref{hr-isom} to the section~$r$ then we obtain an isomorphism $h_r \colon \CH_*(\Cinf/S)\bigl[t\bigr] \xrightarrow{\sim} \CH_*(\Cbul/S)$. By construction, if we identify $\CH_*(\Cinf/S)$ with $\CH_*(\JJ/S)\bigl\langle u\bigr\rangle$ via the isomorphism~$\gamma$ then $h_r$ gives the isomorphism~$\tilde\gamma$. Likewise, if we identify $\CH_*(\Cinf/S)$ with $\CH_*(\JJ/S)\bigl\langle u\bigr\rangle$ via the isomorphism~$\beta$ then $h_r$ gives the isomorphism~$\tilde\beta$.
\end{rem}

\begin{thm}\label{beta-gamma-compare}
\textup{(i)} The isomorphisms $\tilde\beta$, $\tilde\gamma \colon \CH_*(\JJ/S)\bigl[t\bigr]\bigl\langle u\bigr\rangle \xrightarrow{\sim} \CH_*(\Cbul/S)$ are equal modulo~$\psi$.

\textup{(ii)} Write $\partial_t^{[m]}$ and~$\partial_u$ for the operators on $\CH_*(\Cbul/S)$ that correspond, under the isomorphism~$\tilde\gamma$, to the differential operators $\partial_t^{[m]}$ and~$\partial_u$ on $\CH_*(\JJ/S)\bigl[t\bigr]\bigl\langle u\bigr\rangle$. Similarly, write $D_t^{[m]}$ and~$D_u$ for the operators that correspond to $\partial_t^{[m]}$ and~$\partial_u$ under the isomorphism~$\tilde\beta$. Then we have the relations 
\begin{align*}
D_t^{[m]} &= \partial_t^{[m]} = P_{0,1}(\CC)^{[m]}\, ,\\
\partial_u &= P_{0,1}\bigl([p_0]+\psi\bigr)\, ,\\
D_u &= (1+\psi u)^{-1} \cdot \bigl(\partial_u-\psi t\partial_t+\psi P_{1,1}(\CC)\bigr)\\
&= (1+\psi u)^{-1} \cdot \bigl(P_{0,1}([p_0]+\psi)-\psi t P_{0,1}(\CC) + \psi P_{1,1}(\CC)\bigr)\, .
\end{align*}

\textup{(iii)} Let $\Kp^{[n]} :=\Kp \cap \CH_*(\CC^{[n]}/S)$. Then
$$
\im(\tilde{s}^\prime) = \Kp = \dirsum_{n \geq 0}\, \Kp^{[n]}\, ,
\quad\text{and}\quad
\im(\tilde{s}) = \dirsum_{n\ge 0}\; (1+\psi u)^{-n}\cdot \Kp^{[n]}\, .
$$
\end{thm}

\begin{proof}
It will be convenient to set $\Lp := \Im(\tilde s)$, so that $\Lp[t]\langle u\rangle = \CH_*(\Cbul/S)$. Note that
\begin{equation}\label{Ku=Lu}
\Kp\langle u\rangle = \Im(r) = \Lp\langle u\rangle
\end{equation}
as subrings of $\CH_*(\Cbul/S)$. Under~$q_*$ we have $\Im(r) \xrightarrow{\sim} \CH_*(\Cinf/S)$, and as we have seen in the proof of Thm.~\ref{KtuThm} the images of $\Kp$ and~$\Lp$ in $\CH_*(\Cinf/S)$ are equal modulo~$\psi$. Hence $\Kp$ and~$\Lp$ have the same image in $\CH_*(\Cbul/S)/(\psi)$, and this implies~(i).

Next we prove~(ii). It is immediate from~\eqref{Ku=Lu} that $D_t^{[m]} = \partial_t^{[m]}$ for all~$m$. Consider the operator $\wt{D}_u := \partial_u - \psi\cdot t\partial_t + \psi\cdot P_{1,1}(\CC)$. Our goal is to show that $\wt{D}_u = (1+\psi u)\cdot D_u$. We know that $\wt{D}_u$ is a derivation, and it is easy to check that $\wt{D}_u(t) = 0$ and $\wt{D}_u(u) = 1 + \psi u$. Hence we are done if we can show that $\Lp \subset \Ker(\wt{D}_u)$. Using \cite{P-sym}, Theorems 0.1 and~3.2, it is easy to see that $\bigl[\wt{D}_u,\partial_t^{[m]}\bigr] = 0$ for all~$m$.
In particular, $\Lp\langle u\rangle = \cap_{m\geq 1}\, \Ker\bigl(\partial_t^{[m]}\bigr)$ is stable under~$\wt{D}_u$. Also the ideal generated by $(t-1)$ is stable under~$\wt{D}_u$, so $\wt{D}_u$ induces a derivation~$\ol{D}_u$ on $\CH_*(\Cinf/S)$. Under the isomorphism $\Lp\langle u\rangle \xrightarrow{\sim} \CH_*(\Cinf/S)$ the restriction of~$\wt{D}_u$ to~$\Lp\langle u\rangle$ corresponds to the operator~$\ol{D}_u$ on~$\CH_*(\Cinf/S)$. Now we use the identities
$$
\partial_u=P_{0,1}\bigl([p_0]+\psi\bigr)\, ,\quad
\partial_t=P_{0,1}(C)\, ,
\quad\text{and}\quad
t\cdot P_{0,1}\bigl([p_0]\bigr) = P_{1,1}\bigl([p_0]\bigr)\, .
$$
These allow us to rewrite $\wt{D}_u$ as $\wt{D}_u = P_{1,1}\bigl([p_0]+\psi) - (t-1)\cdot \partial_u$. It follows that $\ol{D}_u$ is the operator $\ol{P}_{1,1}\bigl([p_0]+\psi\bigr) = (\xi\cap -)$, and we know that this operator is zero on $\Im(s)$, which is the image of~$\Lp$ in~$\CH_*(\Cinf/S)$.

(iii) Consider the grading of $\CH_*(\Cbul/S)$ for which $\CH_*(\CC^{[n]}/S)$ is placed in degree~$n$. The operators $\partial_t^{[m]} = P_{0,1}(\CC)^{[m]}$ and
$\partial_u = P_{0,1}\bigl([p_0]+\psi\bigr)$ are both homogeneous for this grading, of degrees $-m$ and~$-1$, respectively. Hence $\Kp = \oplus \Kp^{[n]}$. It is easy to check that for $x\in \Kp^{[n]}$ one has $D_u\bigl((1+\psi u)^{-n}\cdot x\bigr)=0$. Hence,
$$
\sum_{n\ge 0}\, (1+\psi u)^{-n}\cdot \Kp^{[n]} \sub \bigcap_{m\geq 1}\, \ker\bigl(\partial_t^{[m]}\bigr) \cap \ker(D_u) = \bigcap_{m\geq 1}\, \ker\bigl(D_t^{[m]}\bigr) \cap \ker(D_u) = \Lp\, .
$$
Since this inclusion becomes an equality modulo~$\psi$, it is an equality.
\end{proof}

\begin{caution}
While the section $s$ is compatible with the operators $[N]_*$ (see Lemma~\ref{sN*=N*s}), this is definitely not true for the section~$r$, and therefore also not for $\tilde{s}$ and~$\tilde{s}^\prime$. In fact, $\Kp$ and $\Lp \subset \CH_*(\Cbul/S)$, are not stable under the operators~$[N]_*$. In Section~\ref{CompatFiltr} we shall prove that there is still a very interesting relation between $\tilde{s}$ (or $\tilde{s}^\prime$) and the operators~$[N]_*$.  
\end{caution}

\begin{cor}\label{top-dim-cor} 
\textup{(i)} We have $\tilde{s}^\prime\bigl(\CH_i(J/S)) \in \oplus_{n\ge i+1}\, \Kp^{[n]}$.

\textup{(ii)} We have $\tilde{s}^\prime\bigl([J]\bigr) \in \Kp^{[2g]}$. This class freely generates $\Kp\cap\CH_g(\Cbul/S)\cong \Z$. 

\textup{(iii)} We have 
$$\tilde{s}\bigl([J]\bigr)=(1+\psi u)^{-2g} \cdot \tilde{s}^\prime\bigl([J]\bigr).$$
\end{cor}

\begin{proof}(i) Recall that the grading $\Kp=\oplus_n\, \Kp^{[n]}$ is compatible with the grading by dimension. It remains to observe that $\Kp^{[n]}$ cannot have nonzero classes of relative 
dimension~$i\leq n$. (For $i<n$ this is obvious; for $i=n$ use that 
$\CH_i(\CC^{[i]}/S) = \Z \cdot u^{[i]}$.)

(ii) Since $\Kp$ projects isomorphically to $\CH_*(J/S)$ we have that
$\Kp\cap\CH_g(\Cbul/S) \cong \CH_g(J/S)\cong \Z$. Hence, $\Kp\cap\CH_g(\Cbul/S)$ coincides
with $\Kp^{[n]}\cap\CH_g(\Cbul/S)$ for some $n$ and is freely generated by 
$\tilde{s}^\prime\bigl([J]\bigr)$. Next, observe that modulo $\psi$ we have 
$\tilde{s}^\prime\bigl([J]\bigr) \equiv \tilde{s}\bigl([J]\bigr) = r(\Gamma)$. Hence, by Remark \ref{Gamma-mod-psi-rem}(i), we obtain that $n\leq 2g$.
On the other hand, since $g!\, \tilde{s}^\prime\bigl([J]\bigr) = \tilde{s}^\prime\bigl([\iota(C)]\bigr)^{*g}$, it follows from~(i) that $n\ge 2g$. Therefore, $n=2g$.

(iii) By part (iii) of Theorem~\ref{beta-gamma-compare} the right-hand side belongs to~$\Lp$. Since its push-forward to $\CH_*(J/S)$ equals $[J]$, this implies our identity.
\end{proof}

\begin{rem}\label{r-alt-descr}
We have defined the section $r \colon \CH_*(\Cinf/S) \to \Kp\langle u\rangle \subset \CH_*(\Cbul/S)$ by taking the inverse of the isomorphism $(q_*)\big|_{\Kp\langle u\rangle} \colon \Kp\langle u\rangle \xrightarrow{\sim} \CH_*(\Cinf/S)$. Under the identification $\Kp[t]\langle u\rangle = \CH_*(\Cbul/S)$ the endomorphism $r \circ q_*$ is the map $F(t,u) \mapsto F(1,u)$, which is the operator $\sum_{n\geq 0}\; (1-t)^n \, \partial_t^{[n]}$. Recall that $t$ acts by the multiplication with $[p_0]$, while $\pa_t^{[n]}$ acts by $P_{0,1}(\CC)^{[n]}$, so only finitely many terms in this sum will be nonzero when acting on $\CH_*(\CC^{[N]}/S)$ for some fixed~$N$.

Thus, if for $x\in\CH_*(\Cinf/S)$ we choose an arbitrary $y \in \CH_*(\Cbul/S)$ with $q_*(y) = x$ then we have
\begin{equation}\label{formula-for-r}
r(x)= \sum_{n\geq 0}\; \bigl(1-[p_0]\bigr)^{*n} * P_{0,1}(\CC)^{[n]}\bigl(y\bigr)\, .
\end{equation}
Note that $r(x)$ is the unique element
$$
\tilde{x} = (\tilde{x}_0,\tilde{x}_1,\tilde{x}_2,\ldots) \in \CH_*\bigl(\Cbul/S\bigr)
$$
with the properties that $q_*(\tilde{x}) = x$ and $P_{0,1}(\CC)^{[m]}\bigl(\tilde{x}_n\bigr) = 0$ for all $m,n \geq 1$.  

We also find a more explicit description of the map $\tilde{s} := r \circ s\colon \CH_*(\JJ/S) \to \CH_*(\Cbul/S)$. Namely, if we lift the class $\Ga \in\CH_g(\Cinf/S)$ to some class $\Ga(N)\in\CH_g(\CC^{[N]}/S)$ then we get, for $z \in \CH_*(\JJ/S)$, the identity
\begin{equation}\label{tilde-s-formula}
\tilde{s}(z) = \sum_{n\geq 0}\; \bigl(1-[p_0]\bigr)^{*n} * P_{0,1}(\CC)^{[n]}\bigl(\sigma_N^*(z)\cap \Gamma(N)\bigr)\, .
\end{equation}
\end{rem}

As an application, we obtain some results on how $\CH_*(\JJ/S)$ embeds into $\CH_*(\Cbul/S)$ via either $\tilde{s}$ or $\tilde{s}^\prime$. First we prove a lemma. Recall that $i_{n-s,n} \colon \CC^{[n-s]} \to \CC^{[n]}$ is the inclusion given by $D \mapsto D+s\cdot p_0$.

\begin{lem}\label{xi-power-lem}
For every $k$ with $0\le k\le n$ one has a relation 
$$
\xi_n^k = \sum_{s=0}^k\, c(n,k,s) \cdot \psi^{k-s}(i_{n-s,n})_* \bigl[\CC^{[n-s]}\bigr]\, ,
$$
in $\CH^k(\CC^{[n]})$, for some integers $c(n,k,s)$.
\end{lem}

\begin{proof}
The case $k=1$ follows from the definition of $\xi_n$:
$$
\xi_n=(i_{n-1,n})_*[\CC^{[n-1]}]+n\psi[\CC^{[n]}]\, .
$$
To deduce the general case use induction on~$k$ together with the fact that
$i_{n-1,n}^*\xi_n=\xi_{n-1}$.
\end{proof}

\begin{prop}\label{range-estimate-1}
As before, consider $\Kp = \im(\tilde{s}^\prime)$ and $\Lp := \im(\tilde{s})$. 

\textup{(i)} The class $\Ga\in\CH_g(\Cinf/S)$ can be lifted to a class in $\CH_g(C^{[N]}/S)$ for 
some $N$ if and only if $\Lp\sub\CH_*(\CC^{[\leq N]}/S)$.

\textup{(ii)} If for some $i$ and $N$ we have $\Lp\cap\CH_i(\Cbul/S)\sub\CH_*(\CC^{[\leq N]}/S)$ then 
also $\Kp\cap\CH_i(\Cbul/S)\sub\CH_*(\CC^{[\leq N]}/S)$.

\textup{(iii)} Both $\Kp$ and $\Lp$ are contained in $\CH_*(\CC^{[\leq 2g+d]}/S)$. Further, if $i \leq g-d$ then under both sections $\tilde{s}$ and~$\tilde{s}^\prime$, the image of $\CH_i(\JJ/S)$ is contained in $\CH_*(\CC^{[\leq g+2d+i]}/S)$.
\end{prop}

\begin{proof}
(i) This follows immediately from \eqref{tilde-s-formula}.

(ii) Given a nonzero element $y \in \Kp^{[j]}\cap\CH_i(\Cbul/S)$, part~(iii) of 
Thm.~\ref{beta-gamma-compare} shows that $(1+\psi u)^{-j} \cdot y$ is in 
$\Lp\cap\CH_i(\Cbul/S)\sub\CH_i(\CC^{[\leq N]}/S)$.
On the other hand, $(1+\psi u)^{-j} \cdot y = y + z$ with $z \in \CH_*(\CC^{[>j]}/S)$. 
Hence, we should have $j\le N$.

(iii) By part~(ii) it is enough to prove the assertion for~$\Lp$ and $\tilde{s}$. Let $n = 2g+d$. 
As we have seen in the proof of Thm.~\ref{Gamma-L-Prop}, we can lift the class $\Ga \in\CH_g(\Cinf/S)$ to a class $\Gamma(n)$ in $\CH_g(\CC^{[n]}/S)$. It then follows from part (i)
that $\Lp \subset \CH_*(\CC^{[\leq n]}/S)$. 

Next consider an irreducible subvariety $Z \subset \JJ$ of dimension $m\leq g$. Note that $[Z] \in \CH_{m-d}(\JJ/S)$. By \eqref{formula-for-Gamma} we have
$$
\sigma_n^*\bigl([Z]\bigr) \cap \Gamma(n) = \sum_{k\geq 0}\; \xi_n^k \cap \sigma_n^*\bigl(c_{n-g-k}(E_n) \cap [Z]\bigr) = \sum_{k\geq g+d-m}\; \xi_n^k \cap \sigma_n^*\bigl(c_{n-g-k}(E_n) \cap [Z]\bigr)
$$
where the second equality holds because $c_{n-g-k}(E_n) \cap [Z]$ can be nonzero only if $n-g-k \leq m$, i.e., $k \geq n-g-m = g+d-m$. Let $y_k := \sigma_n^*\bigl(c_{n-g-k}(E_n) \cap [Z]\bigr)$. Lemma~\ref{xi-power-lem} then gives
$$
\sigma_n^*\bigl([Z]\bigr) \cap \Gamma(n) = \sum_{k\geq g+d-m}\; \sum_{s=0}^k c(n,k,s) \cdot (i_{n-s,n})_* (i_{n-s,n})^*\bigl(\psi^{k-s} y_k\bigr)\, .
$$
But $\psi^{k-s} = 0$ for $k-s > d$, so in the second sum we can restrict to indices $s \geq k-d \geq g-m$. This shows that $s\bigl([Z]\bigr)$ can be lifted to an element in $\CH_*(\CC^{[n-g+m]}/S) = \CH_*(\CC^{[g+d+m]}/S)$, and by our description of the section~$r$ it follows that $\tilde{s}\bigl([Z]\bigr)$ lies in $\CH_*(\CC^{[\leq g+d+m]}/S)$, which is what we wanted to prove.
\end{proof}

In the proof of the proposition we have used the fact
that the class~$\Gamma$ can be lifted to a class in $\CH_g(\CC^{[2g+d]}/S)$. 
The following corollary gives the smallest~$n$ such that $\Gamma$ admits
a lifting in $\CH_g(\CC^{[2g+n]}/S)$.

\begin{cor}\label{Gamma-2g-cor} 
\textup{(i)} The class $\Gamma\in\CH_g(\Cinf/S)_{\Q}$ cannot be lifted to $\CH_g(C^{[2g-1]}/S)_{\Q}$.

\textup{(ii)} The integral class $\Gamma\in\CH_g(\Cinf/S)$ can be realized in
$\CH_g(C^{[2g+n]}/S)$ if and only if $\prod_{i=0}^n(2g+i)\cdot\psi^{n+1}=0$.
\end{cor}

\begin{proof} (i) If we could lift $\Gamma$ to $\CH_g(C^{[2g-1]}/S)_{\Q}$ then by Prop.~\ref{range-estimate-1} this would imply that
$\Kp\sub\CH_*(\CC^{[\leq 2g-1]}/S)_{\Q}$. But we know that $\Kp^{[2g]}_{\Q}\neq 0$ by 
Corollary~\ref{top-dim-cor}(ii).

(ii) Using Prop.~\ref{range-estimate-1}(i) we find that $\Ga$ can be realized in $\CH_*(C^{[2g+n]}/S)$ if an only if $\tilde{s}\bigl([J]\bigr) = r(\Gamma)$ lies in $\CH_*(C^{[\leq 2g+n]}/S)$. By Cor.~\ref{top-dim-cor}(iii), the component of $\tilde{s}\bigl([J]\bigr)$ in $\CH_*(C^{[2g+m]}/S)$ equals
$$
(-1)^m\prod_{i=0}^{m-1}(2g+i)\cdot\psi^m u^{[m]}\cdot \tilde{s}^\prime\bigl([J]\bigr)\, ,
$$
where $\tilde{s}^\prime\bigl([J]\bigr)\in\Kp^{[2g]}$. By \eqref{tu-isom-0}, this expression vanishes if and only if 
$$
\prod_{i=0}^{m-1}(2g+i)\cdot\psi^m \cdot \tilde{s}^\prime\bigl([J]\bigr)=0\, .
$$
Pushing forward to $J$ we see that this is equivalent to the vanishing of
$\prod_{i=0}^{m-1}(2g+i)\cdot\psi^m$.
\end{proof}

\begin{cor} If we work with $\Q$-coefficients, the class $\Gamma$ can be lifted
to $\CH_g(\CC^{[3g-1]}/S)_{\Q}$. Hence, $\Lp_{\Q}$ is contained in $\CH_*(\CC^{[\leq 3g-1]}/S)_{\Q}$.
\end{cor}

\begin{proof} This follows from the well-known vanishing $\psi^g=0$; see Theorem (1.1) of
\cite{Looijenga}.
\end{proof}


\section{Motivic interpretation}

In this section we reformulate the main results on the Chow homology of $\Cbul$ and $\Cinf$ in motivic language. See Thms.~\ref{tilde-beta-gamma-motiv} and \ref{beta-gamma-motiv}. We obtain this motivic interpretation by a slight generalization of the ``Manin Principle''.
\bigskip

As before, $S$ is a smooth quasi-projective variety over~$k$. We write $\Mot(S)$ for the category of Chow motives over~$S$ with respect to graded correspondences, and $\Mot_+(S) \subset \Mot(S)$ for the subcategory of effective motives. See \cite{DeMu} and~\cite{Kunne}, but note that we here consider the theory with $\Z$-coefficients. Writing $\Var(S)$ for the category of smooth projective $S$-schemes, we have a functor $R(-/S) \colon \Var(S)^\opp \to \Mot_+(S)$. If $f \colon X \to Y$ is a morphism in~$\Var(S)$ then we shall usually write $f^* \colon R(Y/S) \to R(X/S)$ for $R(f/S)$.

For $n \in \Z$, define $\unitmot(n) := (S,\id,n)$, which is the unit motive Tate-twisted by~$n$.

On $\Mot(S)$ we have a duality $M \mapsto M^\vee$, where the dual of $M = (X,p,m)$ is $M^\vee = \bigl(X,{}^t p,d-m\bigr)$ if $X$ is of relative dimension~$d$ over~$S$. We write $R_*(X/S) := R(X/S)^\vee$, and we call this the homological motive of~$X$. If $f \colon X \to Y$ is a morphism in $\Var(S)$ then we write $f_* \colon R_*(X/S) \to R_*(Y/S)$ for $R(f/S)^\vee$. (Of course, if $X$ has relative dimension~$d$ then $R_*(X/S) = R(X/S)\bigl(d\bigr)$, so the main difference between the functors $R$ and~$R_*$ is their effect on morphisms.)

Recall that a multiplicative structure on a motive $M$ is a morphism $\alpha \colon M \otimes_S M \to M$. For instance, the group law on the Jacobian~$J$ gives rise to a ``Pontryagin multiplicative structure'' on $R_*(\JJ/S)$; cf.\ \cite{Kunne}, section~(2.5). Given motives $M$ and~$N$ with multiplicative structures $\alpha$ and~$\beta$, respectively, we say that a morphism $F\colon M \to N$ is compatible with the multiplicative structures if $\beta \circ (F \otimes F) = F \circ \alpha$. 

The Chow groups of a motive $M$ are defined by $\Chow^n(M/S) := \Hom_{\Mot(S)}\bigl(\unitmot(-n),M\bigr)$. If $M$ carries a multiplicative structure then this induces the structure of a graded ring on $\Chow^*(M/S) := \oplus_m \, \Chow^m(M/S)$.

As we want to work with ind-schemes and schemes that are not of finite type over~$S$, we need to consider the ind-category $\IndMot(S)$. The $\otimes$-structure on~$\Mot(S)$ naturally extends to one on $\IndMot(S)$.

We now return to the situation as in \ref{NotConv}. We consider the objects in $\IndMot(S)$ defined by
$$
R_*(\Cbul/S) := \oplus_{n \geq 0}\, R_*(\CC^{[n]}/S)
$$
and
$$
R_*(\Cinf/S) := \indcatlim \bigl(R_*(S/S) \xrightarrow{p_{0,*}} R_*(\CC/S) \xrightarrow{i_{2,*}} R_*(\CC^{[2]}/S) \xrightarrow{i_{3,*}} R_*(\CC^{[3]}/S) \cdots\bigr)\, .
$$
(We use the notation $\indcatlim$ to avoid confusion with inductive limits taken in $\Mot(S)$; recall that inductive limits in an ind-category do not, in general, agree with inductive limits in the original category, if they exist.) The addition maps $\alpha_{m,n}$ give rise to multiplicative structures
$$
\alpha \colon R_*(\Cbul/S) \otimes R_*(\Cbul/S) \to R_*(\Cbul/S)
$$
and
$$
\alpha \colon R_*(\Cinf/S) \otimes R_*(\Cinf/S) \to R_*(\Cinf/S)\, ,
$$
referred to as the Pontryagin multiplicative structures.
We have morphism $q_* \colon R_*(\Cbul/S) \to R_*(\Cbul/S)$ and $\sigma_* \colon R_*(\Cinf/S) \to R_*(\JJ/S)$ that are compatible with the multiplicative structures. As before we define $\tilde\sigma_* := \sigma_* \circ q_*$, which is a morphism $R_*(\Cbul/S) \to R_*(\JJ/S)$.

Note that the Chow ring of the motive $R_*(\JJ/S)$ is just the Chow \emph{homology} of $\JJ$ over~$S$ with upper indexing: $\Chow_i(\JJ/S) = \Chow^{-i}\bigl(R_*(\JJ/S)\bigl)$; likewise for $\Cbul/S$ and $\Cinf/S$. The structure of a graded ring that is induced by the Pontryagin multiplicative structures on the motives is of course the one considered before.

Let $\GrAb$ denote the category of $\Z$-graded $\Z$-modules. To a motive $M$ over~$S$ we associate the functor $\omega_M \colon \Var(S)^\opp \to \GrAb$ given by $\omega_M(X) = \Chow^*(R(X/S) \otimes_S M)$. A Yoneda-type argument gives that the functor $M \mapsto \omega_M$ is fully faithful; cf.\ \cite{Scholl}, section~2.2. We need an extension of this ``Manin principle'' to (countable) direct sums of motives. First we extend the definition of~$\omega_M$ to ind-motives: If $M = \indcatlim M_i$ then we define $\omega_M(X) := \varinjlim \omega_{M_i}(X)$.

\begin{lem}[Manin Principle]\label{ManinPrinc-Lem}
Let $M = \oplus_{i\in I}\, M_i$ and $N = \oplus_{j\in J}\, N_j$ be direct sums of motives over~$S$, viewed as objects in $\IndMot(S)$. Then the natural map
\begin{equation}\label{ManinPrinc}
\Hom_{\IndMot(S)}(M,N) \rightarrow \Hom(\omega_M,\omega_N)
\end{equation}
is bijective. In particular, given an isomorphism of functors $f \colon \omega_M \xrightarrow{\sim} \omega_N$ there is a unique isomorphism of ind-motives $f_\mot \colon M \xrightarrow{\sim} N$ with $\omega(f_\mot) = f$.
\end{lem}

\begin{proof}
It suffices to prove that \eqref{ManinPrinc} is bijective if $M$ is an ordinary motive (so $\# I = 1$), as $\Hom(\oplus M_i,N) = \prod \Hom(M_i,N)$ and $\Hom(\oplus \omega_{M_i},\omega_N) = \prod \Hom(\omega_{M_i},\omega_N)$. Further, as we can invert Tate twists, we may assume that $M = (X,p,0)$ for some $X \in \Var(S)$ and some projector $p \in \Corr^0(X,X)$. Then $\Hom_{\IndMot(S)}(M,N)$ is the image of the endomorphism $F \mapsto F \circ p$ of $\Hom_{\IndMot(S)}\bigl(R(X/S),N\bigr)$, and similarly $\Hom(\omega_M,\omega_N)$ is the image of the endomorphism $f \mapsto f \circ \omega(p)$ of $\Hom(\omega_{R(X/S)},\omega_N)$. So we are reduced to the case $M = R(X/S)$, where we may further assume that $X/S$ is of relative dimension~$d$ for some~$d$. In this case the usual Yoneda argument applies. Namely, for any $N \in \IndMot(S)$ we have $\Hom_{\IndMot(S)}(M,N) = \Chow^d\bigl(R(X/S) \otimes_S N\bigr) = \omega_N^d(M)$, and given a morphism of motives $F \colon M \to N$ with associated natural transformation $f = \omega(F) \colon \omega_M \to \omega_N$, we have that $F = f(\id_M)$, viewing $\id_M$ as an element of $\omega_M^d(M)$ and $F$ as an element of $\omega_N^d(M)$.
\end{proof}

Define $R_*(J/S)\bigl[\unitmot\bigr]\bigl\langle \Tate\bigr\rangle \in \IndMot(S)$ to be the direct sum $\oplus_{i,j \geq 0}\, R_*(J/S)\bigl(-j\bigr) \cdot t^i u^{[j]}$, where $t^i$ and~$u^{[j]}$ are formal symbols, inserted for bookkeeping purposes only. We have a multiplicative structure on $R_*(J/S)\bigl[\unitmot\bigr]\bigl\langle \Tate\bigr\rangle$ induced by the Pontryagin structure on $R_*(J/S)$, and with $\bigl(t^i u^{[j]}\bigr) \otimes \bigl(t^k u^{[l]}\bigr) \mapsto \binom{j+l}{j} t^{i+k} u^{[j+l]}$. We think of $R_*(J/S)\bigl[\unitmot\bigr]\bigl\langle \Tate\bigr\rangle$ as the polynomial algebra over $R_*(J/S)$ in two independent ``variables'' $\unitmot$ (the unit motive) and~$\Tate$ (the Tate motive), where $\unitmot$ is an ordinary variable and $\Tate$ is a PD-variable. If $X$ is smooth and projective over~$S$ then $\omega_{R_*(J/S)[\unitmot]\langle \Tate\rangle}(X)$ is just the polynomial ring $\Chow_*(J_X/X)\bigl[t\bigr]\bigl\langle u\bigr\rangle$ on which the grading is given, taking into account the rule $\omega_i := \omega^{-i}$, by the natural grading on $\Chow_*(J_X/X)$ and placing $t$ and~$u$ in (lower) degrees $0$ and~$-1$, respectively.

\begin{thm}\label{tilde-beta-gamma-motiv}
We have isomorphisms
$$
\tilde\beta_\mot, \tilde\gamma_\mot \colon R_*(J/S)\bigl[\unitmot\bigr]\bigl\langle \Tate\bigr\rangle \xrightarrow{\ \sim\ } R_*(\Cbul/S)
$$
in $\IndMot(S)$, compatible with the multiplicative structures, such that the induced isomorphisms $\Chow_*(J_X/X)\bigl[t\bigr]\bigl\langle u\bigr\rangle \xrightarrow{\sim} \Chow_*(\Cbul_X/X)$ are the isomorphisms $\tilde\beta_X$ and $\tilde\gamma_X$ of Theorems \ref{KtuThm} and~\ref{CHJtuThm}, applied to $\CC_X$ over~$X$.
\end{thm}

\begin{proof}
All we need to remark is that the isomorphisms $\tilde\beta$ and $\tilde\gamma$ of Theorems \ref{KtuThm} and~\ref{CHJtuThm} are functorial, i.e., they define natural transformations $\omega_{R_*(J/S)[\unitmot]\langle\Tate\rangle} \to \omega_{R_*(\Cbul/S)}$. Now apply the Manin Principle.
\end{proof}

With obvious notation and similar proof we have an analogous conclusion for~$\Cinf/S$.

\begin{thm}\label{beta-gamma-motiv}
We have isomorphisms
$$
\beta_\mot, \gamma_\mot \colon R_*(J/S)\bigl\langle \Tate\bigr\rangle \xrightarrow{\ \sim\ } R_*(\Cinf/S)
$$
in $\IndMot(S)$, compatible with the multiplicative structures, such that the induced isomorphisms $\Chow_*(J_X/X)\bigl\langle u\bigr\rangle  \xrightarrow{\sim} \Chow_*(\Cinf_X/X)$ are the isomorphisms $\beta_X$ and $\gamma_X$ of Cor.~\ref{[u]-cor} and~\eqref{gamma-isom}, applied to $\CC_X$ over~$X$.
\end{thm}

\begin{rem} 
As in the case of the Chow groups, from the isomorphism $\tilde\gamma_\mot$ we get a new grading on the motive $R_*(J/S)$ which is different from the grading corresponding to Beauville's decomposition. In the case when $S=\Spec(k)$, where $k$ is an algebraically closed field, the corresponding decomposition of the motive of~$J$ \emph{with rational coefficients} coincides with the one constructed by Shermenev in~\cite{Shermenev}.
\end{rem}


\section{Compatibility with PD-structures}\label{PD-compatibility}

In this section we assume that $S=\Spec(k)$ where $k$ is a field. In Section~1 of~\cite{DPCRIFT} we have defined natural PD-structures on ideals of classes of positive dimension in $\CH_*(\Cbul)$ and $\CH_*(J)$. We are going to prove that our isomorphisms $\CH_*(\Cbul) \xrightarrow{\sim} \CH_*(J)[t]\langle u\rangle$ and $\CH_*(\Cinf) \xrightarrow{\sim} \CH_*(J)\langle u\rangle \xrightarrow{\sim} \CH_*(J)\langle x\rangle$ are compatible with the PD-structures.
\bigskip

First we recall the main construction of \cite{DPCRIFT}, Section~1. We consider a commutative graded monoid scheme $M = \oplus_{n \geq 0}\, M_n$ over~$k$ such that each $M_n$ is a quasi-projective $k$-scheme and such that the addition maps $\mu\colon M_m \times M_n \to M_{m+n}$ are proper. The two examples that are relevant for us here are $M = M_0 = J$ and $M = \Cbul$. On $\CH_*(M) := \oplus_{n \geq 0}\, \CH_*(M_n)$ we have a Pontryagin product making $\CH(M)$ into a commutative ring, and $\CH_{>0}(M)$ is an ideal of $\CH_*(M)$ for this ring structure. In \cite{DPCRIFT}, Section~1 we have shown that there is a natural PD-structure $\{\gamma_d\}$ on the ideal $\CH_{>0}(M)$. The idea of the construction is as follows.

If $Z \subset M_n$ is an irreducible closed subvariety then we define $\gamma_d\bigl([Z]\bigr)$, the $d$th divided power of the class~$[Z]$, to be the image of the class of the closed subscheme $\Sym^d(Z) \subset \Sym^d(M_n)$ under the iterated addition map $\Sym^d(M_n) \to M_{dn}$. Next consider a cycle $\zeta = \sum_{j=1}^r\,  n_j [Z_j]$, where the $n_j$ are integers and the $Z_j$ are mutually distinct closed subvarieties of~$M$ of positive dimensions (possibly unequal). For $d \geq 0$ define $\gamma_d(\zeta) \in \CH_*\bigl(\Sym^d(X)\bigr)$ by
$$
\gamma_d(\zeta) := \sum_{d_1 + \cdots + d_r = d}\; n_1^{d_1} \cdots n_r^{d_r} \cdot \gamma_{d_1}(Z_1) * \cdots * \gamma_{d_r}(Z_r)\, .
$$
This gives us maps $\gamma_d \colon \mathcal{Z}_{>0}(M) \to \CH_*(M)$. We prove that these maps descend to maps $\gamma_d \colon \CH_{>0}(M) \to \CH_*(M)$ that define a PD-structure on the ideal $\CH_{>0}(M)$. This PD-structure is functorial with respect to push-forward via homomorphisms. We refer to~\cite{DPCRIFT}, Section~1 for further details.

We now apply this to $\Cbul$ and~$J$. By functoriality, see~\cite{DPCRIFT}, Thm.~1.6, the homomorphism $\tilde{\sigma}_* \colon \CH_*(\Cbul) \to \CH_*(J)$ is a PD-morphism. Recall that $\tilde{\sigma}_*$ factors as a composition of two surjective ring homomorphisms, namely $q_*\colon \CH_*(\Cbul) \to \CH_*(\Cinf)$ and $\sigma_*\colon \CH_*(\Cinf) \to \CH_*(J)$.

\begin{lem}
The ideal $\Ker(q_*) \cap \CH_{>0}(\Cbul)$ is a sub-PD ideal of $\CH_{>0}(\Cbul)$. Hence, the PD-structure on $\CH_{>0}(\Cbul)$ induces a PD-structure $\ol{\gamma}$ on $\CH_{>0}(\Cinf) \subset \CH_*(\Cinf)$ such that the homomorphism $\sigma_*$ is compatible with PD-structures. 
\end{lem}

\begin{proof}
This follows from the observation that $\Ker(q_*) \cap \CH_{>0}(\Cbul)$ is the principal ideal $\bigl([p_0]-1\bigr) * \CH_{>0}(\Cbul)$. 
\end{proof}

\begin{lem}\label{PD-sections-lem}
The sections 
\begin{gather*}
s \colon \CH_*(\JJ)\to \CH_*(\Cinf)\, ,\quad
r \colon \CH_*(\Cinf)\to\CH_*(\Cbul)\, ,\\
\text{and}\quad \tilde{s} = \tilde{s}^\prime \colon \CH_*(J)\to \CH_*(\Cbul)
\end{gather*}
are PD-morphisms.
\end{lem}

Here the PD-ideals we consider are the ideals $\CH_{>0}$. Note that $\tilde{s} = \tilde{s}^\prime$ because we work over a field.

\begin{proof}
Let $(A,I,\gamma)$ be a PD-algebra and let $J \subset A$ be an ideal such that $I \cap J$ is a sub PD-ideal of~$I$. Let $\pi \colon A \to \ol{A} := A/J$ be the quotient homomorphism, and let $\ol{\gamma}$ be the induced PD-structure on $\ol{I} = \pi(I)$. Suppose we have a section $s \colon \ol{A} \to A$ such that $s(\ol{I}) \subseteq I$. Then $s$ is a PD-morphism if and only if its image is a PD-subalgebra of~$A$. Now we apply this to the quotient morphisms $q_* \colon \CH_*(\Cbul)\to \CH_*(\Cinf)$ and $\sigma_* \colon \CH_*(\Cinf) \to \CH_*(\JJ)$
and their sections $r$ and $s$, respectively. Hence to prove that $r$ and~$s$ (and hence also $\tilde{s} = r \circ s$) are PD-morphisms, it remains to prove that $\Im(r)= \cap_{m\geq 1} \ker\bigl(P_{0,1}(C)^{[m]}\bigr)$ and $\Im(s) = \Ker(\xi \cap -)$ are PD-subalgebras. 

Let $\zeta \in \CH_{>0}(\Cbul)$. We are going to prove that for all $m \geq 0$ and $d\ge 1$ we have
\begin{equation}\label{Pm1-gammadzeta}
\mathcal{P}_{m,1}\bigl(\gamma_d(\zeta)\bigr) = \gamma_{d-1}(\zeta) * \mathcal{P}_{m,1}(\zeta)
\end{equation}
and
\begin{equation}\label{P01[m]-gammadzeta}
P_{0,1}(C)^{[m]}\bigl(\gamma_d(\zeta)\bigr) = 
\sum_{\genfrac{}{}{0pt}{2}{d_0+d_1+d_2+\cdots\ =d}{d_1+2d_2+\cdots\ =m}}\; 
\gamma_{d_0}(\zeta) * \gamma_{d_1}\bigl(P_{0,1}(C)(\zeta)\bigr) * \gamma_{d_2}\bigl(P_{0,1}(C)^{[2]}(\zeta)\bigr) * \cdots\; ,
\end{equation}
where $\mathcal{P}_{m,1}$ is the operator defined in~\eqref{calPm1-def}, and where we interpret $P_{0,1}(C)^{[m]}$ on the level of cycles; see just before Lemma~\ref{deriv-lem-NEW}. By \eqref{Pm1-spec-eq} it follows from \eqref{Pm1-gammadzeta} that $P_{m,1}(a)\bigl(\gamma_d(x)\bigr) = \gamma_{d-1}(x) * P_{m,1}(a)\bigl(x\bigr)$ for all $x \in \CH_{>0}(\Cbul)$, which implies that $\Im(s)$ is a PD-subalgebra of $\CH_*(\Cinf)$. (Recall that $\xi\cap -$ is the operator $\ol{P}_{1,1}\bigl([p_0]\bigr)$.) Similarly, \eqref{P01[m]-gammadzeta} implies that $\Im(r)$ is a PD-subalgebra of $\CH_*(\Cbul)$.

To prove \eqref{Pm1-gammadzeta} and \eqref{P01[m]-gammadzeta}, note that $\mathcal{Z}_*(\Cbul)$ has no torsion; so it suffices to prove these identities after multiplication by some nonzero integer. We know that $d! \cdot \gamma_d(\zeta) = \zeta^{*d}$, so it follows from (i) of Lemma~\ref{deriv-lem-NEW} that \eqref{Pm1-gammadzeta} is correct after multiplication by~$d!$. Similarly, since $P_{0,1}(C)$ is a derivation, \eqref{P01[m]-gammadzeta} holds after multiplying both sides with $m!d!$.
\end{proof}

Now we can prove the main result of this section.

\begin{thm}\label{PD-thm}
Assume that $S=\Spec(k)$, where $k$ is a field. Then the isomorphisms 
\begin{alignat*}{2}
&h\colon \CH_*(J)\langle x\rangle \xrightarrow{\sim} \CH_*(\Cinf) 
&&\qquad \text{of Thm.~\ref{KV-thm},}\\
&\beta = \gamma\colon \CH_*(J)\langle u\rangle \xrightarrow{\sim} \CH_*(\Cinf)
&&\qquad \text{of Cor.~\ref{[u]-cor} and \eqref{gamma-isom},}\\
\intertext{and}
&\tilde\beta = \tilde\gamma\colon \CH_*(J)\bigl[t\bigr]\bigl\langle u\bigr\rangle \xrightarrow{\sim} \CH_*(\Cbul)
&&\qquad \text{of Thms.~\ref{CHJtuThm} and \ref{KtuThm},}
\end{alignat*}
are all PD-homomorphisms. Here on the terms on the left we consider the natural PD-structures on the ideals generated by $\CH_{>0}(J)$ and by all $x^{[m]}$ (resp., $u^{[m]}$) for $m \geq 1$.
\end{thm}

\begin{proof} The lemma gives that $s$, $r$, and $\tilde{s} = r\circ s =\tilde{s}^\prime$ are PD-morphisms, so it remains to check that $\bigl[C^{[m]}\bigr]=\gamma_m\bigl([C]\bigr)$ and $h(x^{[m]})=\ol{\gamma}_m(L)$ for $m \geq 1$. The first is immediate from the definitions. The second equality is equivalent to $\ol{\gamma}_m(L) = L^{[m]}$ for $m \geq 1$. This follows from the definitions using Remark~\ref{L-over-field}. 
\end{proof}

\begin{rem}\label{Gamma-form-int}
Our results imply (still working over a field) that the identities \eqref{L-formula} and \eqref{Gamma-L-eq} are valid in $\CH_*(\Cinf)$. In the proof of Cor.~\ref{L-Gamma-form} we can simply replace \eqref{siotaC*g} by the relation $s([\iota(C)])^{[g]} = s\bigl([\iota(C)]^{[g]}\bigr) = s\bigl([J]\bigr) = \Ga$, which gives \eqref{Gamma-L-eq}.
\end{rem}


\section{Beauville decomposition and Fourier duality for $\Cinf$}
\label{BD-FD-Cinf}

In this section we consider Chow groups with $\Q$-coefficients $\CH(?)_\Q:=\CH(?) \otimes_\Z \Q$. The main goal of this section is to give a Beauville decomposition for the Chow homology and Chow cohomology of~$\Cinf$ and to discuss how they are interchanged under Fourier duality.
\bigskip

By \cite{DeMu}, Thm.~(2.19), we have a Beauville decomposition $\CH^*(\JJ/S)_\Q = \oplus_{i,j}\, \CH^i_{(j)}(\JJ/S)_\Q$, where an element $\alpha \in \CH^i(\JJ/S)_\Q$ lies in $\CH^i_{(j)}(\JJ/S)_\Q$ if and only if $[n]^*(\alpha) = n^{2i-j} \cdot \alpha$ for all $n \in \Z$. Further,
\begin{equation}\label{BDec-bounds}
\text{$\CH^i_{(j)}(\JJ/S)_\Q$ can be nonzero only if $\max\{i-g,2(i-g)\} \leq j \leq \min\{2i,i+d\}$}\, .
\end{equation}
It is a conjecture of Beauville that $\CH^i_{(j)}(\JJ)_\Q = 0$ if $j < 0$, at least over a field.

In homological notation, let us define
$$
\CH_{i,(j)}(\JJ/S)_\Q := \CH^{g-i}_{(j)}(\JJ/S)_\Q \, .
$$
Then $\CH_*(\JJ/S)_\Q = \oplus_{i,j}\, \CH_{i,(j)}(\JJ/S)_\Q$, where an element $\alpha \in \CH_i(\JJ/S)_\Q$ lies in the subspace $\CH_{i,(j)}(\JJ/S)_\Q$ if and only if $[n]_*(\alpha) = n^{2i+j} \cdot \alpha$ for all $n \in \Z$. The Fourier transform restricts to a bijection between $\CH_{k,(j)}(\JJ/S)_\Q = \CH^{g-k}_{(j)}(\JJ/S)_\Q$ and $\CH_{g-k-j,(j)}(\JJ/S)_\Q = \CH^{k+j}_{(j)}(\JJ/S)_\Q$.
 
The following result gives a Beauville decomposition for the Chow homology and Chow cohomology of~$\Cinf$.

\begin{thm}\label{Beauv-dec-Cinf}
Let
$$
\CH^i_{(j)}(\Cinf/S)_\Q := 
\bigl\{\alpha \in \CH^i(\Cinf/S)_\Q \bigm| \text{$[N]^* \alpha = N^{2i-j} \cdot \alpha$ for all $N \geq 0$} \bigr\}
$$
and
$$
\CH_{i,(j)}(\Cinf/S)_\Q := 
\bigl\{\alpha \in \CH_i(\Cinf/S)_\Q \bigm| \text{$[N]_* \alpha = N^{2i+j} \cdot \alpha$ for all $N \geq 0$} \bigr\}\, .
$$
Then we have bigradings
$$
\CH^*(\Cinf/S)_\Q = \oplus_{i,j}\; \CH^i_{(j)}(\Cinf/S)_\Q
\quad\text{and}\quad
\CH_*(\Cinf/S)_\Q = \oplus_{i,j}\; \CH_{i,(j)}(\Cinf/S)_\Q\, .
$$
Identifying $\CH^*(\JJ/S)$ with a subring of $\CH^*(\Cinf/S)$ via $\sigma^*$, we have
$$
\CH^i_{(j)}(\Cinf/S)_\Q = \oplus_{n \geq 0}\; \CH^{i-n}_{(j-n)}(\JJ/S)_\Q 
\cdot \xi^n
$$
and $\CH^i_{(j)}(\Cinf/S)_\Q$ can be nonzero only if $i-g \leq j \leq i+d$. Similarly, identifying $\CH_*(\JJ/S)$ with a subring of $\CH_*(\Cinf/S)$ via the homomorphism~$s$ we have
\begin{equation}\label{CHij-relation}
\CH_{i,(j)}(\Cinf/S)_\Q = \oplus_{n \geq 0}\; \CH_{i-n,(j+n)}(\JJ/S)_\Q * L^{* n}
\end{equation}
and $\CH_{i,(j)}(\Cinf/S)_\Q$ can be nonzero only if $-i \leq j \leq g+d-i$.
\end{thm}

\begin{proof}
By (ii) of Lemma~\ref{add-effect-lem} we have $\xi \in \CH^1_{(1)}(\Cinf/S)_\Q$. Using the projection formula, it then easily follows from the definition of $L = L^{[1]}$ that $[N]_* L = N \cdot L$, i.e., $L \in \CH_{1,(-1)}(\Cinf)_\Q$. The first assertions of the theorem then readily follow from Thms.~\ref{Chow-cohom-thm} and \ref{KV-thm}. The restrictions on the pairs $(i,j)$ for which $\CH^i_{(j)}$ and $\CH_{i,(j)}$ can be nonzero follow from~\eqref{BDec-bounds}.
\end{proof}

Next we want to discuss Fourier duality for $\Cinf$. Given the isomorphisms in Theorems \ref{Chow-cohom-thm} and \ref{KV-thm}, it is clear that there is a unique isomorphism $\CH_*(\Cinf/S)_\Q \xrightarrow{\sim} \CH^*(\Cinf/S)_\Q$ that sends $L$ to~$\xi$ and that restricts to the Fourier operator on $\CH(\JJ/S)_\Q$. This isomorphism is in fact induced by an upper correspondence, as in \cite{KV}, Section~3.

\begin{thm}
Define elements $\ell \in \CH^1(\Cinf \times_S \Cinf/S)$ and $\eta \in \CH^2(\Cinf \times_S \Cinf/S)$ by
$$
\ell := (\sigma \times \sigma)^* c_1(\PP_\JJ)
\qquad\text{and}\qquad
\eta := \pr_1^*(\xi) \cdot \pr_2^*(\xi)\, ,
$$
with $\PP_\JJ$ the Poincar\'e bundle on $\JJ \times_S \JJ$.
Then $\Fourier := \exp(\ell + \eta)$ is an upper correspondence in the sense of \cite{KV}, Def.~3.2., that induces an isomorphism of $\CH(S)_\Q$-algebras 
$$
\Fourier \colon \CH_*(\Cinf/S)_\Q \xrightarrow{\sim} \CH^*(\Cinf/S)_\Q\, .
$$ 
We have $\Fourier(L) = \xi$, and if $\Fourier_\JJ$ is the Fourier transform on~$\JJ$ the diagram
$$
\begin{matrix}
\CH_*(\JJ/S)_\Q & \xrightarrow[\ \Fourier_\JJ\quad]{\sim} & \CH^*(\JJ/S)_\Q \\
\llap{$\scriptstyle s$}\Big\downarrow && \Big\downarrow\rlap{$\scriptstyle \sigma^*$} \\[10pt]
\CH_*(\Cinf/S)_\Q & \xrightarrow[\ \Fourier\quad]{\sim} & \CH^*(\Cinf/S)_\Q \\
\end{matrix}
$$
is commutative. For $x$, $y \in \CH_*(\Cinf/S)_\Q$ we have $\Fourier(x * y) = \Fourier(x) \cdot \Fourier(y)$. Further, $\Fourier \circ [N]_* = [N]^* \circ \Fourier$, and $\Fourier$ induces a bijection between the spaces $\CH_{i-n,(j+n)}(\JJ/S)_\Q * L^{* n}$ and $\CH_{(j+n)}^{i+j}(\JJ/S)_\Q \cdot \xi^n$.
\end{thm}

Note that, unlike the situation for the Jacobian, the Fourier transform does not, in general, give a bijection between direct summands $\CH_{i,(j)}(\Cinf/S)_\Q$ and summands $\CH_{(j^\prime)}^{i^\prime}(\Cinf/S)_\Q$ for some $i^\prime$ and~$j^\prime$ depending on $i$ and~$j$.

\begin{proof}
The proof of the theorem is essentially the same as in~\cite{KV}. We omit the details.
\end{proof}

Using the bigrading on $\CH_*(C^{[\infty]}/S)_\Q$ we get a simple interpretation of the class~$L$. Thm.~\ref{Beauv-dec-Cinf} gives that $\CH_{1,(j)}(\Cinf/S)_\Q$ can be nonzero only if $-1 \leq j \leq g+d-1$. If $j$ is in this range then, noting that $\psi \in \CH_{-1,(2)}$ and $L \in \CH_{1,(-1)}$, we obtain from \eqref{[C]-formula} the relation
\begin{equation}\label{[C](j)-formula}
[C]_{(j)} = \frac{\psi^{1+j} * L^{*(2+j)}}{(2+j)!} + \sum_{n\geq 0} s\left(\bigl[\iota(C)\bigr]_{(j-n)}\right) * \frac{\psi^n * L^{* n}}{n!}\, , 
\end{equation}
where $[C]_{(j)}$ denotes the component of~$[C]$ in $\CH_{1,(j)}(\Cinf/S)_\Q$ and similarly for $[\iota(C)]$. We have $[\iota(C)] \in \CH_1(\JJ/S)_\Q = \CH^{g-1}(\JJ/S)_\Q$ and by \eqref{BDec-bounds} only the summands $\CH_{1,(l)}(\JJ/S)_\Q$ with $-1 \leq l \leq \min\{2g-2,g+d-1\}$ can be nonzero. But in fact, the situation is slightly better, as $\CH_{1,(-1)}(\JJ/S)_\Q = \CH^{g-1}_{(-1)}(\JJ/S)_\Q$ is Fourier-dual to $\CH^0_{(-1)}(\JJ/S)_\Q$, which is zero. In particular, taking $j=-1$ in~\eqref{[C](j)-formula} we obtain the following result.

\begin{prop}\label{C-L-prop}
Let $[\CC]=\sum_j\; [\CC]_{(j)}$ be the decomposition of the class $[\CC] \in \CH_1(\Cinf/S)_\Q$, with $[\CC]_{(j)} \in \CH_{1,(j)}(\Cinf/S)_\Q$. Then $[\CC]_{(j)}$ can be nonzero only if $-1 \leq j \leq g+d-1$. Further, $[\CC]_{(-1)} = L$.
\end{prop}


\section{A new grading on $\CH_*(\JJ/S)$, and its relation with Beauville's decomposition}\label{CompatFiltr}

In this section we study the new grading 
$\CH_*(\JJ/S)=\oplus_{n\geq 0} \CH_*^{[n]}(\JJ/S)$ induced by
the grading of $\CH_*(\Cbul/S)$, where $\CH_*(\CC^{[n]}/S)$ is placed in degree~$n$,  
via the isomorphism $\CH_*(\JJ/S) \xrightarrow{\sim} \Kp \subset \CH_*(\Cbul/S)$.
We prove that the associated descending filtration $\Fil^\bullet$ of $\CH_*(\JJ/S)$ is stable under the operators~$[N]_*$, and that $[N]_*$ acts on $\gr^m_{\Fil}$ as multiplication by~$N^m$. It follows that $\Fil^\bullet \otimes \Q$ coincides with the descending filtration obtained from Beauville's decomposition of $\CH_*(\JJ/S)_\Q$. However, even with $\Q$-coefficients the bigrading $\CH_*(\JJ/S) = \oplus_{i,n}\, \CH_i^{[n]}(\JJ/S)$ that we obtain is different, in general, from Beauville's decomposition. Similar results are obtained for the Chow homology of $\Cinf$ over~$S$.
\bigskip

We retain the notation of Section~\ref{CHHomCbul}. As before, write $\Lp = \Image(\tilde{s})$. Given a subspace $V \subseteq \CH_*(\Cbul/S)$, write $V^{[m]} := V \cap \CH_*(\CC^{[m]}/S)$ and $V^{[\geq m]} := V \cap \CH_*(\CC^{[\geq m]}/S)$, where $\CH_*(\CC^{[\geq m]}/S) := \oplus_{n \geq m} \CH_*(\CC^{[n]}/S)$. Note that in general $V^{[\geq m]}$ is much bigger than $\oplus_{n \geq m} V^{[n]}$, but the two are equal if $V$ is a graded subspace of $\CH_*(\Cbul/S)$. 

By Thm.~\ref{beta-gamma-compare}, we have $\Kp = \oplus_{m \geq 0}\, \Kp^{[m]}$ and $\M = \oplus_{m \geq 0}\, \M^{[m]}$. Also, we have isomorphisms $\tilde\sigma_* \colon \Kp \xrightarrow{\sim} \CH_*(\JJ/S)$ and $q_* \colon \M \xrightarrow{\sim} \CH_*(\Cinf/S)$. Therefore, we can transport the above gradings of $\Kp$ and~$\M$ to get new gradings on $\CH_*(\JJ/S)$ and $\CH_*(\Cinf/S)$. Both gradings are compatible with the usual grading by dimension. Note also by Proposition \ref{range-estimate-1}(iii), we get
$\CH_*(\JJ/S)=\oplus_{n=0}^{2g+d}\, \CH_*^{[n]}(\JJ/S)$.

Consider the decreasing filtrations $\Fil^\bullet$ on $\CH_*(\JJ/S)$ and on $\CH_*(\Cinf/S)$ that are induced by this grading. More precisely, we define
$$
\Fil^m\CH_*(\JJ/S) := \tilde\sigma_*\Kp^{[\ge m]}
\qquad\text{and}\qquad
\Fil^m\CH_*(\Cinf/S) := q_* \M^{[\ge m]}\, .
$$
These are filtrations of $\CH(S)$-algebras, in the sense that $\Fil^m * \Fil^n \subseteq \Fil^{m+n}$. Because $\Fil^0\CH_* = \CH_*$, each $\Fil^m$ is an ideal. Furthermore, if we work over a field then by Theorem~\ref{PD-thm} these filtrations are compatible with the PD-structures on these algebras. 

Recall that $r\colon \CH_*(\Cinf/S) \to \CH_*(\Cbul/S)$ is the inverse map to the isomorphism
$q_* \colon \M \to \CH_*(\Cinf/S)$.
Hence for $x \in \CH_*(\Cinf/S)$ we have that $x \in \Fil^m$ if and only if $r(x) \in \M^{[\geq m]}$. 

\begin{lem} 
For all $m\ge 0$ we have $\Fil^m\CH_*(\JJ/S) = \tilde\sigma_*\M^{[\ge m]} = \tilde\sigma_* \Lp^{[\geq m]}$.
\end{lem}

\begin{proof} We have $\Kp\subset \M$, so for the first equality it suffices to verify that 
$\tilde\sigma_*(\M^{[\ge m]}) \sub \tilde\sigma_*\Kp^{[\ge m]}$. 
As $\M^{[m]} = \oplus_{i\geq 0} \Kp^{[m-i]} u^{[i]}$ it suffices 
to check that $\tilde\sigma_*(u^{[i]})$ is in $\Fil^i\CH_*(\JJ/S)$. 
We claim that in fact there is an inclusion $\CH_i(\JJ/S)\sub 
\Fil^i\CH_*(\JJ/S)$. Indeed, we have 
$$
\CH_i(\JJ/S)=\oplus_n\, \tilde\sigma_*\bigl(\Kp \cap \CH_i(\CC^{[n]})\bigr)\, ,
$$
and only terms with $n\ge i$ give a nontrivial contribution.

We have $\Lp\langle u\rangle = \M$, so for the second equality it now suffices to show that $\Fil^m \subseteq \tilde\sigma_* \Lp^{[\geq m]}$. We prove this by descending induction on~$m$. For $m \gg 0$ the claim is immediate from Prop.~\ref{range-estimate-1}. Assume then that $\Fil^{m+1} \subseteq \tilde\sigma_*\Lp^{[\geq m+1]}$. If $x \in \Kp^{[m]}$ then by Thm.~\ref{beta-gamma-compare}(iii) the element $y := (1+\psi u)^{-m} * x$ lies in $\Lp^{[\geq m]}$. Note that $x-y \in \M^{[\geq m+1]}$. Hence, using the first equality in the lemma, $\tilde\sigma_*(x) = \tilde\sigma_*(x-y) + \tilde\sigma_*(y) \in \Fil^{m+1} + \tilde\sigma_* \Lp^{[\geq m]} = \tilde\sigma_* \Lp^{[\geq m]}$. As $\Fil^m$ is spanned by $\Fil^{m+1}$ together with the classes $\tilde\sigma_*(x)$ for $x \in \Kp^{[m]}$ the assertion follows.
\end{proof}

Summing up, for an element $y \in \CH_*(\JJ/S)$ we have
$$
y \in \Fil^m \CH_*(\JJ/S)\quad \Leftrightarrow\quad
\tilde s(y) \in \CH_*(\CC^{[\geq m]}/S) \quad \Leftrightarrow\quad
\tilde s^\prime(y) \in \CH_*(\CC^{[\geq m]}/S)\, .
$$
It follows from the lemma that $\Fil^m \CH_*(\JJ/S) = q_* \Fil^m\CH_*(\Cinf/S)$ for all~$m$.

\begin{prop}
The filtration $\Fil^\bullet$ on $\CH_*(\JJ/S)$ does not depend on the choice of the base point $p_0 \in C(S)$.
\end{prop}

\begin{proof}
Let $p_0^\prime \in C(S)$ be another base point, and let $\sigma^\prime_n \colon C^{[n]} \to J$ be the associated morphism, given on points by $D \mapsto [D - n p_0^\prime]$. By the Lemma we have $\Fil^m\CH_*(\JJ/S) = \tilde\sigma_*\M^{[\ge m]}$. Note that $\M$ does not depend on the choice of the base point. Hence it suffices to show that for all $x \in \M^{[m]}$ we have $\sigma^\prime_m(x) \in \Fil^m \CH_*(J/S)$. But $\sigma^\prime_m(x) = \delta_m * \tilde\sigma_m(x)$, where $\delta_m \in \CH_0(J/S)$ is the class of the section $m\cdot (p_0 -p_0^\prime) \in J(S)$. As $\Fil^m$ is an ideal, this gives the desired conclusion. 
\end{proof}
\medskip

We say that an element $y$ in $\CH_*(\Cinf/S)$ or $\CH_*(\JJ/S)$ has \emph{coweight}~$z$ if $[N]_*(y) = N^z \cdot y$ for all~$N$. Thus, with notation as in Section~\ref{BD-FD-Cinf}, $\CH_{i,(j)}(\Cinf/S)$ and $\CH_{i,(j)}(\JJ/S)$ have coweight $2i+j$. Our goal is to prove the following compatibility between the filtrations $\Fil^\bullet$ just defined and the filtrations given by coweight.

\begin{thm}\label{filtr-Q-thm} 
\textup{(i)} For all $m \geq 0$ we have $\Fil^m \CH_*(\JJ/S)_\Q = \oplus_{2i+j\geq m}\, \CH_{i,(j)}(\JJ/S)_\Q$.

\textup{(ii)} For all $m \geq 0$ we have $\Fil^m \CH_*(\Cinf/S)_\Q = \oplus_{2i+j\geq m}\, \CH_{i,(j)}(\Cinf/S)_\Q$.
\end{thm}

This theorem is an immediate consequence of the following more precise result that takes torsion into account. 

\begin{thm}\label{filtr-Z-thm} 
The filtrations $\Fil^\bullet$ on $\CH_*(\JJ/S)$ and on $\CH_*(\Cinf/S)$ are stable under all operators $[N]_*$, and in both cases $[N]_*$ acts on $\gr^m_{\Fil}$ as multiplication by~$N^m$.
\end{thm}

The key geometric statement that we use in the proof is the following.

\begin{lem}\label{main-filtr-lem}
Let $x\in \M^{[n]}$. Then we have
$P_{0,1}(\CC)^{[j]}\bigl([N]_*x\bigr)=0$ for $j > (N-1)n$, and $P_{0,1}(\CC)^{[(N-1)n]}\bigl([N]_*x\bigr)=N^n\cdot x$.
\end{lem}

\begin{proof}
We divide the proof into some steps. Without loss of generality we may assume the base scheme~$S$ to be irreducible. Recall that $d = \dim(S)$.
\medskip

\noindent
\emph{Step 1.} Given an integer $i$ with $0 \leq i \leq Nn$, define $Y_{n,i}(N)$ by the cartesian diagram
$$
\begin{matrix}
Y_{n,i}(N) & \xrightarrow{\quad } & \CC^{[Nn-i]}\times_S\CC^{[i]} \\
\Big\downarrow && \Big\downarrow\rlap{$\scriptstyle \a_{Nn-i,i}$} \\[6pt]
\CC^{[n]} & \xrightarrow{\ \Delta_N\ } & \CC^{[Nn]}\rlap{\quad .} \\
\end{matrix}
$$
In other words, $Y_{n,i}(N)$ parametrizes triples of effective divisors 
$(D,D_1,D_2)$ on $\CC/S$, of degrees $(n,Nn-i,i)$, such that $ND=D_1+D_2$.
We need some information on irreducible components of $Y_{n,i}(N)$.

The morphism $\alpha_{Nn-i,i}$ is finite flat of degree $\binom{Nn}{i}$; see Remark~1.2 in ~\cite{DPCRIFT}. Hence the map $Y_{n,i}(N) \to \CC^{[n]}$ is finite and flat, too. In particular, all irreducible components of $Y_{n,i}(N)$ have dimension $d+n$.

Given an integer $j$ with $\max(n-i,0)\le j\le n$, consider the natural map
$$
p_{n,i,j}(N)\colon \CC^{[j]} \times_S Y_{n-j,i+j-n}(N-1)\longrightarrow Y_{n,i}(N)
$$
given on points by $(E;D,D_1,D_2)\mapsto (D+E,D_1+NE,D+D_2)$. We are going to prove that if $W \sub Y_{n-j,i+j-n}(N-1)$ is an irreducible component then
$p_{n,i,j}(N)$ maps $\CC^{[j]}\times_S W$ birationally onto some
irreducible component of $Y_{n,i}(N)$, and that all 
irreducible components of $Y_{n,i}(N)$ are obtained in this way.

Given a $\kbar$-valued point $(D,D_1,D_2)$ of $Y_{n,i}(N)$, set
$D' := \gcd(D,D_2)$. Then we can write $D= D'+E$ and $D_2=D'+D'_2$, where $E$ and $D'_2$ are disjoint. The relation $ND=D_1+D_2$ becomes $(N-1)D'+NE=D_1+D'_2$, which implies that $D_1=D'_1+NE$ for some effective divisor~$D_1^\prime$. So 
$$
(D,D_1,D_2)=p_{n,i,j}(N)(E;D',D'_1,D'_2)\, ,
$$
where $j := n - \deg(D')$. Thus, the images of the maps $p_{n,i,j}(N)$ cover $Y_{n,i}(N)$.
It remains to be shown that $p_{n,i,j}(N)$ is birational on every irreducible component. Let $U\sub \CC^{[j]}\times_S Y_{n-j,i+j-n}(N-1)$ be the open subset consisting of $(E;D,D_1,D_2)$ such that $E$ and $D_2$ are disjoint. Then the restriction of $p_{n,i,j}(N)$ to $U$ is an embedding, since on $U$ we have
$D=\gcd(D+E,D+D_2)$. Since for every irreducible component $W \sub Y_{n-j,i+j-n}(N-1)$ the intersection $U\cap \bigl(\CC^{[j]}\times_S W\bigr)$ is non-empty, the assertion follows.
\medskip

\noindent
\emph{Step 2.} Fix $s\ge 0$. The map $\pr_{1,3} \colon (D,D_1,D_2) \mapsto (D,D_2)$ realizes $Y_{n,n-s}(N)$ as a closed subscheme of $\CC^{[n]} \times_S \CC^{[n-s]}$. It readily follows from the definition of the operators $P_{0,1}(\CC)^{[m]}$ that the operator 
$$
P_{0,1}(\CC)^{[(N-1)n+s]}\circ [N]_* \colon \CH_*(\CC^{[n]}/S) \to \CH_*(\CC^{[n-s]}/S)
$$ 
is induced by the fundamental cycle $\bigl[Y_{n,n-s}(N)\bigr]$ of the subscheme $Y_{n,n-s}(N) \subset \CC^{[n]}\times_S \CC^{[n-s]}$, where we view this cycle as a correspondence from $\CC^{[n]}$ to~$\CC^{[n-s]}$.

Let $Z$ be an irreducible component of $Y_{n,n-s}(N)$. We view $Z$ as a reduced subscheme of $\CC^{[n]}\times_S \CC^{[n-s]}$. As shown in the first step, there is some $j$ with $s \leq j \leq n$ and an irreducible component $W$ of $Y_{n-j,j-s}(N-1)$ such that $p_{n,n-s,j}(N)$ gives a birational map from $\CC^{[j]}\times_S W$ to~$Z$. Note that $\pr_1 \circ p_{n,n-s,j}(N) \colon \CC^{[j]}\times_S W \to \CC^{[n]}$ equals the restriction of $\alpha_{j,n-j} \circ \pr_{1,2} \colon \CC^{[j]} \times_S \CC^{[n-j]} \times \CC^{[n-s]} \to \CC^{[n]}$ to $\CC^{[j]}\times_S W \subset \CC^{[j]} \times_S \CC^{[n-j]} \times \CC^{[n-s]}$ and that $\pr_3 \circ p_{n,n-s,j}(N) \colon \CC^{[j]}\times_S W \to \CC^{[n-s]}$
equals the composition of the projections $\CC^{[j]}\times_S W \to W \to \CC^{[n-s]}$. As we have a commutative diagram
$$
\begin{matrix}
\CC^{[n]} & \xleftarrow{\ \alpha_{j,n-j}\ } & \CC^{[j]} \times_S \CC^{[n-j]} & \xleftarrow{\ \pr_{1,2}\ } & \CC^{[j]} \times_S \CC^{[n-j]} \times_S \CC^{[n-s]} & \longhookleftarrow & \CC^{[j]} \times_S W\\
&& \Big\downarrow\rlap{$\scriptstyle \pr_2$} & \square & \Big\downarrow\rlap{$\scriptstyle \pr_{2,3}$} & \square & \Big\downarrow\rlap{$\scriptstyle \pr$} \\[6pt]
&& \CC^{[n-j]} & \xleftarrow[\ \pr_1\ ]{} & \CC^{[n-j]} \times_S \CC^{[n-s]} & \longhookleftarrow & W \\
&&&&&& \Big\downarrow\rlap{$\scriptstyle \pr$} \\[6pt]
&&&&&& \CC^{[n-s]}
\end{matrix}
$$
we find that the operator $[Z]_* \colon \CH_*(\CC^{[n]}/S) \to \CH_*(\CC^{[n-s]}/S)$ given by the correspondence~$Z$ equals the composition 
$$
\CH_*(\CC^{[n]}/S) \xrightarrow{\ P_{0,1}(\CC)^{[j]}\ } \CH_*(\CC^{[n-j]}/S) \xrightarrow{\ [W]_*\ } \CH_*(\CC^{[n-s]}/S)\, .
$$ 
But $\M = \cap_{m > 0}\, \Ker\bigl(P_{0,1}(\CC)^{[m]}\bigr)$, so if $j > 0$ we find that $[Z]_*$ is zero on~$\M^{[n]}$. If $s > 0$ then, recalling that $j \geq s$, this applies to all components~$Z$ of $Y_{n,n-s}(N)$, so we obtain the required vanishing $P_{0,1}(\CC)^{[(N-1)n+s]}\bigl([N]_*\M^{[n]}\bigr)=0$ for $s>0$.
\medskip

\noindent
\emph{Step 3.} In the case $s=0$ there is a unique irreducible component of $Y_{n,n}(N)$ that gives a nonzero contribution to our operator, namely $p_{n,n,0}\bigl(Y_{n,0}(N-1)\bigr)\cong \CC^{[n]}$, which gives the identity correspondence from $\CC^{[n]}$ to $\CC^{[n]}$. So we only
have to check that the multiplicity of $Y_{n,n}(N)$ at the diagonal component is equal to~$N^n$. It follows from the definition of $Y_{n,n}(N)$ that this multiplicity is the number of branches of the finite covering $\alpha_{(N-1)n,n} \colon \CC^{[(N-1)n]} \times_S \CC^{[n]} \to \CC^{[Nn]}$ that over the closed subscheme $\Delta_N(\CC^{[n]})\subset\CC^{[Nn]}$ lie inside the diagonal component. Write $\CC^{[n],\gen} \subset \CC^{[n]}$ for the open subscheme of divisors consisting of $n$ distinct points. Take any point $Q \in \CC^{[Nn],\gen}$ and specialize to a point $P \in \Delta_N(\CC^{[n],\gen})$. The points in the fibre of $\alpha_{(N-1)n,n}$ over~$Q$ are indexed by the subsets $I \subset \{1,\ldots,Nn\}$ with $\# I = n$. If $\tilde{Q}_I$ is the point corresponding to~$I$ then it is clear that its specialization to the fibre over~$P$ lies in the diagonal component if and only if $\#\bigl(I \cap \{kN+1,kN+2,\ldots,(k+1)N\}\bigr) = 1$ for all $k \in \{0,\ldots,n-1\}$. There are $N^n$ such sets~$I$.
\end{proof}

\begin{proof}[Proof of Theorem \ref{filtr-Z-thm}]
We have $q \circ [N] = [N] \circ q$, and the filtration $\Fil^\bullet$ on $\CH_*(\JJ/S)$ is the one induced by $\Fil^\bullet$ on $\CH_*(\Cinf/S)$ under the quotient map $q_* \colon \CH_*(\Cinf/S) \to \CH_*(\JJ/S)$. So it suffices to prove the result for $\Cinf/S$.

By definition, 
$\Fil^m \CH_*(\Cinf/S)$ is spanned by elements of the form $x=\tilde{\sigma}_*\tilde{x}$ where 
$\tilde{x}\in \M^{[n]}$ with $n\geq m$. 
Then $[N]_* \tilde{x} \in \CH_*(\CC^{[Nn]}/S)$ lifts $[N]_* x$, so by \eqref{formula-for-r} we have
$$
r\bigl([N]_* x\bigr) = \sum_{i \geq 0}\, \bigl(1-[p_0]\bigr)^{*i} * P_{0,1}(\CC)^{[i]}\bigl([N]_* \tilde{x}\bigr)\, .
$$
By Lemma~\ref{main-filtr-lem}, all the nonzero terms in the right-hand side lie in 
$\CH^{[\geq n]}$ and the component in $\CH^{[n]}$ equals $N^n\cdot\tilde{x}$.
This immediately implies the result.
\end{proof}

\begin{cor}\label{filtr-Z-cor}
Let ${\Kp}_i := \Kp \cap \CH_i(\Cbul/S)$. Then $({\Kp}_i \otimes \Q) \subset \CH_i(\CC^{[\leq \nu]}/S)$ with $\nu := \min\{2g,g+d+i\}$. In particular, $(\Kp\otimes \Q) \sub \CH_*(\CC^{[\leq 2g]}/S)_\Q$. Further, $\Kp^{[0]} = \CH_*(\CC^{[0]}/S)$ is a free $\CH(S)$-module of rank~$1$ with generator $\tilde{s}^\prime\bigl([0]\bigr)$, and $\Kp^{[2g]}_\Q$ is a free $\CH(S)_\Q$-module of rank~$1$ with generator $\tilde{s}^\prime\bigl([\JJ]\bigr)$.
\end{cor}

\begin{proof} 
The first two assertions follow immediately from the fact that the possible eigenvalues of $[N]_*$ on $\CH_i(\JJ/S)_{\Q}$ are $N^z$, where $0\le z\le \min\{2g,g+d+i\}$. (Cf.~\eqref{BDec-bounds}.) 

It follows directly from the definition of~$\Kp$ that $\CH_*(\CC^{[0]}/S) \subset \Kp^{[0]}$. The opposite inclusion is obvious, so $\Kp^{[0]} = \CH_*(\CC^{[0]}/S)$. We have $\tilde{s}^\prime\bigl([0]\bigr) = [S] = \bigl[\CC^{[0]}\bigr]$ because $[0] \in \CH_*(\JJ/S)$ is the identity element for the $*$-product. So indeed $\Kp^{[0]}$ is free of rank~$1$ over $\CH(S)$ with generator $\tilde{s}^\prime\bigl([0]\bigr)$.

For the last assertion, first remark that Thm.~\ref{filtr-Q-thm}, together with the mentioned bounds on the coweights, implies that $\tilde{s}^\prime$ restricts to an isomorphism
$$
\oplus_{2i+j=2g}\, \CH_{i,(j)}(\JJ/S)_\Q \xrightarrow{\ \sim\ } \Kp^{[2g]}_\Q
$$
from the coweight $2g$ subspace of $\CH_*(\JJ/S)_\Q$ to $\Kp^{[2g]}_\Q$. If $y \in \CH_*(\JJ/S))_\Q$ has coweight $2g$ then it has weight~$0$, which means that $[N]^*(y) = y$ for all~$N$. Taking $N=0$ this gives $y = \rho^*\bigl(0^*(y)\bigr) = 0^*(y) * [\JJ]$. (For the last identity, see Remark~\ref{psi-times-rem}.) This proves that $\oplus_{2i+j=2g}\, \CH_{i,(j)}(\JJ/S)_\Q \subseteq \CH(S) * [\JJ]$. As the opposite inclusion is clear we obtain the stated result.
\end{proof}

Let us now summarize the main conclusions of this section for $\CH_*(\JJ/S)$:
\smallskip

\noindent
(1) Transporting the grading on $\Kp$ via the isomorphism $\tilde\sigma_* \colon \Kp \xrightarrow{\sim} \CH_*(\JJ/S)$ we have a decomposition
\begin{equation}\label{newgrading}
\CH_*(\JJ/S) = \oplus_{m=0}^{2g+d}\; \CH_*^{[m]}(\JJ/S)\, .
\end{equation}
Together with the grading by relative dimension, we obtain a bigrading 
$$
\CH_*(\JJ/S) = \bigoplus_{\genfrac{}{}{0pt}{2}{-d \leq i\leq g}{0\leq m\leq \min\{g+2d+i,2g+d\}}}\, \CH_i^{[m]}(\JJ/S)\, .
$$
(For the bounds on~$m$ see Prop.~\ref{range-estimate-1}.) The decreasing filtration $\Fil^\bullet\CH_*(\JJ/S)$ associated to \eqref{newgrading} is stable under all the operators $[N]_*$, and $[N]_*$
acts on $\gr^n_{\Fil}$ as multiplication by~$N^n$. Hence, after tensoring with $\Q$ we just get
the filtration by coweight.
\smallskip

\noindent
(2) The subspaces $\CH_i^{[m]}(\JJ/S)$ with $m > \min\{2g,g+d+i\}$ are torsion. In fact, we can prove that there is a bound on the torsion depending only on $g$ and~$d$ but we shall not give the details here.
\smallskip

\noindent
(3) Tensoring \eqref{newgrading} with $\Q$ we obtain a new grading on $\CH_*(\JJ/S)_\Q$.
This grading is \emph{different}, in general, from the Beauville's decomposition. Indeed, for $0$-cycles our new grading coincides with the one obtained in \cite{P-sym}, Thm.~0.3, and as shown in loc.\ cit., Section~1, for a general curve of genus $\geq 3$ it is not the same as Beauville's.


\section{Tautological classes}\label{taut-cl-sec}

In this section we study the tautological subalgebras $\TCH_*(?)_{\Q}\sub\CH_*(?)_{\Q}$
of the Chow homology of $\Cbul$, $\Cinf$ and $\JJ$ over~$S$ with rational coefficients. By definition, these are obtained as the smallest $\CH(S)_\Q$-subalgebras that contain the image of $\CH(\CC)$ and are stable under all operators $[N]_*$. The main result is that the isomorphisms $\tilde{\beta}$ and $\tilde{\gamma}$ of Section~\ref{CHHomCbul} give rise to isomorphisms $\TCH_*(\JJ/S)_\Q\bigl[t,u\bigr] \xrightarrow{\sim} \TCH_*(\Cbul/S)_\Q$. In our calculations, certain ``modified diagonal classes'' $\Gamma_n(a)$ play a key role.
\bigskip

\begin{defi}
The tautological subrings 
\begin{align*}
\TCH_*(\Cbul/S)_\Q &\subset \CH_*(\Cbul/S)_\Q\\
\TCH_*(\Cinf/S)_\Q &\subset \CH_*(\Cinf/S)_\Q\\
\TCH_*(\JJ/S)_\Q &\subset \CH_*(\JJ/S)_\Q
\end{align*}
are defined to be the smallest $\CH(S)_\Q$-subalgebras (with respect to the Pontryagin product) that are stable under all operators~$[N]_*$ and that contain the image of $\CH_*(\CC)_\Q$ under the inclusion $\CH_*(\CC)_\Q \subset \CH_*(\Cbul/S)_\Q$, the natural map $\CH_*(\CC)_\Q \to \CH_*(\Cinf/S)_\Q$, and $\iota_* \colon \CH_*(\CC)_\Q \to \CH_*(\JJ/S)_\Q$, respectively.
\end{defi}

\begin{rems}\label{Taut-def-rems}
(i) To avoid any confusion, note that we consider Chow \emph{homology}, so the ring multiplication is the Pontryagin product. The tautological rings $\TCH_*(\Cbul/S)_\Q$ and $\TCH_*(\Cinf/S)_\Q$ are generated, as $\CH(S)_\Q$-algebras, by the classes $[n]_*(a)$ for $n \geq 1$ and $a \in \CH_*(\CC)$. To see this, note that $[N]_*$ commutes with Pontryagin product, and that $[N]_*\bigl([n]_*(a)\bigr) = [Nn]_*(a)$.

(ii) Over a field, the ring $\TCH_*(\JJ/S)_\Q$ defined here is not \emph{a priori} the same as the one defined in \cite{P-Lie} or~\cite{Mo}. However, it follows from \cite{P-sym}, Thm.~4.2, that they are the same. 

(iii) It follows from the definitions that we have surjective homomorphisms of $\CH(S)_\Q$-algebras 
$$
\TCH_*(\Cbul/S)_\Q \twoheadrightarrow \TCH_*(\Cinf/S)_\Q \twoheadrightarrow \TCH_*(\JJ/S)_\Q\, .
$$
(Cf.\ \cite{P-sym}, part~(iv) of Theorem~4.2.)
\end{rems}

Our main goal here is to prove analogues (with rational coefficients) of Theorems \ref{KV-thm} and~\ref{KtuThm} for tautological rings. Recall that we have $\CH_*(\Cbul/S) = \Kp[t]\langle u\rangle$ with $t = [p_0]$ and $u^{[m]} = \bigl[\CC^{[m]}\bigr]$. Consider the operators $\Pi_t$ and~$\Pi_u$ on $\CH_*(\Cbul/S)_\Q$ given by
\begin{align}
\Pi_t &:= \sum_{n \geq 0}\; (-1)^n \, t^n \, \partial_t^{[n]} = 1 - t\pa_t + t^2\pa_t^{[2]} - t^3 \pa_t^{[3]} + \cdots\, ,\label{Pit-def}\\
\Pi_u &:= \sum_{n \geq 0}\; (-1)^n\, u^{[n]} \partial_u^n = 1 - u\pa_u + u^{[2]}\pa_u^2 - u^{[3]}\pa_u^3 + \cdots\, .\notag
\end{align}
Then $\Pi_t$ is the projector onto $\cap_{n \geq 1}\, \Ker(\partial_t^{[n]})$ along the ideal $(t)$; in other words, it is the operator $F(t,u) \mapsto F(0,u)$. Similarly, $\Pi_u$ is the operator $F(t,u) \mapsto F(t,0)$. Since these two operators commute, $\Pi_u\circ \Pi_t$ is again a projector. Its image is precisely the subring $\Kp \subset \CH_*(\Cbul/S)$.

\begin{defi}\label{Gamman(a)-def}
For $n \geq 0$ and $a \in \CH_*(\CC/S)$, define classes $\Gamma_n(a)$ and~$\Gamma^\natural_n(a)$ in $\CH_*(\CC^{[n]}/S)$ by
$$
\Gamma_n(a) := \sum_{k = 0}^n \; (-1)^k \, \binom{n}{k}\, (t+\psi u)^k\, \Delta_{n-k,*}(a)
$$
and
$$
\Gamma_n^\natural(a) := \Gamma_n(a) + (-u)^n \psi^{n-1} p_0^*(a)\, ,
$$
with the convention that $\Gamma_0^\natural(a) = \Gamma_0(a) = \pi_*(a)$. We use the same notation $\Gamma_n(a)$ for the image of this class in $\CH_*(\Cbul/S)_\Q$. It is clear from the definitions that this class is tautological, i.e., $\Gamma_n(a) \in \TCH_*(\Cbul/S)_\Q$.
\end{defi}

For example,
\begin{align*}
\Ga_0(a) &= \pi_*(a)\, ,\\
\Ga_1(a) &= a - \bigl([p_0]+\psi[\CC]\bigr) * \pi_*(a)\, ,\\
\Gamma_2(a) &= \Delta_{2,*}(a) - 2 \bigl([p_0]+\psi[\CC]\bigr) * a + \bigl([p_0]+\psi[\CC]\bigr)^2 * \pi_*(a)\, .
\end{align*}
In the case when the base is a point the classes $\Gamma_n(\CC)$ are the
\emph{modified diagonal classes}; see \cite{GS}, and~\cite{P-sym}. For example, $\Ga_1(\CC) = [\CC]$, and modulo~$\psi$ we have
\begin{align*}
\Ga_2(\CC) &\equiv \De_{2*}\bigl([\CC]\bigr) - 2 [p_0]*[\CC]\, ,\\
\Gamma_3(\CC) &\equiv \Delta_{3,*}\bigl([\CC]\bigr) - 3[p_0]*\Delta_{2,*}\bigl([\CC]\bigr) + 3 [p_0]^{*2} * [\CC]\, .
\end{align*}

\begin{lem}\label{Deltan(a)-project}
Let $n \geq 0$. For $a \in \CH_*(\CC/S)$ we have the identity
$$
(\Pi_u\circ \Pi_t)\bigl(\Delta_{n,*}(a)\bigr) =  \Gamma^\natural_n(a)
$$
in $\CH_*(\Cbul/S)$. In particular, $\Gamma_n^\natural(a) \in \Kp$ for all $n \geq 0$ and all $a \in \CH_*(\CC/S)$. 
\end{lem}

\begin{proof}
By Thm.~3.2 of \cite{P-sym}, we have the following relation between operators on $\CH_*(\Cbul/S)$:
$$P_{0,1}(\CC)^{[m]}P_{n,0}(a)=\sum_{i\ge 0}\binom{n}{i}P_{n-i,0}(a)P_{0,1}(\CC)^{[m-i]}.$$
Hence,
$$
\pa_t^{[m]}\bigl(\De_{n,*}(a)\bigr) = P_{0,1}^{[m]}(\CC)P_{n,0}(a)(1)=
\binom{n}{m} P_{n-m,0}(a)(1) = \binom{n}{m} \De_{n-m,*}(a)\, .
$$
Thus, if we set
$$
\wt{\Ga}_n(a) := \sum_{i\ge 0}\; (-1)^i\, \binom{n}{i}\, t^i\, \De_{n-i,*}(a)
$$
then we find that $\Pi_t\bigl(\De_{n,*}(a)\bigr) = \wt{\Ga}_n(a)$. Next we observe that on $\Kp\langle u\rangle = \cap_{m\geq 1} \Ker(\partial_t^{[m]})$ we have $\pa_u|_{\Kp\langle u\rangle} = P_{0,1}\bigl([p_0]\bigr)|_{\Kp\langle u\rangle}$. Recall that $P_{0,1}\bigl([p_0]\bigr)$ is a derivation such that 
\begin{align*}
&P_{0,1}\bigl([p_0]\bigr)\bigl(\De_{n,*}(a)\bigr) = p_0^*(a)\cdot nt^{n-1}\\
\intertext{(see Example~(c) in \ref{Pij-exa}), and}
&P_{0,1}\bigl([p_0]\bigr)(t) = (\pa_u-\psi\pa_t)\bigl(t\bigr) = -\psi\, .
\end{align*}
Hence, for $n>1$ we have
\begin{align*}
\pa_u\wt{\Ga}_n(a) &= P_{0,1}\bigl([p_0]\bigr)\bigl(\wt{\Ga}_n(a)\bigr)\\
&= p_0^*(a)\cdot\sum_{i\ge 0}\; (-1)^i\, \binom{n}{i}\, (n-i)\, t^{n-1}
- \psi\cdot \sum_{i\ge 1}\; (-1)^i\, \binom{n}{i}\, i\, t^{i-1}\, \De_{n-i,*}(a)\\
&=\psi n\cdot \wt{\Ga}_{n-1}(a)\, .
\end{align*}
On the other hand, $\pa_u \wt{\Ga}_0(a) = 0$ and $\pa_u\wt{\Ga}_1(a) = p_0^*(a) + \psi\cdot \wt{\Ga}_{e,0}(a)$. It follows that for $n\ge 0$ we have
$$
(\Pi_u \circ \Pi_t)\bigl(\De_{n,*}(a)\bigr) =
(-u)^n \psi^{n-1} p_0^*(a) + \sum_{i\ge 0}\; (-1)^i\, \binom{n}{i}\, (\psi u)^i\wt{\Ga}_{n-i}(a)\, .
$$
The last sum equals $\Gamma_n(a)$, as one easily checks.
\end{proof}

We have $\Delta_{m,*}(a) = \sum_{n=0}^m \, \binom{m}{n}\, (t+\psi u)^{m-n} \Gamma_n(a)$, as one verifies by direct calculation. So,
$$
\Delta_{m,*}(a) = \left(-\sum_{n=1}^m\, \binom{m}{n}\, (t+\psi u)^{m-n} (-u)^n \psi^{n-1}\right) \cdot p_0^*(a) + \sum_{n=0}^m\, \binom{m}{n}\, (t+\psi u)^{m-n} \cdot \Gamma_n^\natural(a)
$$
gives the expression of $\Delta_{m,*}(a)$ as a polynomial in $t$ and~$u$ with coefficients in~$\Kp$.

\begin{thm}\label{Gamma-thm}
\textup{(i)} Write $\TKp := \Kp_\Q \cap \TCH_*(\Cbul/S)_\Q$, where the intersection is taken inside $\Kp_\Q[t,u] = \CH_*(\Cbul/S)_\Q$. Then $\TKp$ is the $\CH(S)_\Q$-subalgebra of $\CH_*(\Cbul/S)_\Q$ that is generated by the classes $\Gamma^\natural_n(a)$ for $n \geq 1$ and $a \in\CH(\CC)$. The isomorphism $\tilde\sigma_* \colon \Kp_\Q \xrightarrow{\sim} \CH_*(\JJ/S)_\Q$ restricts to an isomorphism $\TKp \xrightarrow{\sim} \TCH_*(\JJ/S)_\Q$.

\textup{(ii)} We have
$$
\TKp\bigl[t,u\bigr] = \TCH_*(\Cbul/S)_\Q\, .
$$
Hence, the isomorphism $\tilde\gamma \colon \CH_*(\JJ/S)_\Q\bigl[t,u\bigr] \xrightarrow{\sim} \CH_*(\Cbul/S)_\Q$ of Thm.~\ref{KtuThm} restricts to an isomorphism
$$
\Tgammatilde \colon \TCH_*(\JJ/S)_\Q\bigl[t,u\bigr] \xrightarrow{\ \sim\ } \TCH_*(\Cbul/S)_\Q\, .
$$

\textup{(iii)} As before, let $\Lp := \im(\tilde{s}) \subset \CH_*(\Cbul/S)$. Define $\TL := \Lp_\Q \cap \TCH_*(\Cbul/S)_\Q$, where the intersection is taken inside $\Lp_\Q[t,u] = \CH_*(\Cbul/S)_\Q$. Then $\tilde{s} \colon \CH_*(\JJ/S)_\Q \xrightarrow{\sim} \Lp_\Q$ restricts to an isomorphism $\TCH_*(\JJ/S)_\Q \xrightarrow{\sim} \TL$, and $\TL[t,u] = \TCH_*(\Cbul/S)_\Q$. Hence, the isomorphism $\tilde\beta \colon \CH_*(\JJ/S)_\Q\bigl[t,u\bigr] \xrightarrow{\sim} \CH_*(\Cbul/S)_\Q$ of Thm.~\ref{CHJtuThm} restricts to an isomorphism
$$
\Tbetatilde \colon \TCH_*(\JJ/S)_\Q\bigl[t,u\bigr] \xrightarrow{\ \sim\ } \TCH_*(\Cbul/S)_\Q\, .
$$
\end{thm}

\begin{proof}
Write $\ol{\TCH}_*(\JJ/S)_\Q := \TCH_*(\JJ/S)_\Q/\psi \cdot \TCH_*(\JJ/S)_\Q$. Write $\delta_n(a)$ for the class in $\ol{\TCH}_*(\JJ/S)_\Q$ represented by the image of $\Gamma^\natural_n(a)$ under~$\tilde\sigma_*$.  Let $U \subseteq \ol{\TCH}_*(\JJ/S)_\Q$ be the $\CH(S)$-subalgebra (with identity) generated by the classes $\delta_n(a)$, for $n \geq 1$ and $a \in \CH(\CC)$. We are going to show that $U = \ol{\TCH}_*(\JJ/S)_\Q$. This implies the assertions in~(i), because we know that $\tilde\sigma_* \colon \TKp \to \TCH_*(\JJ/S)_\Q$ is injective.  

Given a class $a\in \CH^*(\CC)_\Q$ and an integer~$k$, define $\cow_k(a)$ to be the component of $\iota_*(a)$ in coweight~$k$. So by definition we have $\iota_*(a) = \sum_k\, \cow_k(a)$ in $\CH_*(\JJ/S)_\Q$, with $[N]_*\bigl(\cow_k(a)\bigr) = N^k \cdot \cow_k(a)$ for all~$N$. With this notation, $\tilde\sigma_*\bigl([N]_*(a)\bigr) = [N]_*\bigl(\iota_*(a)\bigr) = \sum_{m=0}^{2g}\; N^m \cdot \cow_m(a)$. Calculating modulo~$\psi$ we then find
$$
\delta_n(a) =
\begin{cases}
{\displaystyle -p_0^*(a) \bigl[\iota(\CC)\bigr] + \sum_{m=0}^{2g}\; S(m,1) \cdot \cow_m(a) }& \text{if $n=1$,}\\
{\displaystyle n! \cdot \sum_{m=0}^{2g}\; S(m,n) \cdot \cow_m(a) }& \text{if $n > 1$,}
\end{cases}
$$
where $S(m,n)$ denotes the Stirling number of the second kind. Note that $S(m,n) = 0$ if $n>m$, and $S(m,m) = 1$. Letting $n$ run from $2g$ down to~$2$ we obtain that all classes $\cow_k(a)$ with $2\leq k\leq 2g$ are in~$U$. 

As we have seen just before Prop.~\ref{C-L-prop} (or by Corollary \ref{top-dim-cor}(i)), 
the class $\bigl[\iota(C)\bigr]$ has no components in coweight $< 2$. So taking $a = [\CC]$ we find that $\bigl[\iota(\CC)\bigr] = \sum_{k=2}^{2g}\, \cow_k(\CC) \in U$. Hence, taking $n=1$, we find that also all classes $\cow_1(a)$ are in~$U$. Finally, if $x \in \CH_*(\JJ/S)_\Q$ is any class in coweight~$0$ then this means that $[N]_*(x) = x$ for all~$N$, so in particular $x = [0]_*(x) = 0_*\bigl(\rho_*(x)\bigr)$ is in the image of $\CH(S)$. Hence also all classes $\cow_0(a)$ are in~$U$. 

This proves that all classes $\cow_k(a)$ are in~$U$. But then $U$ contains the image of $\iota_*$ and is stable under all operators~$[N]_*$. So $U = \ol{\TCH}_*(\JJ/S)_\Q$, and it follows that $\TKp$ is generated by the classes $\Gamma^\natural_n(a)$ and that $\tilde\sigma_*$ restricts to an isomorphism $\TKp \xrightarrow{\sim} \TCH_*(\JJ/S)_\Q$.

For part (ii) note that the classes $t = [p_0]$ and $u=[\CC]$ are clearly tautological. The claim that $\TKp\bigl[t,u\bigr] = \TCH_*(\Cbul/S)_\Q$ is then an immediate consequence of \cite{P-sym}, Thm.~4.2(i), which gives that $\TCH_*(\Cbul/S)_\Q$ is stable under the operators $\partial_t = P_{0,1}(\CC)$ and $\partial_u = P_{0,1}\bigl([p_0]+\psi\bigr)$. The rest of~(ii) is straightforward.

For (iii) it suffices to show that $\TL[t,u] = \TKp[t,u]$. The inclusion ``$\subseteq$'' is clear. For the opposite inclusion, write $\TKp^{[n]} := \TKp \cap \Kp^{[n]}$. By (iii) of Thm.~\ref{beta-gamma-compare} we have $(1 + \psi u)^{-n} \cdot \TKp^{[n]} \subseteq \TL$. Hence $\TKp^{[n]} \subseteq \TL[t,u]$ for all~$n$, and we are done.
\end{proof}

\begin{cor}\label{KV-taut-thm}
The isomorphisms $\beta$, $\gamma \colon \CH_*(\JJ/S)[u] \xrightarrow{\sim} \CH_*(\Cinf/S)$ restrict to isomorphisms $\Tbeta$, $\Tgamma \colon \TCH_*(\JJ/S)_\Q[u] \xrightarrow{\sim} \TCH_*(\Cinf/S)_\Q$. 
\end{cor}

The following result gives an important connection between the modified diagonal classes $\Gamma_n(C)$ and the homomorphism $\tilde{s}\colon \CH_*(J/S)\to\CH_*(\Cbul/S)$.

\begin{prop}\label{Gamma-comp-prop} 
One has the following relation in $\CH_*(\Cbul/S)_{\Q}$:
$$
\tilde{s}\left(\frac{\log\bigl(1+\psi[\iota(C)]\bigr)}{\psi}\right)= 
\sum_{n\geq 2}\, (-1)^n(1+\psi u)^{-n}\, \frac{\Gamma_n^\natural(C)}{n}\, .
$$
\end{prop}

As $\Gamma^\natural_n(C) \in \Kp^{[n]}$ for all $n \geq 2$, the sum is finite. Modulo~$\psi$ (e.g., working over a field) we find that in $\CH_*(\Cbul/S)_{\Q}/(\psi)$ we have
$$
\tilde{s}\bigl[\iota(C)\bigr] \equiv \sum_{n\geq 2}\, (-1)^n\, \frac{\Gamma_n(C)}{n} \bmod (\psi)\, .
$$

\begin{proof} By Theorem \ref{beta-gamma-compare}(iii), 
the right-hand side belongs to $\Lp_\Q$, so it is enough to check that
$$
\sum_{n\geq 2}\, (-1)^n(1+\psi [\iota(C)])^{-n}\, 
\frac{\tilde{\sigma}_*\bigl(\Gamma_n(C)\bigr)+\bigl(-[\iota(C)])^n\psi^{n-1}}{n}=
\frac{\log\bigl(1+\psi[\iota(C)]\bigr)}{\psi}
$$
in $\CH_*(J/S)_{\Q}$. Let us write $c :=\bigl[\iota(C)\bigr]$ and  $d_n :=\tilde{\sigma}_*\bigl(\Gamma_n(C)\bigr)$. Since
$$
\sum_{n\geq 2}\, (-1)^n(1+\psi c)^{-n}\, \frac{(-c)^n\psi^{n-1}}{n}=
-\frac{\log\bigl(1-\frac{\psi c}{1+\psi c}\bigr)}{\psi}-\frac{c}{1+\psi c}=\frac{\log(1+\psi c)}{\psi}-\frac{c}{1+\psi c}\, ,
$$
the identity that we want to establish is equivalent to
$$
\sum_{n\geq 2}\, (-1)^n(1+\psi c)^{-n}\, \frac{d_n}{n}=\frac{c}{1+\psi c}\, ,
$$
which we can rewrite as
\begin{equation}\label{C-Gamma-push-eq}
\sum_{n\geq 1}\, (-1)^n(1+\psi c)^{-n}\, \frac{d_n}{n}=0\, .
\end{equation}
Recalling the definition of $\Gamma_n(C)$, we can write the generating series for the classes~$d_n$ in the following form (where $x$ is a formal variable):
$$
\sum_{n\geq 1}\, d_n\, \frac{x^n}{n!} =
\sum_{k\ge 0,m\ge 1}\, (-1)^k(1+\psi c)^k\cdot [m]_*(c)\, \frac{x^{k+m}}{k!m!}=
\exp\bigl(-x(1+\psi c)\bigr)\cdot \bigl(\sum_{m\ge 1}\, [m]_*(c)\, \frac{x^m}{m!}\bigr)\, .
$$
Hence, we have
$$
\bigl(\sum_{n\geq 1}\, d_n\, \frac{x^n}{n!}\bigr)\cdot \exp\bigl(x(1+\psi c)\bigr)=
\sum_{m\ge 1}\, [m]_*(c)\, \frac{x^m}{m!}\, .
$$
Comparing the coefficients of $x^m$ we get
$$
[m]_*(c) =\sum_{n\geq 1}\, \binom{m}{n} (1+\psi c)^{m-n}\cdot d_n\, .
$$
Note that both sides are polynomials in~$m$. 
Since the class $c=\bigl[\iota(C)\bigr]$ has no components in coweight $<2$ 
(by Corollary \ref{top-dim-cor}(i)), the coefficient of~$m$ in $[m]_*(c)$ is zero. Calculating the coefficient of~$m$ in the right-hand side of the above equality we get \eqref{C-Gamma-push-eq}.
\end{proof}


\section{Some relations between tautological classes}

In this section we show how our techniques can be applied to derive explicit relations between tautological classes on the Jacobian. In particular, under the assumption that the curve has a $g^r_d$ we obtain some vanishing relations, both for classes on $\Cbul$ and on the Jacobian. On the Jacobian we recover, working modulo algebraic equivalence, results of Herbaut~\cite{Her} and van der Geer-Kouvidakis \cite{GK} that extend an earlier result of Colombo-van Geemen \cite{ColGeem}. Working modulo rational equivalence this result was obtained by one of us in~\cite{Mo}. The nice feature of our approach is that the assumption about the existence of a~$g^r_d$ can be translated directly into a statement about classes on~$\Cbul$ (Lemma~\ref{GK-lem}), from which the vanishing result follows by  
a short calculation. 
\bigskip

We first introduce some notation. Throughout, we assume that $S = \Spec(k)$, 
where $k$ is a field. Given an effective divisor $D$ of degree~$d$ on~$C$, let $[D] \in \CH_0(C^{[d]})$ be the class of the point in $C^{[d]}$ corresponding to~$D$. We define
$$
e_i(D) := P_{0,1}(C)^{[d-i]}\bigl([D]\bigr) \in \CH_0(C^{[i]})\, .
$$
Note that if $D = p_1+\cdots+p_d$ for some points $p_i \in \CC(k)$, the element $e_i(D)$ is simply the $i$th elementary symmetric function of the classes $\bigl([p_i]\bigr)$ with respect to the Pontryagin product on $\CH_*(\Cbul)$. Indeed, this follows immediately from Lemma~\ref{deriv-lem-NEW}, part~(ii). In general, we have $e_0(D)=1$ and $e_d(D)=[D]$.
We also set
\begin{equation}\label{eibardef}
\ov{e}_i(D) := \Pi_t\bigl(e_i(D)\bigr) = \sum_{j=0}^i\, (-1)^j \binom{d-i+j}{j} e_{i-j}(D)t^j\in\CH_0(C^{[i]})\, ,
\end{equation}
where $\Pi_t$ is the projector defined in \eqref{Pit-def}. In the case $D=p_1+\cdots+p_d$, the class $\ov{e}_i(D)$ is the $i$th elementary
symmetric function of the classes $([p_i]-[p_0])_{i= 1,\ldots,d}$. 

If $D_1$ and $D_2$ are rationally equivalent then viewed as points of $\CC^{[d]}$ they lie in the same fibre of the map $\sigma_d \colon \CC^{[d]} \to \JJ$. But the fibres of~$\sigma_d$ are projective spaces, so $[D_1] = [D_2]$. Hence the classes $e_i(D)$ and $\ov{e}_i(D)$ only depend on the rational equivalence class of~$D$. 

Note that the isomorphism $\CH_*(\Cbul) \cong \Kp[t]\langle u\rangle$ restricts to an isomorphism $\CH_0(\Cbul) \cong \Kp_0[t]$ with $\Kp_0 := \Kp \cap \CH_0(\Cbul)$. In particular, $\ov{e}_i(D)\in\Kp$. Further it is easy to see that in $\CH_0(C^{[d]})$ we have
\begin{equation}\label{sym-fun-eq}
[D] = \sum_{i=0}^d\, \ov{e}_i(D)t^{d-i}\, ;
\end{equation}
so this gives another way to think of the classes $\ov{e}_i(D)$.

Finally, we denote by $\eps_i(D)$ the $i$th component of $\tilde{\sigma}_*\ov{e}_i(D) \in \CH_0(J)$ with respect to Beauville's decomposition $\CH_0(J)_{\Q} = \oplus_j\, \CH_{0,(j)}(J)$. 
If $D=p_1+\cdots+p_d$ for points $p_i \in \CC(k)$ then $\eps_i(D)$ is the $i$th elementary symmetric function of
the classes $\bigl(\alpha(p_i)\bigr)$, where $\alpha(p) \in \CH_{0,(1)}(J)$ is the component of $[\iota(p)] - [0] \in \CH_0(J)$. 

\begin{thm}\label{GK-thm} 
Let $C$ be a curve over a field~$k$ with a $k$-rational point~$p_0$.
Assume that $D\in C^{[d]}(k)$ is an effective divisor of degree~$d$ with $h^0(D)> 1$. Write $r(D) := h^0(D) -1$.

\textup{(i)} For integers $\nu$ and $s$, define $U(\nu,s) \in \CH_s(\CC^{[\nu]})_\Q$ by 
$$
U(\nu,s) := \sum_{\genfrac{}{}{0pt}{2}{n_1, n_2,\ldots,n_s \geq 2}{n_1+n_2+\cdots+ n_s = \nu}}\,  \frac{\Gamma_{n_1}(\CC) * \cdots * \Gamma_{n_s}(\CC)}{n_1\cdots n_s}\, ,
$$
with the convention that $U(0,0) = [\Spec(k)]$. Then for all $s\leq r(D)$ and $N > d-r(D)+s$,
\begin{equation}\label{GK-id1}
\sum_{i=0}^{N-2s}\;  (-1)^i\, \ov{e}_i(D) * U(N-i,s) = 0\qquad \text{in $\CH_*(\Cbul)_{\Q}$.}
\end{equation}

\textup{(ii)} For integers $\nu$ and $s$, define an element $\Upsilon(\nu,s) \in \CH_*(\JJ)_\Q$ by
$$
\Upsilon(\nu,s) := \sum_{\genfrac{}{}{0pt}{2}{n_1, n_2,\ldots,n_s \geq 2}{n_1+n_2+\cdots+ n_s = \nu}}\,  (n_1-1)!\cdots (n_s-1)! \; \cow_{n_1}(\CC) * \cdots * \cow_{n_s}(\CC)\, ,
$$
with the convention that $\Upsilon(0,0) = [S]$, and where we recall that $\cow_n(\CC)$ is the component of $\bigl[\iota(C)\bigr]$ in $\CH_{1,(n-2)}(\CC)$. (So $\cow_n(\CC)$ has coweight~$n$.) Then for all $s \leq r(D)$ and $N > d-r(D)+s$, we have
\begin{equation}\label{GK-id2}
\sum_{i=0}^{N-2s}\, (-1)^i\, \eps_i(D) * \Upsilon(N-i,s) = 0\qquad \text{in $\CH_*(\JJ)_\Q$.}
\end{equation}
\end{thm}

Note that for $d \leq 2g-2$ (the only case of interest) we have $N-2s > d-2r(D) \geq 0$ by Clifford's theorem. The idea for the proof is to use the representation of the divided powers of the class $L\in\CH_1(\Cinf)$ by projective spaces in the symmetric powers of~$C$.

\begin{lem}\label{GK-lem} 
Let $C$ and $D$ be as in Theorem~\ref{GK-thm}. Then the class 
$[D]*L^{[r(D)]}\in \CH_*(\Cinf)$ can be realized in $\CH_*(C^{[d]})$.
\end{lem}

\begin{proof} 
For $N \gg 0$ we know that $L^{[m]}$ is represented by the class of an $m$-dimensional linear subspace in the fiber of the map $\sigma_N\colon C^{[N]}\to J$ over $0\in J$. (See 
Remark \ref{L-over-field}.) 
Hence, $[D] * L^{[r(D)]}$ is the class of an $r(D)$-dimensional linear subspace in the fiber of $\sigma_N$ over $\sigma_d(D)\in J$. But $|D| \sub C^{[d]}$, viewed as a subvariety of $\CC^{[N]}$ via the embedding $i_{d,N}\colon \CC^{[d]} \hookrightarrow \CC^{[N]}$, is such a subspace.
\end{proof}

\begin{proof}[Proof of Theorem \ref{GK-thm}]
Recall from \eqref{[C]-formula} that (over a field) we have the relation $L=[C]-s[\iota(C)]$ in $\CH_*(\Cinf)$. To avoid notational confusion, let us set $m = r(D)$. Lemma \ref{GK-lem} gives that the class $[D] * \bigl([C]-s[\iota(C)]\bigr)^{*m} \in \CH_*(\Cinf)$ can be realized in $\CH_*(C^{[d]})$. By Remark~\ref{r-alt-descr}, it follows that the image of this class under $r\colon \CH_*(\Cinf)\to \CH_*(\Cbul)$ is an element of $\CH_*(C^{[\leq d]})$.

Viewing $[D]$ as an element of $\CH_*(\Cinf)$ we have
$$
r\bigl([D]\bigr) = \sum_{n\ge 0}\, (1-t)^n\cdot \partial_t^{[n]}[D]
= \sum_{n=0}^d\, e_{d-n}(D)\cdot (1-t)^n = \sum_{i=0}^d\, \ov{e}_i(D)\, ,$$
where for the first equality we use \eqref{formula-for-r}. Hence,
\begin{equation}\label{expr-in-Cleqd}
\sum_{i=0}^d\, \ov{e}_i(D) * \bigl(u-\tilde{s}[\iota(C)]\bigr)^{*m}
\end{equation}
is an element of $\CH_*(C^{[\leq d]})$.

It follows from Prop.~\ref{Gamma-comp-prop} that the component of $\tilde{s}[\iota(C)]^{*a}$ in $\CH_*(C^{[b]})$ equals $(-1)^b \cdot U(b,a)$. With $s$ and~$N$ as in the statement of the theorem, we know that the component of the expression~\eqref{expr-in-Cleqd} in $\CH_*(\CC^{[N+m-s]})$ is zero. Direct calculation gives that this component equals
$$
\sum_{j=0}^m \sum_{i=0}^d\, (-1)^{N-i} \binom{m}{j}\, \ov{e}_i(D) * u^{m-j} * U(N+j-s-i,j)\, .
$$
But all elements $\ov{e}_i(D)$ and $U(N+j-s-i,j)$ lie in $\Kp[t]$ and the powers of~$u$ are linearly independent over $\Kp[t]$. Hence,
$$
\sum_{i=0}^d\, (-1)^{N-i} \binom{m}{j} \ov{e}_i(D) * U(N+j-s-i,j) = 0 \qquad \text{for all $j \leq m$.}
$$
Taking $j=s$ and noting that $U(N-i,s) = 0$ if $N-i < 2s$, gives (i) of the theorem.

By Thm.~\ref{filtr-Q-thm}, $\tilde{\sigma}_* \ov{e}_i(D)$ only has components in coweight $\geq i$. By definition, the component in coweight~$i$ is $\eps_i(D)$. Further, as we have seen in the proof of Thm.~\ref{Gamma-thm}, $\tilde{\sigma}_* \Gamma_n(\CC) = n! \cdot \sum_{l=0}^{2g}\, S(l,n) \cow_n(C)$. Now we again use that $S(l,n) = 0$ if $n>l$ and $S(n,n) = 1$. It follows that $\tilde{\sigma}_* U(\nu,s)$ has components only in coweight $\geq \nu$, and that the component in coweight~$\nu$ is exactly $\Upsilon(\nu,s)$. With these remarks, (ii) follows from~(i) by pushing forward to the Jacobian and taking the component in coweight~$N$.
\end{proof}

\begin{rem}
Part (ii) of the theorem is Thm.~4.6 of \cite{Mo}. (The result is stated there in Fourier-dual form.)
\end{rem}

Let $A_*(J)_{\Q}$ denote the quotient of $\CH_*(J)_{\Q}$ modulo algebraic equivalence. Considering the identity \eqref{GK-id2} for $s=r(D)$ and using the fact that $\eps_i(D) \sim_\alg 0$ for $i>0$ we recover the following relations obtained in~\cite{GK}. (Equivalent identities were first derived by Herbaut in~\cite{Her}.)

\begin{cor}[Herbaut, van der Geer-Kouvidakis]\label{GK-cor1}
Let $C$ be a curve over a field $k$ with a $k$-rational point~$p_0$.
Assume that $C$ admits a $g^r_d$ defined over~$k$. Then for every $N>d$ one has
$$
\sum_{n_1+\cdots+n_r=N}(n_1-1)!\cdots (n_r-1)!\; \cow_{n_1}(\CC) *\cdots * \cow_{n_r}(\CC) = 0
\qquad \text{in $A_*(J)_{\Q}$.}
$$
\end{cor}

On the other hand, for $s=0$ \eqref{GK-id1} and \eqref{GK-id2} give the following vanishing result.

\begin{cor}\label{GK-cor2} 
Let $C$ be a curve of genus $\ge 1$. With the same assumptions as in Theorem~\ref{GK-thm}, for $d-r(D) < N \le d$ one has $\ov{e}_N(D)=0$ in $\CH_0(C^{[N]})_{\Q}$ and $\eps_N(D)=0$ in $\CH_0(J)_{\Q}$.
\end{cor}

The case $N=d$ gives the following result that generalizes the well-known property of Weierstrass points on a hyperelliptic curve (see~\cite{BV}, Prop.~3.2).

\begin{cor}\label{GK-cor3} 
Let $C$ be a curve over $k$ that has a $k$-rational point $p_0 \in \CC(k)$. Assume that $h^0(p_1+\cdots+p_d)>1$ for some $k$-rational points $p_1,\ldots,p_d$. Then
$$
\bigl([p_1]-[p_0]\bigr) * \cdots * \bigl([p_d]-[p_0]\bigr) = 0
\qquad\text{in $\CH_0(C^{[d]})_{\Q}$.}
$$
\end{cor}

\begin{rem} If $k$ is algebraically closed then the vanishing statements in the previous two corollaries hold integrally by Rojtman's theorem; see \cite{Roitman}, \cite{Milne}.
 \end{rem}

It is also instructive to rewrite some of the identities \eqref{GK-id1} in terms of the classes
$\De_{n,*}(C)$. (Cf.\ Prop.~3.4 in~\cite{ColGeem} for a similar result modulo algebraic equivalence.)

\begin{cor}\label{GK-cor4} 
Let $C$ be a curve over $k$. Assume that $h^0(D)>1$ for some $k$-rational divisor~$D$ of degree~$d$ on~$C$. Then one has
$$
\sum_{j=1}^{d+1}\, (-1)^j e_{d+1-j}(D) * \frac{\De_{j,*}(C)}{j}=0
\qquad\text{in $\CH_1(C^{[d+1]})_{\Q}$.}
$$
\end{cor}

\begin{proof} Since we use rational coefficients, we can pass to a finite extension of $k$ and assume that $C$ has a $k$-rational point~$p_0$. Now take $s=1$ and $N=d+1$ in~\eqref{GK-id1}. We get
$$
\sum_{i=0}^{d-1}\, (-1)^i\, \ov{e}_i(D) * \frac{\Gamma_{d+1-i}(C)}{d+1-i} = 0\, .
$$
As $\ov{e}_d(D) = 0$ by Corollary \ref{GK-cor2}, we can extend the range of the index $i$ to $\{0,\ldots,d\}$. Substituting the definition of the classes
$\Gamma_n(C)$ in terms of classes $\Delta_{j,*}(C)$, see Def.~\ref{Gamman(a)-def}, we get
\begin{align*}
0 &= 
\sum_{i=0}^{d} \sum_{k=0}^{d+1-i}\, (-1)^{i+k}\, \ov{e}_i(D)\, \binom{d+1-i}{k}\, \frac{1}{d+1-i}\, t^k\, \De_{d+1-i-k,*}(\CC) \\
&= \sum_{i=0}^{d} \sum_{k=0}^{d-i}\, (-1)^{i+k}\, \ov{e}_i(D)\, \binom{d-i}{k}\, t^k\, \frac{\De_{d+1-i-k,*}(C)}{d+1-i-k}\, .
\end{align*}
Now we use the identity
$$
\sum_{i=0}^M\, \binom{d-i}{M-i}\, \ov{e}_i(D) t^{M-i} = e_M(D)
$$
that can be easily checked using \eqref{eibardef}. We then get
$$
0 = \sum_{M=0}^d\, (-1)^M \, e_M(D) * \frac{\Delta_{d+1-M,*}(\CC)}{d+1-M}\, ,
$$
which, setting $j=d+1-M$, gives what we wanted to prove.
\end{proof}

\end{document}